\documentclass[11pt, reqno] {amsart}

\usepackage{amssymb}
\usepackage{hyperref} 
\usepackage{mathrsfs,bbold}
\usepackage{amsmath}
\usepackage{amsthm}
\usepackage{enumerate}
\usepackage{dsfont}
\usepackage{color}
\usepackage[a4paper]{geometry}
\usepackage[utf8]{inputenc}
\usepackage{stmaryrd}
\usepackage{titlesec}
\usepackage{mathabx}  
\usepackage{appendix}
\usepackage{todonotes, fancyhdr}
\usepackage{chngcntr}
\usepackage{cancel}

\usepackage{setspace}
\renewcommand{\baselinestretch}{0.99}

\allowdisplaybreaks[1]

\usepackage[all]{xy}
\numberwithin{subsection}{section}
\numberwithin{subsubsection}{subsection}
\numberwithin{equation}{section} 

\renewcommand{\labelenumi}{\textsf{(\theenumi)}}

\newenvironment{Dem}[1][\unskip]{%
    \begin{list}{\hspace{0.5cm}{\sf \textbf{Proof #1 --}}}{%
        \setlength{\topsep}{0pt}%
        \setlength{\leftmargin}{0pt}%
        \setlength{\rightmargin}{0pt}%
        \setlength{\listparindent}{0pt}%
        \setlength{\itemindent}{0pt}%
        \setlength{\parsep}{0pt}%
        \addtolength{\leftmargin}{20pt}%
        \addtolength{\rightmargin}{0pt}%
    } \item }{\hfill $\rhd$\end{list}\smallskip}

\renewcommand\thesection       {\arabic{section}}
\renewcommand\thesubsection    {\thesection{\boldmath $.$}\arabic{subsection}}
\renewcommand\thesubsubsection    {\thesection{\boldmath $.$}\arabic{subsection}{\boldmath $.$}\arabic{subsubsection}}

\titleformat{\section}[block]
{\filcenter\normalfont\sffamily\bfseries\Large}
{{\hspace{-0.7cm}}\thesection \hspace{0.2em} --\vspace{0.3cm}}{0.5em}{}

\titleformat{\subsection}[block]
{\filcenter\normalfont\sffamily\bfseries\large}  						  
{\hspace{-0.7cm}\thesubsection \hspace{0.5em} \vspace{0.3cm}}{.5em}{}  
\titlespacing{\subsection}{-0pc}{1.5ex plus .1ex minus .2ex}{0pc}

\titleformat{\subsubsection}[block]
{\filcenter\normalfont\sffamily\bfseries}					  
{\hspace{-0.7cm}\thesubsubsection \hspace{0.5em} \vspace{0.3cm}}{.5em}{}  
\titlespacing{\subsection}{-0pc}{1.5ex plus .1ex minus .2ex}{0pc}

\newtheoremstyle{mystyle}
{3pt}               
{3pt}               
{\it }                      
{}                      
{\sffamily\bfseries}             
{}                      
{0.5em}                 
{#1 #2{\Large$.$}  }

\theoremstyle{mystyle}

\newtheorem{thm}{Theorem}
\newtheorem*{thm*}{Theorem}

\newtheorem{cor}[thm]{\hspace{-0.15cm}  {Corollary} }
\newtheorem{lem}[thm]{\hspace{-0.2cm}  {Lemma} }
\newtheorem{prop}[thm]{\hspace{-0.13cm} {Proposition}}
\newtheorem{defn}[thm]{ \hspace{-0.3cm} {Definition}}

\newtheoremstyle{mystyle2}
{3pt}               
{3pt}               
{\it }                      
{}                      
{\sffamily\bfseries}             
{}                      
{0.5em}                 
{\llap{#2 }#1{\hspace{0.2cm}--}}

\theoremstyle{mystyle2}

\newtheorem*{definition*}{Definition}
\newtheorem*{theorem*}{Theorem}
\newtheorem*{Remark*}{Remark}
\newtheorem*{lem*} {Lemma}
\newtheorem*{defn*} {Definition}
\newtheorem*{prop*} {Proposition}
\newtheorem*{cor*} {Corollary}

\newcommand{\iden}{\text{\rm Id}}

\newcommand{\unit}{\mathbf{1}}

\newcommand{\Brac}[1]{[\![#1]\!]}

\newcommand{\trino}[1]{|\!|\!|#1|\!|\!|}

\newcommand{\bsf}{\boldsymbol{f}}
\newcommand{\bsh}{\boldsymbol{h}}

\newcommand{\rap}{\mathrm{rap}}
\newcommand{\slow}{\mathrm{slow}}
\newcommand{\BHZ}{\mathrm{BHZ}}

\newcommand{\bfP}{\mathbf{P}}
\newcommand{\bfQ}{\mathbf{Q}}

\newcommand{\reR}{\textsf{\textbf{R}}}

\newcommand{\sfM}{\mathsf{M}}
\newcommand{\sfP}{\mathsf{P}}
\newcommand{\sfR}{\mathsf{R}}
\newcommand{\sff}{\mathsf{f}}
\newcommand{\sfg}{\mathsf{g}}
\newcommand{\sfPi}{\mathsf{\Pi}}

\newcommand{\bsF}{\boldsymbol{F}}

\newcommand{\bsT}{\boldsymbol{T}}
\newcommand{\bsmcB}{\boldsymbol{\mathcal{B}}}

\newcommand{\scrM}{\mathscr{M}}

\newcommand{\scrT}{\mathscr{T}}

\newcommand{\frkL}{\mathfrak{L}}

\newcommand{\frkt}{\mathfrak{t}}
\newcommand{\frkn}{\mathfrak{n}}
\newcommand{\frko}{\mathfrak{o}}
\newcommand{\frke}{\mathfrak{e}}
\newcommand{\frkf}{\mathfrak{f}}

\newcommand{\pord}{\trianglelefteq}
\newcommand{\spord}{\triangleleft}

\newcommand{\ssk}{\smallskip}

\renewcommand{\epsilon}{\varepsilon}

\newcommand{\bbN}{\mathbb{N}}

\newcommand{\bbR}{\mathbb{R}}

\newcommand{\bbZ}{\mathbb{Z}}

\newcommand{\mcA}{\mathcal{A}}
\newcommand{\mcB}{\mathcal{B}} 
\newcommand{\mcC}{\mathcal{C}} 
\newcommand{\mcD}{\mathcal{D}}

\newcommand{\mcF}{\mathcal{F}}
\newcommand{\mcG}{\mathcal{G}}

\newcommand{\mcK}{\mathcal{K}}
\newcommand{\mcL}{\mathcal{L}}

\newcommand\mcS{\mathcal{S}}

\newcommand{\bfJ}{\mathbf{J}}

\begin{document}

\begin{center}
{\LARGE\sffamily{Paracontrolled calculus and regularity structures II   \vspace{0.5cm}}}
\end{center}

\begin{center}
{\sf I. BAILLEUL} \& {\sf M. HOSHINO}
\end{center}

\vspace{1cm}

\begin{center}
\begin{minipage}{0.8\textwidth}
\renewcommand\baselinestretch{0.7} \scriptsize \textbf{\textsf{\noindent Abstract.}} We prove a general equivalence statement between the notions of models and modelled distributions over a regularity structure, and paracontrolled systems indexed by the regularity structure. This takes in particular the form of a parametrisation of the set of models over a regularity structure by the set of reference functions used in the paracontrolled representation of these objects. A number of consequences are emphasized. The construction of a modelled distribution from a paracontrolled system is explicit, and takes a particularly simple form in the case of the regularity structures introduced by Bruned, Hairer and Zambotti for the study of singular stochastic partial differential equations. 
\end{minipage}
\end{center}

\vspace{0.6cm}

\section{Introduction}
\label{SectionIntro}

The set of singular stochastic partial differential equations (PDEs) is characterized by the appearance in each equation of this class of ill-defined products, typically the product of a distribution with a function that is not sufficiently regular. The parabolic Anderson model equation
$$
(\partial_t-\Delta) u = u\zeta,
$$
on  the two dimensional torus is a typical example of singular PDE. The space white noise $\zeta$ has almost surely parabolic H\"older regularity $\alpha-2$, for any $\alpha<1$, and $u$ cannot be expected to have better regularity than being $\alpha$-H\"older. So the product $u\zeta$ does not make sense, since $\alpha+(\alpha-2)<0$. Two different sets of tools for the study of singular stochastic PDEs have emerged recently, under the form of Hairer's theory of regularity structures \cite{Hai, BHZ, CH16, BCCH18} and paracontrolled calculus \cite{GIP, BBF15, BB16}, after Gubinelli, Imkeller and Perkowski' seminal work. Both of them implement the same mantra: Make sense of the equation in a restricted space of functions/distributions whose elements look like the linear combination of reference random quantities, for which the ill-defined terms that come from the analysis of the product problems can be defined using probabilistic tools. Within the setting of regularity structures, Taylor-like pointwise expansions and jet-like objects are used to make sense of what it means to look like a linear combination of reference quantities
$$
f (\cdot) \sim \sum_\tau f_\tau(z)\big({\sf \Pi}^{\sf g}_z\tau\big)(\cdot), \quad \textrm{near $z$, for all spacetime points $z$}.
$$
In the paracontrolled approach, one uses paraproducts to implement this mantra
$$
f \sim \sum_\tau {\sf P}_{f_\tau}[\tau].
$$
Each term ${\sf P}_ab$ is a function or a distribution. This approach is justified at an intuitive level by the fact that ${\sf P}_{f_\tau}[\tau]$ can be thought of as a modulation of the reference function/distribution $[\tau]$. The two options seem technically very different from one another.

\ssk

While Hairer's theory has now reached the state of a ready-to-use black box for the study of singular stochastic PDEs, like Cauchy-Lipschitz well-posedness theorem for ordinary differential equations, the task of giving a self-contained treatment of renormalisation matters within paracontrolled calculus remains to be done. It happens nonetheless to be possible to compare the two languages, independently of their applications to the study of singular stochastic PDEs. This task was initiated in Gubinelli, Imkeller, Perkowski' seminal work \cite{GIP} and Martin and Perkowski's work \cite{MP18}, and in our previous work \cite{BH}, where we proved that the set of admissible models $\sf M=(g,\Pi)$ over a concrete regularity structure $\scrT=\big((T^+,\Delta^+), (T,\Delta)\big)$ equipped with an abstract integration map is parametrised by a \emph{paracontrolled representation of} $\sf \Pi$ on the set of elements $\tau$ with non-positive homogeneity. (Admissible models play a crucial in the regularity structures approach to the study of singular stochastic PDEs.) Theorem 21 in \cite{BH} says indeed that given \emph{any} family $\big(\Brac{\tau}\in\mcC^{\vert\tau\vert}\big)_{\vert\tau\vert\leq 0}$, with $\tau$ in a linear basis of $T$, there exists a unique admissible model $\sf (g,\Pi)$ on $\mathscr{T}$ such that one has
\begin{equation} \label{EqPCRepresentationPi}
{\sf \Pi} \tau = \sum_{\sigma<\tau} {\sf P}_{{\sf g}(\tau/\sigma)}\Brac{\sigma} + \Brac{\tau},
\end{equation}
for all $\tau\in T$ in the basis, with non-positive homogeneity. (All notations and words are explained below.) This result provides a parametrisation of the \emph{nonlinear} set of admissible models by a \emph{linear} space, providing for instance a natural notion of tangent space to the space of admissible models. The distribution $\Brac{\tau}$ appears in \eqref{EqPCRepresentationPi} as `the' part of ${\sf\Pi}\tau$ of regularity $\vert\tau\vert$ in this decomposition, while the paraproducts ${\sf P}_{{\sf g}(\tau/\sigma)}\Brac{\sigma}$ have regularity $\vert\sigma\vert<\vert\tau\vert$, for $\sigma<\tau$.

\medskip

To understand the practical relevance of this linear parametrization of the space of admissible models on $\mathscr{T}$, assume $\mathscr{T}$ stands for the Bruned, Hairer, and Zambotti's regularity structure \cite{BHZ} associated with a singular stochastic PDE and ${\sf M}^\epsilon = ({\sf g}^\epsilon,{\sf \Pi}^\epsilon)$ stands for the naive interpretation model associated with a smoothened noise in the equation, with regularization parameter $\epsilon$. The BPHZ renormalization procedure for the model involves a real-valued map $k$ acting on a side space $T^-$, which also defines a homogeneity-preserving linear map $\widetilde{k}$ from $T$ into itself. It follows from Theorem 21 in \cite{BH} that the bracket data associated with the renormalised model $^k{{\sf M}^\epsilon}$ is simply given by the $\Brac{\widetilde{k}(\tau)}$, for $\tau$ of negative homogeneity. The convergence of renormalised admissible models has thus a direct counterpart in terms of bracket data. This answers one of the problems mentioned at the end of Tapia and Zambotti's work \cite{TZ18} on the parametrization problem for the set of branched rough paths, in the present general setting.

\smallskip

Here is another illustration of the use of the parametrization result of admissible models proved in \cite{BH} that will be developed in Section \ref{SubsectionPCtoModels} and Section \ref{SubsectionCorollaries}. Consider the elementary setting of branched rough paths; they are admissible models on particular examples of regularity structures. Theorem 21 in \cite{BH} gives a direct proof of Lyons' extension theorem, saying that a branched H\"older $p$-rough path has a unique extension into a branched H\"older $q$-rough path, for any $q>p$. (Recall weak geometric rough paths are branched rough paths.) This result allows to define the signature of a branched rough path. Similarly, let $\mathscr{T}$ be a regularity structure built from integration operators, with elements of arbitrary large positive homogeneity. It follows from Theorem 21 in \cite{BH} that an admissible model defined on the quotient space of $\mathscr{T}$, modulo elements of a given positive homogeneity $\alpha$, has a unique extension into an admissible model over the regularity structure $\mathscr{T}$ quotiented by the elements of homogeneity $\beta$, for any $\beta>\alpha$. This allows to define the signature of an admissible model.

\ssk

Such statements are concerned with admissible models on regularity structures associated with singular stochastic PDEs. We step back in the present work and prove a general result giving a parametrization of the \emph{nonlinear} space of arbitrary models $\sf M=(g,\Pi)$ on any reasonable concrete regularity structure, by a \emph{linear} space, in terms of representations of the maps $\sf g$ and $\sf \Pi$ by paracontrolled systems, similar to identity \eqref{EqPCRepresentationPi}. (The set of models on any given regularity structure is always nonempty, as it contains the element ${\sf M}_0 = ({\sf g}_0, {\sf \Pi}_0)$, with ${\sf g}_0$ the character on $T^+$ that sends any basis element of $T^+$ on $1$, and ${\sf \Pi}_0$ the null map. The nonlinearity of the space of models can be seen from the analytical constraints that they need to satisfy, that involves nonlinear operations on $\sf g$.) Being reasonable means here satisfying assumptions \textbf{\textsf{(A-C)}} from Section \ref{SectionRStoPC} and Section \ref{SectionPCtoRS}. We insist here on the fact that these assumptions are not related to any kind of singular stochastic PDE or any dynamics or structure that could be modelled with such a regularity structure. As we shall see, the regularity structures used for the study of singular stochastic PDEs enjoy these properties, so all our results hold for them.

\ssk

The result takes the following form. Given a concrete regularity structure 
$$
\scrT=\big((T^+,\Delta^+), (T,\Delta)\big),
$$
denote by $\scrM_\rap(\scrT, \bbR^d)$ the space of models on $\bbR^d$ decreasing rapidly at infinity. Once again, all terms will be properly defined below.

\medskip

\begin{thm}\label{mainthm1}
Let $\scrT$ be a concrete regularity structure satisfying assumptions \textbf{\textsf{(A-C)}}. Then one can construct a locally Lipschitz continuous map
\begin{equation} \begin{split}   \label{EqMapModelsToRemainders}
\scrM_\rap(\scrT, \bbR^d) &\to
\prod_{\sigma\in\mcB^+\setminus\mcB_X^+} \mcC_\rap^{|\sigma|}(\bbR^d) \times \prod_{\tau\in\mcB\setminus\mcB_{\underline{X}}} \mcC_\rap^{|\tau|}(\bbR^d)   \\
{\sf (g,\Pi)} &\mapsto \Big(\Brac{\sigma}^{\sf M}, \Brac{\tau}^{\sf g}\,;\, \sigma\in\mcB^+\setminus\mcB_X^+, \tau\in\mcB\setminus\mcB_{\underline{X}} \Big)
\end{split} \end{equation}
by giving paracontrolled representations of $\sf g$ and $\sf \Pi$, for $({\sf g}, {\sf \Pi})\in \scrM_\rap(\scrT, \bbR^d)$. Furthermore, $\scrM_\rap(\scrT, \bbR^d)$ is locally bi-Lipschitz homeomorphic to the product space
\begin{align} \label{EqParametrizationThm1}   
\prod_{\sigma\in\mcG_\circ^+} \mcC_\rap^{|\sigma|}(\bbR^d) \times \prod_{\tau\in\mcB_\bullet,\,|\tau|<0} \mcC_\rap^{|\tau|}(\bbR^d).
\end{align}
\end{thm}

\medskip

The first claim in Theorem \ref{mainthm1} is part of Theorem 1 in \cite{BH}; see formulas \eqref{def:bracket g} and \eqref{def:bracket Pi} below for an explicit description of the map \eqref{EqMapModelsToRemainders}. The sets $\mcB$ and $\mcB^+$ are fixed linear bases of the spaces $T$ and $T^+$, respectively, consisting of homogeneous vectors. The set $\mcB_\bullet$ in \eqref{EqParametrizationThm1}parametrizes part of the basis $\mcB$, while the set $\mcG_\circ^+$ parametrizes part of the basis $\mcB^+$. The letter $\mcG$ stands for `generator'. In the present setting of a general concrete regularity structure, the space $T^+$ is not related to $T$, unlike what happens with the special regularity structures used for the study of singular stochastic PDEs. It is thus not surprising that there is some freedom in the construction/parametrization of the map $\sf g$. The degrees of freedom are parametrized by the set $\mcG^+_\circ$, described in assumption \textbf{\textsf{(C)}}. Assumption \textbf{\textsf{(A)}} is a harmless requirement on how polynomials sit within $T$ and $T^+$. Assumption \textbf{\textsf{(B)}} is a very mild requirement on the splitting map $\Delta : T\rightarrow T\otimes T^+$, and assumption \textbf{\textsf{(C)}} is a structure requirement on $T^+$ and $\Delta^+$ that provides a fundamental induction structure. The three assumptions are met by all concrete regularity structures built for the study of singular stochastic PDEs. 

This type of parametrization is not entirely new as Tapia and Zambotti described in \cite{TZ18} a free transitive action of a product of H\"older spaces on the space of branched rough paths, a particular example of model over a particular regularity structure. This action was not proved to be continuous however. In relation with the renormalization problem of stochastic models, Theorem \ref{mainthm1} describes precisely the freedom that we have to tweak a divergent family of models and turn it into a convergent family of models. The renormalization process needs to give converging bracket data $\Brac{\sigma}^{\sf M}, \Brac{\tau}^{\sf g}$. See a forthcoming work.

\medskip

We single out here two direct consequences of Theorem \ref{mainthm1} about density and extension questions on the space of models.

\medskip

$\bullet$ The set of models with rapid decrease is equipped with a family of norms ${\sf M\mapsto \|M\|}_a$, indexed by positive exponents $a$. Smooth functions are known to be dense in any H\"older space $\mcC^\beta_a(\bbR^d)$, with growth exponent $a$, if one sees the latter as a subset of $\mcC^{\beta-\epsilon}_a(\bbR^d)$, for any positive $\epsilon$. Theorem \ref{mainthm1} provides as a consequence a direct proof of the following density result, proved in Section \ref{SubsectionCorollaries} -- see Singh and Teichmann's work \cite{SinghTeichmann} for a similar result, proved therein from an explicit mollification procedure on models.

\begin{cor} \label{CorDensity}
Given any positive exponent $\epsilon$, the set of smooth models with rapid decrease
is dense in the set of models with finite $\mcC_{\rap}^{|\tau|}$-norms, for the topology induced by the $\mcC_{\rap}^{|\tau|-\epsilon}$-norms.
\end{cor} 

\ssk

$\bullet$ Branched rough paths are models on a finite time interval $[0,T]$, over a particular example of concrete regularity structure of the form $\mathscr{T}^+=\big((T^+,\Delta^+),(T^+,\Delta^+)\big)$, that satisfies assumptions \textbf{\textsf{(A-C)}}. These models are entirely determined by their $\sf g$-maps, and the elements of $T^+$ are planted rooted trees with decorations on the nodes in a finite set $\{1,\dots,\ell\}$; edges are not decorated. Defining a branched rough path above an $\ell$-dimensional control $h=(h^1,\dots,h^\ell)$ means defining a $\sf g$-map over $(T^+,\Delta^+)$ such that $\sf g$ assigns $h^j$ to the tree with only one node with decoration $j$, for all $1\leq j\leq \ell$, and no edge. The first proof that this is possible for any choice of H\"older control $h$ was found by Lyons and Victoir \cite{LyonsVictoir}, for geometric rough paths, using the axiom of choice. This unexpected device stimulated further explorations of this questions, and different proofs not using the axiom of choice were given subsequently \cite{Unterberger, Hai, TZ18, LPT18}. Unterberger constructs in \cite{Unterberger} a rough path above $h$ using paraproduct-like tools. Hairer uses in \cite{Hai} the reconstruction theorem for that purpose, while Liu, Pr\"omel and Teichmann use in \cite{LPT18} a version of the reconstruction theorem for Sobolev models and the notion of Sobolev rough path to extend Lyons-Victoir extension in their setting. Tapia and Zambotti used in \cite{TZ18} an explicit form of the Baker-Campbell-Hausdorff formula to bypass the use of the axiom of choice in the construction of a lift, and gave a parametrization of the set of all branched rough paths above $h$. Theorem \ref{mainthm1} provides a direct access to such an extension result, in so far as the family of trees with only one node with decoration $j$, for $1\leq j\leq\ell$, is a subset of the generator set $\mcG^+_\circ$ in that setting. 

\begin{cor} \label{CorLyonsVictoir}
\textbf{\textsf{(Lyons-Victoir's extension theorem)}} -- Given any $\bbR^\ell$-valued H\"older control $h$ on the time interval $[0,T]$, there exists a branched rough path above $h$.
\end{cor}

\medskip

The above branched rough path is said to be a lift of the control $h$. Like in \cite{TZ18}, Theorem \ref{mainthm1} actually gives a parametrization of the set of all branched rough paths above $h$. One can formulate the same extension problem for the models on the class of concrete regularity structures introduced for the study of singular stochastic PDEs by Bruned, Hairer and Zambotti in \cite{BHZ} -- we talk of BHZ regularity structures. We refer the reader to Section \ref{SubsectionBHZ-RS} for basics on BHZ regularity structures, and simply mention here that as in the case of branched rough paths, the elements of BHZ regularity structures $\mathscr{T}$ are rooted decorated trees. Their roots may have decorations, and we denote by $(\bullet_j)_{1\leq j\leq \ell}$ the family of one node trees with no edges, and decoration $j$.

\begin{cor} \label{CorExtension}
\textbf{\textsf{(Extension result for models on BHZ regularity structures)}} -- Given a multidimensional noise $\zeta=(\zeta_1,\dots,\zeta_\ell)$, with $\zeta_i\in \mcC^{\vert\bullet_i\vert}_\textrm{\emph{rap}}(\bbR^m)$, for all $1\leq i\leq\ell$, there exists a model $\sf M=(g,\Pi)$ on the BHZ regularity structure $\mathscr{T}$ such that ${\sf \Pi}(\bullet_j)=\zeta_j$, for all $1\leq j\leq\ell$.
\end{cor}

\medskip

The above model is said to be a lift of the noise $\zeta$. The parametrization of admissible models proved in \cite{BH} shows that one can further impose to the extension that it is an admissible model. (Recall the notion of admissibility is related to a peculiar feature of the regularity structures used for the study of singular stochastic PDEs.) As for the case of branched rough paths, Theorem \ref{mainthm1} actually gives a description of the set of all models above the $\ell$-dimensional noise $\zeta$. Corollary \ref{CorLyonsVictoir} and Corollary \ref{CorExtension} are proved in Section \ref{SubsectionCorollaries}. Corollary \ref{CorLyonsVictoir} and Corollary \ref{CorExtension} are actually previously unnoticed consequences of Theorem 21 in \cite{BH}. The general extension result stated in Corollary \ref{CorExtension} is outside the scope of Theorem 21 in \cite{BH}.

Enough for the consequences of Theorem \ref{mainthm1}; we now turn to the problem of the parametrization of the space of modelled distributions associated with a given model.

\medskip

Given a model $\sf M=(g,\Pi)$ on a concrete regularity structure, natural regularity spaces are given by the H\"older-type spaces $\mcD^\gamma(T,\sf g)$, with generic element
$$
\bsf = \sum_{\tau\in\mcB,\, \vert\tau\vert < \gamma} f_\tau \tau.
$$
For $\sfM=(\sfg,\sfPi)\in\scrM_\rap(\scrT, \bbR^d)$, there is an associated notion of rapidly decreasing space of modelled distributions taking values in the vector space $T$, with regularity exponent $\gamma$, denoted by $\mcD_\rap^\gamma(T, \sfg)$. The parametrization of $\mcD_\rap^\gamma(T,\sf g)$ by data in paracontrolled representations of elements of that space requires in general a {\it structure condition} on these data reminiscent of a similar condition introduced by Martin and Perkowski in \cite{MP18}; it is stated in Theorem \ref{thm:PD to MD}. This non-trivial structure condition has a clear meaning in terms of an extension problem for the map $\sf g$ from the Hopf algebra $T^+$ to a larger Hopf algebra; an interesting technical point on its own. The structure condition happens nonetheless to take a very simple form for special concrete regularity structures satisfying assumption \textbf{\textsf{(D)}}.

\medskip

\begin{thm} \label{mainthm2}
Let a concrete regularity structure $\scrT$ satisfy assumptions \textbf{\textsf{(A-D)}}. Pick $\gamma\in\bbR\setminus\{0\}$ such that $\gamma-|\tau|\notin\bbN$ for any basis element $\tau$ of $T$ with $|\tau|<\gamma$, and $\sfM=(\sfg,\sfPi)\in\scrM_\rap(\scrT, \bbR^d)$. Then one can construct a locally Lipschitz continuous map
$$
\mcD_\rap^\gamma(T,\sfg) \to \prod_{\tau\in\mcB,\, |\tau|<\gamma} \mcC_\rap^{\gamma-|\tau|}(\bbR^d)
$$
by giving a paracontrolled representation of elements in $\mcD_\rap^\gamma(T,\sfg)$. Furthermore, $\mcD_\rap^\gamma(T, \sfg)$ is locally bi-Lipschitz homeomorphic to the product space
$$
\prod_{\tau\in\mcB_\bullet,\, |\tau|<\gamma} \mcC_\rap^{\gamma-|\tau|}(\bbR^d).
$$
\end{thm}

\medskip

See formula \eqref{eq:para system} for the paracontrolled representation of a modelled distribution in $\mcD_\rap^\gamma(T,\sfg)$. 
Similarly, we can see the further homeomorphism result
$$
\scrM_{\rap}\ltimes\mcD_{\rap}^\gamma \simeq 
\prod_{\sigma\in\mcG_\circ^+} \mcC_\rap^{|\sigma|}(\bbR^d) \times \prod_{\tau\in\mcB_\bullet,\,|\tau|<0} \mcC_\rap^{|\tau|}(\bbR^d)
\times \prod_{\tau\in\mcB_\bullet,\, |\tau|<\gamma} \mcC_\rap^{\gamma-|\tau|}(\bbR^d),
$$
where $\scrM_{\rap}\ltimes\mcD_{\rap}^\gamma$ is the space of all pairs $\big(({\sf g},{\sf\Pi}),\bsf\big)$ of models $({\sf g},{\sf\Pi})\in\scrM_{\rap}(\scrT,\bbR^d)$ and modelled distributions $\bsf\in\mcD_{\rap}^\gamma(T,{\sf g})$. Following Corollary \ref{CorDensity}, say here that given a model $\sf M$ on a concrete regularity structure $\mathscr{T}$, the set of modelled distributions with rapid decrease is equipped with a family of norms $\bsf \mapsto \|\bsf\|_a$, indexed by a positive growth exponent $a$. The following result is obtained as a direct consequence of Theorem \ref{mainthm2} and the density of smooth functions in any H\"older space $\mcC_a^\beta(\bbR^d)$, equipped with the weaker $\mcC_a^{\beta-\epsilon}(\bbR^d)$-topology, for any positive exponent $\epsilon$. As pointed out in Section 2 of Singh and Teichmann's work \cite{SinghTeichmann}, one can use the reconstruction theorem to define a mollification operator on modelled distributions and obtain as a consequence a density statement for the set of smooth modelled distributions. Theorem \ref{mainthm2} shows that any mollification operation on H\"older spaces induces a mollification operation on the space of modelled distributions; this result is independent of the reconstruction theorem. See Section \ref{SubsectionCorollaries} for a proof.  

\begin{cor} \label{CorDensitySmoothModelled}
Let a concrete regularity structure $\scrT$ satisfy assumptions \textbf{\textsf{(A-D)}}. Fix a model on $\mathscr{T}$. Given any exponents $\gamma\in\bbR$ as in Theorem \ref{mainthm2} and $\epsilon>0$, the set of smooth elements $\big(({\sf g},{\sf\Pi}),\bsf\big)$ in $\scrM_{\rap}\ltimes\mcD_{\rap}^\gamma$ is dense in the same space but with the topology induced by the $\mcC_{\rap}^{|\tau|-\epsilon}$-norms and the $\mcD_{\rap}^{\gamma-\epsilon}$-norm.
\end{cor}

\medskip

Unlike the other assumptions, assumption \textbf{\textsf{(D)}} is fundamentally a requirement on a linear basis of $T$, not on the concrete regularity structure itself. It may then happen that one basis of $T$ satisfies it whereas another does not. Satisfying assumption \textbf{\textsf{(D)}} thus means the existence of a linear basis satisfying this assumption. It happens that the class of concrete BHZ regularity structures introduced by Bruned, Hairer and Zambotti in \cite{BHZ} for the study of singular stochastic PDEs all satisfy assumption \textbf{\textsf{(D)}}, despite the fact that their canonical bases do not satisfy it. We refer the reader to Section \ref{SubsectionBHZ-RS} for the notations $\mathfrak{t}\in\mathfrak{L}$ and $|\mathfrak{t}|$.

\medskip

\begin{thm}  \label{mainthm3}
Assume that the set $\{|\mathfrak{t}|\}_{\mathfrak{t}\in\mathfrak{L}}\cup\{1\}$ is rationally independent. Then the \emph{BHZ} concrete regularity structures satisfy assumptions \textbf{\textsf{(A-D)}}.
\end{thm}

\medskip

\noindent \textbf{\textsf{$\bullet$ BHZ regularity structures vs general regularity structures.}} Readers familiar with the use of regularity structures for the study of singular stochastic PDEs may feel unconfortable at the idea of working regularity structures that do not come from a singular stochastic PDE and with models where the maps $\sf g$ and $\sf\Pi$ are unrelated, unlike in the former setting. This freedom is useful, and Hoshino showed for instance in \cite{HoshinoCommutator1, HoshinoCommutator2} how this leads to a clear understanding of a number of fundamental continuity results for iterated correctors introduced in Bailleul \& Bernicot's work \cite{BB16} on high order paracontrolled calculus, from a regularity structures point of view. As a further illustration of the use of this freedom, let us see how Theorem \ref{mainthm1} gives back a proof of the continuity of the product map $(a,b)\in \mcC^\alpha(\bbR^d)\times \mcC^\beta(\bbR^d) \rightarrow ab\in \mcC^\beta(\bbR^d)$, for $\alpha\in(0,1)$, $\beta<0$, and $\alpha+\beta>0$; this is another formulation of Proposition 4.14 in \cite{Hai}.
Indeed, consider the concrete regularity structure $\mathscr{T}=\big((T^+,\Delta^+), (T,\Delta)\big)$ with
$$
T^+ = \textrm{span}({\bf 1}_+, A), \qquad T = \textrm{span}(B,C),
$$
with $\vert{\bf1}_+\vert:=0,\, \vert A\vert:=\alpha,\, \vert B\vert:=\beta,\, \vert C\vert:=\alpha+\beta$, and splitting maps
$$
\Delta^+{\bf 1}_+ = {\bf 1}_+\otimes {\bf 1}_+, \qquad \Delta^+A = A\otimes{\bf 1}_+ + {\bf 1}_+\otimes A, 
$$
and 
$$
\Delta B = B\otimes {\bf 1}_+, \qquad \Delta C = C\otimes{\bf 1}_+ + B\otimes A. 
$$
Theorem \ref{mainthm1} tells us that the model $\sf (g,\Pi)$, with ${\sf g}({\bf 1}_+)=1$, is uniquely characterized by the two inputs
$$
{\sf g}(A) = a\in \mcC^\alpha(\bbR^d), \qquad {\sf \Pi}B = b\in \mcC^\beta(\bbR^d).
$$
The distribution $c:={\sf \Pi}C$ is in particular determined by $a$ and $b$. We see that $\sf \Pi$ provides an extension of the product map $(a,b)\mapsto ab$, by noting that for smooth inputs $a,b$, the identity $({\sf \Pi}_x^{\sf g}C)(x)=0$, implies in that case $c(x)=a(x)b(x)$, for all $x\in\bbR^d$.

\smallskip

As a matter of fact, working with models with unrelated $\sf g$ and $\sf \Pi$ should somehow be easier than working with admissible models, where $\sf g$ and $\sf \Pi$ are entangled with one another so as to satisfy the admissibility condition.

\smallskip

As far as working with general regularity structures rather than just working regularity structures associated with singular stochastic PDEs is concerned, we would like to encourage the reader to think about general regularity structures as mathematical models of rough 'media' within which one still has a calculus. Rough medias have no reason to be associated with any PDE on a general basis.

\medskip

The following additional remarks put further our results in perspective.
\begin{itemize}
	\item In the theory of regularity structures, the solution map of a singular stochastic PDE has the following structure
$$
\scrM_\rap(\scrT,\bbR^d) \to \mcD_\rap^\gamma(T,\sfg) \xrightarrow{\reR} \mcC^{\beta_0}(\bbR^d).
$$
The first arrow associates to a model the solution in $\mcD_\rap^\gamma(T,\sfg)$ of the regularity structure counterpart of the equation; the second arrow involves the model-dependent reconstruction map $\reR$. The composition of these two maps defines a locally Lipschitz map. Theorem \ref{mainthm1} implies that the solution map actually has the structure
\begin{align}\label{Introduction: solution map}
\prod_{\sigma\in\mcG_\circ^+} \mcC_\rap^{|\sigma|}(\bbR^d) \times \prod_{\tau\in\mcB_\bullet,\,|\tau|<0} \mcC_\rap^{|\tau|}(\bbR^d)
\to \mcC^{\beta_0}(\bbR^d).
\end{align}
The map \eqref{Introduction: solution map} is a general form of the solution maps constructed in the previous works \cite{GIP, BB16} on paracontrolled calculus. Since the ansatz on solutions were given by hand in those papers, it was very hard to extend the argument to a whole class of equations. Our results reveal the relation between such handmade ansatz and the sophisticated algebraic structure in Hairer's theory, showing that it is possible to apply paracontrolled calculus to more general equations in an automatic way, like the works \cite{Hai, BHZ, CH16, BCCH18}.   \vspace{0.15cm}
	
	\item The map \eqref {Introduction: solution map} provides interesting insights on parts of the theory of regularity structures. For example, one of the difficult part of the theory is the continuity result for the model-dependent multi-level extension
$$
\mcK^{\sf M}: \mcD^\gamma(T,{\sf g}) \to \mcD^{\gamma+2}(T,{\sf g}),
$$
of the resolution map $\mcL^{-1}$, with the property that $\reR^{\sf M}(\mcK^{\sf M}\bsf)=\mcL^{-1}(\reR^{\sf M}\bsf)$, for any modelled distributions $\bsf\in\mcD^\gamma(T,{\sf g})$ -- its very definition is non-obvious, see \cite[Section 5]{Hai}. From the paracontrolled
point of view, we take profit from the fact that the classical resolution map $\mcL^{-1}$ preserves the paracontrolled structure 
$$
\mcL^{-1}:\sum \sfP_{f_\tau}\Brac{\tau}+\Brac{\bsf}\mapsto \sum \overline{\sfP}_{f_\tau}\big(\mcL^{-1}\Brac{\tau}\big) + \mcL^{-1}\Brac{\bsf},
$$
up to the introduction of the modified paraproduct $\overline{\sfP}_fg := \mcL^{-1}\sfP_f(\mcL g)$ -- see \cite{BB16}. The main results in the present paper can be applied to such a modified paraproduct. The map $\mcK^{\sf M}$ can be obtained directly from Theorem \ref{mainthm2} by giving first a paracontrolled representation of an element of $\mcD^\gamma(T,{\sf g})$, then applying $\mcL^{-1}$, using the modified paraproduct, and finally using Theorem \ref{mainthm2} again to get back an element of $\mcD^{\gamma+2}(T,{\sf g})$. We do not give the details here and leave it to a future work.   \vspace{0.15cm}
	
	\item The local Lipschitz parametrizations of the sets of models and modelled distributions from Theorem \ref{mainthm1} and Theorem \ref{mainthm2} offer the possibility to define dynamics in these spaces by solving ordinary (or controlled/rough) differential equations driven by vector fields on the parametrization spaces. In the setting of pathspace analysis on manifolds, this kind of pathwise dynamics provided a clean understanding of Driver's flow equation on pathspace, in relation with quasi-invariance questions for Wiener measure on pathspace over a compact Riemannian manifold \cite{DriverFlow, LyonsQian, BailleulRegularity}. One may also make sense of classical stochastic PDEs on the space of models or modelled distributions, as in Liu, Pr\"omel and Teichmann's work \cite{LPT18}.
\end{itemize}

\medskip

Notice that we considered a function space whose elements decrease rapidly at infinity mainly for a technical reason. Our assumption is related to localizing a singular PDE. Indeed, we can consider a class of models on a bounded domain vanishing on the boundary, via a diffeomorphism between $\mathbb{R}^d$ and that domain. Instead, the following modifications are also possible.
\begin{itemize}
   \item Theorems \ref{mainthm1} and \ref{mainthm2} above hold also with the spaces $\scrM_{\slow}$ and $\mcD_{\slow}^\gamma$ of slowly growing models and modelled distributions, respectively, with the spaces $\mcC_{\rap}^\alpha$ replaced by $C_{\slow}^\alpha$. Such a modification is important because temporally or spatially stationary models belong to $\scrM_{\slow}$, but not to $\scrM_{\rap}$. More details can be found in Appendix \ref{AppendixSlow}.   \vspace{0.15cm}
   
   \item If the elements in $\mcG_\circ^+$ and $\mcB_\bullet$ all have homogeneities smaller than $1$, then Theorem \ref{mainthm1} and Theorem \ref{mainthm2} above hold for the unweighted spaces $\scrM$ and $\mcD^\gamma$, with the spaces $\mcC_{\rap}^\alpha$ replaced by usual H\"older spaces $\mcC^\alpha$. An important example is the space of branched rough paths. As said above, Tapia and Zambotti proved in \cite{TZ18} an analogue of Theorem \ref{mainthm1} for the space of branched rough paths by a different approach.
\end{itemize}

\medskip

Like in our previous work \cite{BH}, we work here with the usual isotropic H\"older space rather than with anisotropic spaces. All results given here hold true in that more general setting, with identical proofs. The reader will find relevant technical details in the work \cite{MP18} of Martin and Perkowski.

\ssk

Section \ref{SectionFunctionalSetting} is dedicated to describing different functional spaces and operators. Section \ref{SectionRStoPC} is dedicated to giving paracontrolled representations of models and the reconstruction of modelled distributions in terms of data in paracontrolled systems, proving part of Theorem \ref{mainthm1}. The later is proved in Section \ref{SectionPCtoRS}, where the main work consists in providing a parametrization of $\sf g$-maps by paracontrolled representations, Theorem \ref{mainresult:reconst of g}. Theorem \ref{mainthm2} and Theorem \ref{mainthm3} are proved in Section \ref{SubsectionPCtoModelledDistributions} and \ref{SubsectionBHZ-RS}, respectively. Appendix \ref{SectionAppendixConcreteRS} gives back the setting of concrete regularity structures introduced in \cite{BH}, while Appendix \ref{SectionAppendix} gives a number of technical details that are variations on corresponding results from \cite{BH}.

\bigskip

\noindent \textsf{\textbf{Notations}} $\bullet$ \textit{We use exclusively the letters $\alpha, \beta, \gamma$ to denote real numbers that play the role of regularity exponents, and use the letters  $\sigma, \tau, \mu, \nu$ to denote elements of $T$ or $T^+$.}

$\bullet$ \textit{We agree to use the shorthand notation $\frak{s}^{(+)}$ to mean both the statement $\frak{s}$ and the statement $\frak{s}^+$.}

$\bullet$ \textit{We use the pairing notation $\langle \cdot,\cdot\rangle$ for duality between a finite dimensional vector space and its dual space.}

$\bullet$ \textit{We adopt the notations and terminology of the work \cite{BH}, and write in particular ${\sf \Pi}^{\sf g}_x$ and $\widehat{{\sf g}_{yx}}$, for what is denoted by $\Pi_x$ and $\Gamma_{xy}$ in Hairer's terminology.}

\bigskip

\section{Functional setting}
\label{SectionFunctionalSetting}

We describe in this section different function spaces we shall work with and introduce a modified paraproduct. For $x\in\bbR^d$, set
$$
|x|_* := 1+|x|,\quad x\in\bbR^d.
$$
The weight function $|x|_*$ satisfies the inequalities
$$
|x+y|_*\le|x|_*|y|_*, \qquad |x/\lambda|_*\le|x|_*,
$$
for any $\lambda\ge1$.

\medskip

Let $(\rho_i)_{-1\leq i<\infty}$ be a dyadic decomposition of unity on $\bbR^d$, i.e. $\rho_i:\bbR^d\to[0,1]$ is a compactly supported smooth radial function with the following properties.
\begin{itemize}
\item $\text{\rm supp}(\rho_{-1})\subset\{x\in\bbR^d\,;\,|x|<\frac43\}$ and $\text{\rm supp}(\rho_0)\subset\{x\in\bbR^d\,;\,\frac34<|x|<\frac83\}$.
\item $\rho_i(x)=\rho_0(2^{-i}x)$ for any $x\in\bbR^d$ and $i\ge0$.
\item $\sum_{i=-1}^\infty\rho_i(x)=1$ for any $x\in\bbR^d$.
\end{itemize}
We define the \emph{Littlewood-Paley blocks} $(\Delta_i)_{-1\leq i<\infty}$ by $\Delta_if :=: \rho_i(\nabla)f:= \mcF^{-1}(\rho_i\mcF f)$, where $\mcF$ is a Fourier transform on $\bbR^d$ and $\mcF^{-1}$ is its inverse. For $j\geq -1$, set 
$$
S_j := \sum_{i<j-1}\Delta_i.
$$ 
Denote by $Q_i$ and $P_j$ the integral kernels associated with $\Delta_i$ and $S_j$
$$
\Delta_if(x) := \int_{\bbR^d} Q_i(x-y)f(y)dy, \qquad S_jf(x) := \int_{\bbR^d} P_j(x-y)f(y)dy.
$$

\medskip

\begin{itemize}
	\item[-] For any measurable function $f:\bbR^d\to\bbR$, set
$$
\|f\|_{L_a^\infty(\bbR^d)} := \big\||\cdot|_*^a f\big\|_{L^\infty(\bbR^d)},
$$
and define the corresponding space $L_a^\infty(\bbR^d)$ of functions with finite $\|\cdot\|_{L_a^\infty(\bbR^d)}$-norm. Set
$$
L_{\rap}^\infty(\bbR^d) := \bigcap_{a=1}^\infty L_a^\infty(\bbR^d), \qquad L_{\slow}^\infty(\bbR^d) := \bigcup_{a=1}^\infty L_{-a}^\infty(\bbR^d).
$$

	\item[-] For any distribution $\xi\in\mcS'(\bbR^d)$, set
$$
\|\xi\|_{\mcC_a^\alpha(\bbR^d)} := \sup_{j\ge-1}2^{j\alpha}\|\Delta_j\xi\|_{L_a^\infty(\bbR^d)}.
$$
and define the corresponding space $\mcC_a^\alpha(\bbR^d)$ of functions with finite $\|\cdot\|_{\mcC_a^\alpha(\bbR^d)}$-norm. We have $\mcC^\alpha_0(\bbR^d) = \mcC^\alpha(\bbR^d)$, with the usual definition of the H\"older space $\mcC^\alpha(\bbR^d)$ as the Besov space $\mcB^\alpha_{\infty,\infty}(\bbR^d)$ -- see e.g. Bahouri, Chemin and Danchin's book \cite{BCD11}. Set

$$
\mcC_\rap^\alpha(\bbR^d) := \bigcap_{a=1}^\infty \mcC_a^\alpha(\bbR^d), \qquad \mcC_\slow^\alpha(\bbR^d) := \bigcup_{a=1}^\infty \mcC_{-a}^\alpha(\bbR^d).
$$

	\item[-] For any two-parameter function $F:\bbR^d\times\bbR^d\to\bbR$ and $\alpha>0$, set
$$
\trino{F}_{\mcC_{(2),a}^\alpha(\bbR^d\times\bbR^d)} := \sup_{x,y\in\bbR^d}\big(|x|_*^a\wedge|y|_*^a\big)\frac{\big|F(x,y)\big|}{|x-y|^\alpha}.
$$
and define the corresponding space $\mcC_{(2),a}^\alpha(\bbR^d\times\bbR^d)$ of functions with finite $\|\cdot\|_{\mcC_{(2),a}^\alpha(\bbR^d\times\bbR^d)}$-norm. Set also
$$
\mcC_{(2)}^\alpha(\bbR^d\times\bbR^d) := \mcC_{(2),0}^\alpha(\bbR^d\times\bbR^d), \qquad \mcC^\alpha_{(2),\rap}(\bbR^d\times\bbR^d) := \bigcap_{a=1}^\infty \mcC^\alpha_{(2),a}(\bbR^d\times\bbR^d).
$$

	\item[-] For any $\bbR^d$-indexed family of distributions $\Lambda=(\Lambda_x)_{x\in\bbR^d}\subset\mcS'(\bbR^d)$ on $\bbR^d$, and $\alpha\in\bbR$, set
\begin{align*}
\trino{\Lambda}_{D_a^\alpha}:=\sup_{x\in\bbR^d}\sup_{j\ge-1}|x|_*^a2^{j\alpha}\big| \langle\Lambda_x, P_j(x-\cdot)\rangle\big|.
\end{align*}
Set 
$$
D^\alpha:=D_0^\alpha,   \qquad   D_\rap^\alpha := \bigcap_{a=1}^\infty D^\alpha_a.
$$
\end{itemize}
(In Hairer' seminal work \cite{Hai}, models are assumed to satisfy a ($\lambda, \varphi$)-uniform regularity condition 
$$
\big|({\sf \Pi}^{\sfg}_x\tau)(\varphi_x^\lambda)\big|\lesssim\lambda^{|\tau|},
$$ 
locally uniformly in $x$. Requiring $({\sf \Pi}^{\sfg}_x\tau)_{x\in\bbR^d}\in D^{|\tau|}$ is equivalent to the above uniform estimate -- see e.g. Lemma 6.6 of Gubinelli, Imkeller and Perkowski' seminal work \cite{GIP} on paracontrolled distributions.)

\medskip

For any distributions $f,g\in\mcS'(\bbR^d)$, we define the paraproduct
$$
\sfP_fg:=\sum_{j=1}^\infty (S_jf) (\Delta_jg),
$$
and resonant operator
$$
\sfPi(f,g):=\sum_{|i-j|\le1}(\Delta_if)(\Delta_jg).
$$
For any $g\in\mcS'(\bbR^d)$, set
\begin{equation} \label{EqDefnSmoothS}
\textsc{S} g := g - \sfP_1g = (\Delta_{-1} + \Delta_0)g \in C^\infty(\bbR^d).
\end{equation}
(The letter $\textsc{S}$ is chosen for `smooth'.) The following continuity result is an elementary variation on the classical continuity results for the paraproduct and resonant operators. We refer for instance the reader to Lemma 2.1.34 in J. Martin's thesis \cite{MartinThesis} for a reference.

\medskip

\begin{prop}\label{prop:continuity of the paraproduct}
Let $\alpha,\beta\in\bbR$, $a,b\in\mathbb{Z}$.
\begin{itemize}
\item If $\alpha\neq0$, then $\mcC_a^\alpha(\bbR^d)\times \mcC_b^\beta(\bbR^d)\ni(f,g)\mapsto\sfP_fg\in \mcC_{a+b}^{\alpha\wedge0+\beta}(\bbR^d)$, is continuous.
\item If $\alpha+\beta>0$, then $\mcC_a^\alpha(\bbR^d)\times \mcC_b^\beta(\bbR^d)\ni(f,g)\mapsto\sfPi(f,g)\in \mcC_{a+b}^{\alpha+\beta}(\bbR^d)$, is continuous.
\item If $\alpha,\beta\neq0$ and $\alpha+\beta>0$, then $\mcC_a^\alpha(\bbR^d)\times \mcC_b^\beta(\bbR^d)\ni(f,g)\mapsto f\cdot g\in \mcC_{a+b}^{\alpha\wedge\beta}(\bbR^d)$, is continuous.
\end{itemize}
\end{prop}

\medskip

As a consequence of the last item, the product $fg$, of $f\in\mcS(\bbR^d)$ and $g\in \mcC^\alpha(\bbR^d)$, belongs to $\mcC_\rap^\alpha(\bbR^d)$, for any $\alpha\in\bbR$ -- so the space $\mcC_\rap^\alpha(\bbR^d)$ is in particular not empty.

\medskip

We use a modified paraproduct in Section 3.1.3. Note that 
$$
|\nabla|^m f:=\mcF^{-1}\big(|\cdot|^m\mcF f\big),
$$ 
for $m\in\mathbb{Z}$, is well-defined for functions $f\in\mcS(\bbR^d)$ whose Fourier transform have support in an annulus. For $m\in\bbN$ and $\alpha\in\bbR$, the map $|\nabla|^m$ sends continuously $\mcC_\rap^\alpha(\bbR^d)$ into $\mcC_\rap^{\alpha-m}(\bbR^d)$. For $m\in\bbN$, we define the \textsf{\textbf{modified paraproduct}}
$$
\sfP_f^m g := |\nabla|^m\big(\sfP_f|\nabla|^{-m}g\big) = \sum_{j=1}^\infty|\nabla|^m\big(S_jf\cdot|\nabla|^{-m}\Delta_jg\big).
$$
Note that $\sfP^0=\sfP$. The first item of Proposition \ref{prop:continuity of the paraproduct} also holds for the modified paraproduct $\sfP^m$. This modified paraproduct will play a pivotal role in the proof of Lemma \ref{lem:model on T^M}, along the proof of Theorem \ref{mainresult:reconst of g}. The latter provides the construction of a $\sf g$-map from bracket data.

\bigskip

\section{From regularity structures and models to paracontrolled systems}
\label{SectionRStoPC}

This section sets the scene and contains a proof of the first part of Theorem \ref{mainthm1}. We work in the setting of concrete regularity structures, a special case of regularity structures introduced in \cite{BH}. Their definition is recalled in Appendix \ref{SectionAppendixConcreteRS}. As we work a priori with the most general concrete regularity structures, we need to identify a number of conditions that serve our purpose in Section \ref{SubsectionBasicAssump}. Assumption \textbf{\textsf{(A)}} is a harmless assumption on how polynomials sit inside $T$ and $T^+$. Assumption \textbf{\textsf{(B)}} is a very mild requirement on the splitting map $\Delta : T\rightarrow T\otimes T^+$. Both assumptions are met by the regularity structures used in the study of singular PDEs. This is all we need to get a representation of models and reconstructions of modelled distributions by paracontrolled systems. Before embarking on the journey, recall from \cite{BH} that we use the notations
$$
\Delta\sigma = \sum_{\mu\leq\sigma} \mu\otimes\sigma/\mu,  \qquad \Delta^+\tau = \sum_{\nu\leq^+\tau} \nu\otimes\tau/^+\nu
$$
to denote the action of the splitting map $\Delta$ on $T$ and the coproduct $\Delta^+$ on $T^+$ -- see the comments following assumption \textbf{\textsf{(A)}}. The notation $\mu<\sigma$ will mean $\mu\leq \sigma$ and $\mu\neq\sigma$; we shall make a similar use of the expression $\nu<^+\tau$. 

\medskip

We shall introduce along the way three assumptions \textbf{\textsf{(A), (B), (C)}} on general regularity structures. Their meaning is to be understood in the light of what regularity structures are useful for: They encode the algebra at hand in the pointwise description of `irregular' functions. One will for instance read assumption \textbf{\textsf{(A)}} as saying that the classicaly regular part of functions behave as in the classical Taylor calculus. Interpretations of assumptions \textbf{\textsf{(B)}} and \textbf{\textsf{(C)}} are given after their statement.

\bigskip

\subsection{A basic assumption}
\label{SubsectionBasicAssump}

Appendix \ref{SectionAppendixConcreteRS} recalls elementary properties of concrete regularity structures. Let $\scrT=\big((T^+,\Delta^+),(T,\Delta)\big)$ be a concrete regularity structure with $T^+ = \bigoplus_{\alpha\in A^+}T_\alpha^+$ and $T=\bigoplus_{\beta\in A} T_\beta$. Write ${\bf 1}_+$ for the unit of the algebra $T^+$. Recall that we agree to use the shorthand notation $\frak{s}^{(+)}$ to mean both the statement $\frak{s}$ and the statement $\frak{s}^+$. 

\medskip

\noindent \textbf{\textsf{Assumption (A) --}} {\it The spaces $T^+$ and $T$ have linear bases $\mcB^+$ and $\mcB$, respectively, with the following properties.   \vspace{0.1cm}
\begin{enumerate}
\renewcommand{\labelenumi}{\textbf{\textsf{(\arabic{enumi})}}}
   \item $\mcB^+$ is a commutative monoid freely generated by a finite set $\mcB_\circ^+$ and Taylor monomials $X_1,\dots,X_d$. 
Each element $\tau\in\mcB_\circ^+$ has a positive homogeneity. For general elements in $\mcB^+$, homogeneities are defined by $|X_i|=1$, and multiplicativity
$$
|\tau\sigma|=|\tau|+|\sigma|.
$$

   \item The action of $\Delta^+$ on polynomials is characterised by its action on the monomials
\begin{equation} \label{EqPreLeibniz}
\Delta^+X_i = X_i\otimes {\bf 1}_+ + {\bf 1}_+\otimes X_i,
\end{equation}
that holds for all $1\leq i\leq d$. Denote by $\mcB_X^+$ the submonoid generated by $X_1,\dots,X_d$.   \vspace{0.15cm}

   \item There exists a subset $\mcB_\bullet\subset\mcB$, such that $\mcB$ is in bijection with $\bbN^d\times\mcB_\bullet$. An element $(k,\sigma)\in\bbN^d\times\mcB_\bullet$, is denoted by $\underline{X}^k\sigma$, and assigned a homogeneity
$$
|\underline{X}^k\sigma| := |k|+|\sigma|.   \vspace{0.1cm}
$$

   \item If $\mcB_\bullet$ contains an element ${\bf 1}$ with homogeneity $0$, then it is unique and satisfies the identity
$$
\Delta{\bf 1} = {\bf 1}\otimes{\bf 1}_+.
$$
Write $\underline{X}^k$ for $\underline{X}^k{\bf1}$. Set 
$$
\mcB_{\underline{X}} := \{\underline{X}^k\}_{k\in\bbN^d}\subset\mcB.
$$
The coproduct $\Delta$ on $\underline{X}^k$ is characterised by its action on the monomials
$$
\Delta \underline{X}_i = \underline{X}_i\otimes {\bf 1}_+ + {\bf 1}\otimes X_i,
$$
where $\underline{X}_i=\underline{X}^{e_i}$ and $e_i=(\delta_{ij})_{j=1}^d\in\bbN^d$, that holds for all $1\leq i\leq d$, and by requiring multiplicativity on $\mcB_{\underline{X}}$. For general elements, one has the multiplicative formula
$$
\Delta(\underline{X}^k\sigma) = (\Delta \underline{X}^k)\,(\Delta\sigma).
$$
\end{enumerate}}

\medskip

For later use, denote by $\{\tau'\}_{\tau\in\mcB}$ the dual basis of $\mcB$. Following \cite{BH}, for $\sigma,\tau\in\mcB^{(+)}$, write $\sigma\le^{(+)}\tau$, if $\sigma$ appears in the left hand side of the tensor products in the optimal expansion of $\Delta^{(+)}\tau$, so we have the unique representation
$$
\Delta^{(+)}\tau = \underset{\sigma\le^{(+)}\tau}{\sum_{\sigma\in\mcB^{(+)}}} \sigma\otimes(\tau/^{(+)}\sigma),
$$
where $\tau/^{(+)}\sigma\in T^+\setminus\{0\}$. 
The relation $\le$ needs not to be transitive. Using the coassociativity 
$$
\big(\Delta^{(+)}\otimes\text{\rm Id}\big)\Delta^{(+)}\tau = \big(\text{\rm Id}\otimes\Delta^+\big)\Delta^{(+)}\tau,
$$ 
we obtain the chain formula
\begin{align}\label{chain formula}
\Delta^{+}(\tau/^{(+)}\mu)=\sum_{\mu\le^{(+)}\nu\le^{(+)}\tau}(\nu/^{(+)}\mu)\otimes(\tau/^{(+)}\nu).
\end{align}

Write $\sigma<^{(+)}\tau$, if $\sigma\le^{(+)}\tau$ and $\sigma\neq\tau$. 
{\it The notations $\tau/^{(+)}\sigma$ and $\sigma<^{(+)}\tau$ are only used for $\tau$ and $\sigma$ in $\mcB^{(+)}$.} Be careful! The notations $\leq, <$, etc. are basis-dependent -- like the matrix of a linear map. 

\medskip

The following structural assumption simplifies some arguments in this paper.

\medskip

\noindent \textbf{\textsf{Assumption (B)}} \textit{
\begin{enumerate}
\renewcommand{\labelenumi}{\textbf{\textsf{(\arabic{enumi})}}}
\item For each $\tau,\sigma\in\mcB$ with $\sigma<\tau$, either $\tau/\sigma\in\text{\rm span}(\mcB_X^+)$ or $\tau/\sigma\in\text{\rm span}(\mcB^+\setminus\mcB_X^+)$.
\item For any $\tau\in\mcB^+\setminus\mcB_X^+$ and $\sigma\in\mcB_X^+$, $\tau/^+\sigma\in\text{\rm span}(\mcB^+\setminus\mcB_X^+)$.
\item For any $\tau\in\mcB\setminus\mcB_{\underline{X}}$ and $\sigma\in\mcB_{\underline{X}}$, $\tau/\sigma\in\text{\rm span}(\mcB^+\setminus\mcB_X^+)$.
\end{enumerate}
}

\medskip

Assumption \textbf{\textsf{(B)}} is about the distinction between polynomial and non-polynomial elements.
Assumption \textbf{\textsf{(B-1)}} means that, in the expansion of $\Delta\tau$, there is no term of the form $\sigma\otimes (\eta+X^k)$, with $\eta\in\text{\rm span}(\mcB^+\setminus\mcB_X^+)$ and $k\in\bbN^d$. It is used to justify the quantity $\Brac{\mu/\tau}^{\sfg}$ in the formula \eqref{eq:para system} below.
Assumptions \textbf{\textsf{(B-2)}} and \textbf{\textsf{(B-3)}} is needed in the proof of Theorem \ref{PropDefnBracket}.
BHZ regularity structure satisfies assumption \textbf{\textsf{(B)}}, since polynomial and non-polynomial elements are obviously distinguished by the number of their edges. See Section \ref{SubsectionBHZ-RS} for details.

\medskip

A natural way to ensure \textbf{\textsf{(B)}} is to give homogeneities not in $\bbN$ for the non-polynomial elements.
Hence the following is one of the sufficient conditions.

\medskip

\noindent \textbf{\textsf{Assumption (B') --}} \textit{Homogeneities of elements in $\mcB^+\setminus\mcB_X^+$ and $\mcB\setminus\mcB_{\underline{X}}$ are not nonnegative integers.}

\medskip

This is a kind of natural assumption on regularity structures associated with PDEs; Assumption 5.3 in \cite{Hai} is a part of assumption \textbf{\textsf{(B')}} for elements in $\mcB\setminus\mcB_{\underline{X}}$. Under assumption \textbf{\textsf{(B)}} write, for $\tau\in\mcB^+$,
\begin{equation}\label{def of D^ktau} 
\begin{split}
\Delta^+\tau &= \sum_{\sigma\in\mcB^+\setminus\mcB_X^+}\sigma\otimes(\tau/^+\sigma) + \sum_{k}X^k\otimes(\tau/^+X^k)   \\
&=: \sum_{\sigma\in\mcB^+\setminus\mcB_X^+}\sigma\otimes(\tau/^+\sigma) + \sum_{k}\frac{X^k}{k!}\otimes D^k\tau
\end{split} 
\end{equation}
Extend by linearity the map $D^k$ from $T_\alpha^+$ to $T_{\alpha-|k|}^+$, for all $\alpha\in A$. 

\medskip

\begin{lem}\label{lem:recursive rules}
Under assumptions \textbf{\textsf{(A)}} and \textbf{\textsf{(B)}} one has, for all $\sigma,\tau\in T^+$ and all $k,\ell\in\bbN^d$,
\begin{enumerate}
\setlength{\itemsep}{0.07cm}

	\item $D^0\tau = \tau$,
	
	\item if $\tau\in\mcB^+\setminus\mcB_X^+$, then $D^k\tau\in\text{\rm span}(\mcB^+\setminus\mcB_X^+)$,
	moreover, $D^k\tau\neq 0$ only if $|k|<|\tau|$,
		
	\item $D^{k}D^\ell\tau = D^{k+\ell}\tau$,

	\item $D^{k}X^\ell = \mathbf{1}_{k\le \ell}\, \frac{\ell!}{(\ell-k)!} \, X^{\ell-k}$,
	
	\item $D^k(\tau\sigma) = \sum_{k'} \binom{k}{k'} D^{k'}\tau D^{k-k'}\sigma$ -- {\sl Leibniz rule},
\end{enumerate}
\end{lem}

\medskip

\begin{Dem}
Item \textsf{\textit{(a)}} comes from the property \eqref{EqDefnDeltaPlus} satisfied by coproducts in Hopf algebras, recalled in Appendix \ref{SectionAppendixConcreteRS}. The former part of item \textsf{\textit{(b)}} is a consequence of assumption \textbf{\textsf{(B-2)}}.
The latter part comes from the property \eqref{EqDefnDeltaPlus}. Since $|\tau|>0$ by the definition of concrete regularity structures, the term $\frac{X^k}{k!}\otimes D^k\tau$ appears in the expansion of $\Delta^+\tau$ only if $|k|<|\tau|$.
Item \textsf{\textit{(c)}} is a consequence of the coassociativity property
$$
(\Delta^+\otimes\textrm{Id})\Delta^+ = (\textrm{Id}\otimes\Delta^+)\Delta^+
$$
of the coproduct $\Delta^+$. Expanding both sides at $\tau\in\mcB^+$, we have
\begin{align*}
&\sum_{\sigma,\eta\in\mcB^+\setminus\mcB_X^+}\eta\otimes(\sigma/^+\eta)\otimes(\tau/^+\sigma)
+\sum_{\sigma\in\mcB^+\setminus\mcB_X^+,\,k\in\bbN^d}\frac{X^k}{k!}\otimes D^k\sigma\otimes(\tau/^+\sigma)
+\sum_{k,\ell\in\bbN^d}\frac{X^k}{k!}\otimes\frac{X^\ell}{\ell!}\otimes D^{k+\ell}\tau\\
&=\sum_{\sigma\in\mcB^+\setminus\mcB_X^+}\sigma\otimes\Delta^+(\tau/^+\sigma)
+\sum_{k\in\bbN^d}\frac{X^k}{k!}\otimes \Delta^+D^k\tau.
\end{align*}
It gives indeed the identity
\begin{align}
\label{defi:D^ktau}
\Delta^+D^k\tau &= \sum_{\sigma\leq^+\tau,\, \sigma\notin\mcB_X^+}D^k\sigma\otimes(\tau/^+\sigma) + \sum_{\ell\in\bbN^d}\frac{X^{\ell}}{\ell!}\otimes D^{k+\ell}\tau,
\end{align}
this means \textsf{\textit{(c)}}. Item \textsf{\textit{(d)}} is a direct consequence of the Leibniz formula for the polynomials, which follows from identity \eqref{EqPreLeibniz} giving the action of $\Delta^+$ on $X_i$ and the multiplicativity property of $\Delta^+$. Item \textsf{\textit{(e)}} is again a consequence of the multiplicativity property of $\Delta^+$.
\end{Dem}

\bigskip

\subsection{From models to paracontrolled systems}
\label{SubsectionFromModelsToPCSystems}

We recall in this section some of the results proved in \cite{BH}, stated here in the slightly more general setting of the present work. 
The proofs of these extensions are given in Appendix \ref{SectionAppendix}. These results are proved in Sections 2 and 3 in \cite{BH} without any extra assumptions about `bounded polynomials' and interaction between $T^+$ and $T$. Hence the proofs are completely parallel to the proofs in \cite{BH}, except for the use here of the modified paraproduct and the weight.

\ssk

Given Fr\'echet spaces $E$ and $F$, denote by $L(E,F)$ the space of continuous linear maps from $E$ into $F$. Recall $G^+$ stands for the set of characters of the Hopf algebra $T^+$. Given maps
$$
\sfg:\bbR^d\to G^+,\quad \sfPi\in L\big(T,\mcS'(\bbR^d)\big),
$$
and $x,y\in\bbR^d$, set 
$$
\sfg_{yx} := (\sfg_y\otimes\sfg_x^{-1})\Delta^+\in G^+,
$$
and 
$$
{\sf \Pi}^{\sf g}_x:=(\sfPi\otimes\sfg_x^{-1})\Delta\in L\big(T,\mcS'(\bbR^d)\big).
$$
Set 
$$
\beta_0:=\min A,
$$ 
where $A$ is a homogeneity set of $T=\bigoplus_{\beta\in A}T_\beta$.

\medskip

\begin{defn}  \label{defi:weighted model}
Let a concrete regularity structure $\scrT$ satisfying assumption \textbf{\textsf{(A)}} be given. We denote by 
$$
\scrM_\rap(\scrT,\bbR^d),
$$ 
the set of pairs of maps
$$
\sfg : \bbR^d\to G^+,\quad \sfPi\in L\big(T,\mcS'(\bbR^d)\big),
$$
such that 
\begin{enumerate}
	\item one has $\sfg_x(X^k)=x^k$, for all $x\in\bbR^d, k\in\bbN^d$;   \vspace{0.1cm}

	\item for any $\tau\in\mcB_\circ^+$, the function $x\mapsto\sfg_x(\tau)$ belongs to $L_\rap^\infty(\bbR^d)$, and the function 
$$
(x,y)\mapsto\sfg_{yx}(\tau),
$$ 
belongs to $\mcC_{(2),\rap}^{|\tau|}(\bbR^d\times\bbR^d)$;   \vspace{0.1cm}

	\item one has $(\sfPi \underline{X}^k\sigma)(x)=x^k(\sfPi\sigma)(x)$ and $(\sfPi{\bf1})(x) = 1$;   \vspace{0.1cm}

	\item for any $\tau\in\mcB_\bullet\setminus\{{\bf1}\}$, one has $\sfPi\tau\in \mcC^{\beta_0}_\rap(\bbR^d)$, and the $\bbR^d$-indexed family of distributions $({\sf \Pi}^{\sf g}_x\tau)_{x\in\bbR^d}$ belongs to $D_\rap^{|\tau|}$.   \vspace{0.1cm}
\end{enumerate}
The pair $\sf (g,\Pi)$ is called a \textbf{\textsf{rapidly decreasing model}} on the concrete regularity structure $\scrT$.
\end{defn}

\medskip

We define metrics on the space of rapidly decreasing models on $\scrT$ setting

\begin{align*}
d_a\big(({\sf g},{\sf\Pi}),({\sf g}',{\sf\Pi}')\big)
&:= \sup_{\tau\in\mcB_\circ^+} 
\bigg(\|({\sf g}_\cdot-{\sf g}_\cdot')(\tau)\|_{L_a^\infty(\bbR^d)} 
+ \big\|({\sf g}_{\cdot\cdot}-{\sf g}_{\cdot\cdot}')(\tau)\big\|_{\mcC^{\vert\tau\vert}_{(2),a}(\bbR^d\times\bbR^d)}\bigg)\\
&\quad+\sup_{\sigma\in\mcB_\bullet}
\bigg(\|({\sf\Pi}-{\sf\Pi}')\sigma\|_{\mcC^{\beta_0}_a(\bbR^d)} 
+ \big\|\big(({\sf\Pi}_\cdot^{\sf g}-({\sf\Pi}_\cdot')^{{\sf g}'})\sigma\big)(\cdot) \big\|_{D^{\vert\sigma\vert}_{a}}\bigg).
\end{align*}
With a slight abuse of notations, we write
$$
\sfg_x(\tau)\in L_\rap^\infty(\bbR^d),\quad \sfg_{yx}(\tau)\in \mcC_{(2),\rap}^{|\tau|}(\bbR^d\times\bbR^d).
$$

The definition of a model depends on the choice of subspaces $\text{\rm span}(\mcB^+\setminus\mcB_X^+)$ and $\text{\rm span}(\mcB\setminus\mcB_{\underline{X}})$, but not on the choice of their bases.
Indeed, since
$$
\sfg_x(X^k)\in L_\slow^\infty(\bbR^d),\quad \sfg_{yx}(X^k)\in\mcC_{(2)}^{|k|}(\bbR^d\times\bbR^d)
$$
and since $L_\slow^\infty(\bbR^d)\cdot L_\rap^\infty(\bbR^d)\subset L_\rap^\infty(\bbR^d)$ and $\mcC_{(2)}^\alpha(\bbR^d\times\bbR^d)\cdot\mcC_{(2),\rap}^\beta(\bbR^d\times\bbR^d)\subset\mcC_{(2),\rap}^{\alpha+\beta}(\bbR^d\times\bbR^d)$, for all non-negative $\alpha,\beta$, condition \textsf{\textit{(b)}} holds for any $\tau\in\mcB^+\setminus\mcB_X^+$. Recall that the set $\mcB_\bullet$ in item \textit{\textsf{(d)}} of Definition \ref{defi:weighted model} stands for the index parametrizing the non-polynomial part of the basis of $T$. It is not so obvious to see whether condition \textsf{\textit{(d)}} holds or not for any $\tau\in\mcB\setminus\mcB_{\underline{X}}$; however, the following lemma holds.

\medskip

\begin{lem}
Assume \textbf{\textsf{(A)}} and \textbf{\textsf{(B)}}.
Under the condition \textsf{\textit{(c)}}, the estimate $({\sf\Pi}_x^{\sfg}\tau)_{x\in\bbR^d}\in D_{\rap}^{|\tau|}$ holds for any $\tau\in\mcB\setminus\mcB_{\underline{X}}$.
\end{lem}

\medskip

\begin{Dem}
We prove the estimate for $\underline{X}^k\tau$, with $\tau\in\mcB_\bullet$, $k\in\bbN^d\setminus\{0\}$. Because of the multiplicative property in assumption \textbf{\textsf{(A-4)}} and item \textsf{\textit{(c)} in Definition \ref{defi:weighted model}}, we have
$$
{\sf\Pi}_x^{\sfg}\underline{X}^k\tau=(\cdot-x)^k{\sf\Pi}_x^{\sfg}\tau.
$$
Recalling the notations at the beginning of Section \ref{SectionFunctionalSetting}, we have
\begin{align*}
S_j\big({\sf\Pi}_x^{\sfg}\underline{X}^k\tau\big)(x)
= \int P_j(x-y)(y-x)^k({\sf\Pi}_x^{\sfg}\tau)(y)\,dy
= \int P_j^{k}(x-y)({\sf\Pi}_x^{\sfg}\tau)(y)\,dy,
\end{align*}
where $P_j^k(x):=(-x)^kP_j(x)$. Hence $P_{j+1}*P_j^k=P_j^k$ by the property of support of Fourier transform, and we have
\begin{align*}
S_j\big({\sf\Pi}_x^{\sfg}\underline{X}^k\tau\big)(x) = P_{j+1}*P_j^k*({\sf\Pi}_x^{\sfg}\tau)(x) = \int P_j^k(x-y)S_{j+1}({\sf\Pi}_x^{\sfg}\tau)(y)\,dy.
\end{align*}
By induction, we assume the required estimate for any $\sigma\le\tau$ (hence $\sigma=\tau$ or $\sigma=\underline{X}^k\eta$ with some $k\neq0$ and $\eta\in\mcB_\bullet$ such that $|\eta|<|\tau|$).
Since either of $\sigma$ or $\tau/\sigma$ is non-polynomial by assumption \textbf{\textsf{(B)}}, for any $a\in\bbN$,
\begin{align*}
\big|S_{j+1}({\sf\Pi}_x^{\sfg}\tau)(y)\big| &= \big|S_{j+1}({\sf\Pi}_y^{\sfg}\widehat{\sfg}_{yx}\tau)(y) \big|
\le\sum_{\sigma\le\tau} \big|\sfg_{yx}(\tau/\sigma)\big| \big|S_{j+1}({\sf\Pi}_y^{\sfg}\sigma)(y)\big|   \\
&\lesssim \big(|x|_*^{-a}+|y|_*^{-a}\big) \sum_{\sigma\le\tau}|y-x|^{|\tau|-|\sigma|}\,2^{-j|\sigma|}.
\end{align*}
By using the scaling property $P_j^k(x)=2^{j(d-|k|)}P_j^k(2^jx)$ and by a similar argument to Lemma \ref{lem:weighted scaling} in Appendix \ref{SectionAppendix}, we can conclude that
\begin{align*}
\big|S_j({\sf\Pi}_x^{\sfg}\underline{X}^k\tau)(x)\big| \lesssim |x|_*^{-a}\,2^{-j(|\tau|+|k|)};
\end{align*}
hence $\big({\sf\Pi}_x^{\sfg}\underline{X}^k\tau\big)_{x\in\bbR^d}\in D_{\rap}^{|\underline{X}^k\tau|}$.
\end{Dem}

\medskip

The next statement is a variation on Proposition 12 of \cite{BH}, where we use now the usual polynomials and polynomial weights, and the modified paraproducts ${\sf P}^m$ instead of the bounded polynomials, no weights and the usual paraproduct $\sf P$. Its proof is given in Appendix \ref{SectionAppendix}.

\medskip

\begin{thm}\label{PropDefnBracket}
Let $\scrT$ stand for a regularity structure satisfying assumptions \textbf{\textsf{(A)}} and \textbf{\textsf{(B)}}. Pick $m\in\bbN$. For any model $\sfM=(\sfg,\sfPi)\in\scrM_\rap(\scrT,\bbR^d)$, there exists a family 
$$
\bigg(\Big(\Brac{\tau}^{m,\sfg}\in \mcC^{\vert\tau\vert}_\rap(\bbR^d)\Big)_{\tau\in\mcB^+\setminus\mcB_X^+}, \, \Big(\Brac{\sigma}^{m,\sfM}\in \mcC^{\vert\sigma\vert}_\rap(\bbR^d)\big)_{\sigma\in\mcB\setminus \mcB_{\underline{X}}}\bigg)
$$ 
such that one has, for any $\tau\in\mcB^+\setminus\mcB_X^+$ and $\sigma\in\mcB\setminus \mcB_{\underline{X}}$, the identities
\begin{align}
\label{def:bracket g}
\sfg(\tau) &= \underset{\nu\in\mcB^+\setminus\mcB_X^+}{\sum_{\unit_+<^+\nu<^+\tau}} \sfP^m_{\sfg(\tau/^+\nu)}\Brac{\nu}^{m,\sfg} + \Brac{\tau}^{m,\sfg},   \\
\label{def:bracket Pi}
\sfPi\sigma &= \underset{\mu\in\mcB\setminus \mcB_{\underline{X}}}{\sum_{\mu<\sigma}} \sfP^m_{\sfg(\sigma/\mu)}\Brac{\mu}^{m,\sfM} + \Brac{\sigma}^{m,\sfM}.
\end{align}
Moreover, the mapping
$$
\sfM \mapsto \bigg(\Big(\Brac{\tau}^{m,\sfg}\in \mcC^{\vert\tau\vert}_\rap(\bbR^d)\Big)_{\tau\in\mcB^+\setminus\mcB_X^+}, \, \Big(\Brac{\sigma}^{m,\sfM}\in \mcC^{\vert\sigma\vert}_\rap(\bbR^d)\Big)_{\sigma\in\mcB\setminus \mcB_{\underline{X}}}\bigg)
$$
is locally Lipschitz continuous.
\end{thm}

\medskip

This version of the statement, with $m\geq 1$, will be used in the proof of Theorem \ref{mainresult:reconst of g} given in Section \ref{SubsectionPCtoModels}. Write $\Brac{\tau}^\sfg$ and $\Brac{\sigma}^\sfM$ instead of $\Brac{\tau}^{m,\sfg}$ and $\Brac{\sigma}^{m,\sfM}$, when $m=0$. 

\medskip

Given a model ${\sf M}\in\scrM_\rap(\scrT,\bbR^d)$ on a regularity structure $\scrT$, and $\gamma\in\bbR$, define the space $\mcD^\gamma_\rap(T,\sf g)$ of \textbf{\textsf{rapidly decreasing modelled distributions}} as the set of functions 
$$
{\bsf} : \bbR^d\to\bigoplus_{\beta<\gamma} T_\beta,
$$
such that, for each $\tau\in\mcB$, the function $\langle\tau',{\bsf}(\cdot)\rangle$ belongs to $L^\infty_\rap(\bbR^d)$, and the function
$$
(x,y) \mapsto \big\langle\tau',{\bsf}(y) - \widehat{{\sf g}_{yx}}\bsf(x)\big\rangle
$$
belongs to $\mcC^{\gamma-\vert\tau\vert}_{(2),\rap}(\bbR^d\times\bbR^d)$. The reconstruction $\reR\bsf$ of $\bsf\in\mcD_{\rap}^\gamma(T,{\sf g})$ is an element of $\mcS'(\bbR^d)$ satisfying the condition
$$
\Big(\textsf{\textbf{R}}\bsf - {\sf \Pi}^{\sfg}_x\bsf(x)\Big)_{x\in\bbR^d}\in D^\gamma_\rap.
$$
If $\gamma>0$, there exists exactly one reconstruction. If $\gamma<0$, there are infinitely many reconstructions and two reconstructions are equal modulo $\mcC_{\rap}^\gamma(\bbR^d)$. (This is a key point to prove Proposition \ref{modifiedCorUniqueExtension}.) In what follows, we assume $\gamma\neq0$ and denote $\reR\bsf$ by the one defined in Corollary \ref{canonical reconstruction} in Appendix \ref{SectionAppendix}. If $\gamma=0$, existence of the reconstruction is not ensured in general. See Example 5.5 in \cite{CZ}.

The next statement was proved in \cite{BH}, Theorem 14, in the unweighted setting; its extension to the present setting is given in Appendix \ref{SectionAppendix}.

\medskip

\begin{thm}   \label{thm:MD to PD}
Let $\scrT$ be a regularity structure satisfying assumptions \textbf{\textsf{(A)}} and \textbf{\textsf{(B)}}. Let a regularity exponent $\gamma\in\bbR\setminus\{0\}$ and a model $\sfM=(\sfg,\sfPi)\in\scrM_\rap(\scrT,\bbR^d)$ on $\scrT$ be given. For any modelled distribution 
$$
\bsf=\sum_{|\sigma|<\gamma}f_\sigma\sigma\in\mcD_\rap^\gamma(T,\sf g),
$$ 
each coefficient $f_\sigma$ has a paracontrolled representation
\begin{align}\label{eq:para system} 
f_\sigma = \underset{\mu/\sigma\,\in\,\textrm{\emph{span}}(\mcB^+\setminus\mcB_X^+)}{\sum_{\sigma<\mu}} \sfP_{f_\mu}\Brac{\mu/\sigma}^\sfg + \Brac{f_\sigma}^\sfg,
\end{align}
where $\Brac{f_\sigma}^\sfg\in \mcC_\rap^{\gamma-|\sigma|}(\bbR^d)$. (The quantity $\Brac{\mu/\sigma}^{\sfg}$ is defined as a linear extension of the symbols $\Brac{\tau}^{\sfg}$ in Theorem \ref{PropDefnBracket}.) Moreover, there exists a distribution $\Brac{\bsf}^\sfM\in \mcC_\rap^\gamma(\bbR^d)$ such that
\begin{align}\label{eq:para reconst}
\reR\bsf=\sum_{\sigma\in\mcB\setminus\mcB_{\underline{X}}} \sfP_{f_\sigma}\Brac{\sigma}^\sfM+\Brac{\bsf}^\sfM.
\end{align}
The mapping
$$
\Big(\bsf\in\mcD_\rap^\gamma(T,{\sf g})\Big)\mapsto \left( \Big(\Brac{\bsf}^\sfM,\big(\Brac{f_\sigma}^\sfg\big)_{\sigma\in\mcB}\Big) \in \mcC^\gamma_\rap(\bbR^d)\times\prod_{\sigma\in\mcB} \mcC^{\gamma-\vert\sigma\vert}_\rap(\bbR^d) \right)
$$
is locally Lipschitz continuous.
\end{thm}

\medskip

A similar statement with ${\sf P}^m$ used in place of $\sf P$ holds true. We end this section with three useful formulas involving $\sf g$, that will be used in the proof of Theorem \ref{mainresult:reconst of g}. The reader can skip this statement now and come back to it at the moment where it is needed. Recall $D^k\tau=0$, for $\vert k\vert>\vert\tau\vert$. Let $P_X: T^+\rightarrow T^+_X$, stand for the canonical projection map on $T^+_X$, and for $\tau\in\mcB^+$ set
\begin{equation*} \begin{split}
\sff_x(\tau) :=&\; - \big({\sf g}_x\otimes {\sf g}_x^{-1}\big)(P_X\otimes\textrm{Id})(\Delta^+\tau)   \\
 		 =&\; - \sum_{\ell}\frac{x^\ell}{\ell!} \, \sfg_x^{-1}(D^\ell\tau).
\end{split} \end{equation*}
For $\tau\neq {\bf 1}_+$, we also have
\begin{equation*} \begin{split}
\sff_x(\tau) &:= ({\sf g}_x\otimes {\sf g}_x^{-1}) \big((\textrm{Id}-P_X)\otimes\textrm{Id}\big) 
(\Delta^+\tau)  \\
					 &= \sum_{\sigma\leq^+\tau,\, \sigma\notin\mcB_X^+} \sfg_x(\sigma)\,\sfg_x^{-1}(\tau/^+\sigma).
\end{split} \end{equation*}

\medskip

\begin{lem} \label{LemUsefulIdentities}
Let $\scrT$ be a regularity structure satisfying assumption \textbf{\textsf{(A)}} and \textbf{\textsf{(B)}}. For any $\tau\in\mcB^+\setminus\mcB_X^+$ and any $k\in\bbN^d$, we have 
\begin{align}  \label{formula:g_x and f_x}
\sfg_x(D^k\tau) &= \sum_{\sigma\le^+\tau,\,\sigma\notin\mcB_X^+} \sfg_x(\tau/^+\sigma)\sff_x(D^k\sigma).
\end{align}
and
\begin{align}   \label{formula:g_yx and f}
\sfg_{yx}(D^k\tau) &= \sum_{\sigma\le^+\tau,\,\sigma\notin\mcB_X^+} \sfg_{yx}(\tau/^+\sigma)\sff_y(D^k\sigma) - \sum_{\ell}\frac{(y-x)^\ell}{\ell!}\sff_x(D^{k+\ell}\tau), 
\end{align}
and 
\begin{equation}   \label{EqThirdFormula} \begin{split}
\sff_x(D^k\tau) &= {\bf1}_{|k|<|\tau|}\partial_y^k \Big\{({\sfg}_y\otimes\sfg_x^{-1}) \big((\text{\rm Id}-P_X)\otimes\text{\rm Id}\big) \Delta^+\tau\Big\}\big|_{y=x}   \\
		       &= {\bf1}_{|k|<|\tau|}\partial_y^k \bigg\{\sum_{\sigma\le^+\tau,\, \sigma\notin\mcB_X^+} \sfg_y(\sigma)\,\sfg_x^{-1}(\tau/^+\sigma)\bigg\}\Big|_{y=x}.
\end{split} \end{equation}
\end{lem}

\medskip

Note that one cannot interchange in \eqref{EqThirdFormula} the derivative operator with the sum, as a given function ${\sfg}_y(\sigma)$ may not be sufficiently regular to be differentiated $k$ times. Note that formula \eqref{formula:g_yx and f} does not have the classical feature of a Taylor-type expansion formula, which would rather involve an $x$-dependent term in front of $\sfg_{yx}(\tau/^+\sigma)$, in the first term of the right hand side.

\medskip

\begin{Dem}
$\bullet$ Note first that formula \eqref{defi:D^ktau} for $\Delta^+(D^k\tau)$ gives 
\begin{equation} \label{EqUsefulOne}
{\sff}_x(D^k\tau) = ({\sf g}_x\otimes {\sf g}_x^{-1})\big((\textrm{Id}-P_X)\otimes\textrm{Id}\big)\Delta^+D^k\tau = \sum_{\nu\leq^+\tau,\, \nu\notin\mcB_X^+} {\sfg}_x(D^k\nu)\,{\sfg}_x^{-1}(\tau/^+\nu).
\end{equation}
Formula \eqref{formula:g_x and f_x} is an inversion formula for the preceding identity. One obtains the former from the latter by writing

\begin{align*}
\sum_{\sigma\le^+\tau,\,\sigma\notin\mcB_X^+}\sfg_x(\tau/^+\sigma)\sff_x(D^k\sigma) &= \sum_{\nu\le^+\sigma\le^+\tau,\,\sigma,\nu\notin\mcB_X^+} \sfg_x(\tau/^+\sigma)\sfg_x^{-1}(\sigma/^+\nu)\sfg_x(D^k\nu)   \\
&= \sum_{\nu\le^+\sigma\le^+\tau,\,\nu\notin\mcB_X^+} \sfg_x(\tau/^+\sigma)\sfg_x^{-1}(\sigma/^+\nu)\sfg_x(D^k\nu)   \\
&= \sum_{\nu\le^+\tau,\,\nu\notin\mcB_X^+}(\sfg_x^{-1}\otimes\sfg_x)(\tau/^+\nu)\sfg_x(D^k\nu)   \\
&= \sfg_x(D^k\tau).
\end{align*}
(In the second equality, we can remove the condition ``$\sigma\notin\mcB_X^+$" because $\nu\le^+X^k$ implies that $\nu\in\mcB_X^+$. In the last equality, we use the property of the antipode.)

\ssk

$\bullet$ Applying $\sfg_y\otimes\sfg_x^{-1}$ to \eqref{defi:D^ktau}, we have
\begin{align}\label{eq:f and g}
\begin{aligned}
\sfg_{yx}(D^k\tau)
&= \sum_{\mu\le^+\tau,\,\mu\notin\mcB_X^+} \sfg_y(D^k\mu)\,\sfg_x^{-1}(\tau/^+\mu) + \sum_\ell\frac{y^\ell}{\ell!}\,\sfg_x^{-1}(D^{k+\ell}\tau)   \\
&= \sum_{\mu\le^+\nu\le^+\tau,\,\mu\notin\mcB_X^+} \sfg_y(D^k\mu)\,\sfg_y^{-1}(\nu/^+\mu)\,\sfg_{yx}(\tau/^+\nu) - \sum_{\ell'}\frac{(y-x)^{\ell'}}{\ell'!}\,\sff_x(D^{k+\ell'}\tau),
\end{aligned}
\end{align}
where we use the formula \eqref{chain formula} in the expansion of $\sfg_x^{-1}(\tau/^+\mu)$. Identity \eqref{formula:g_yx and f} follows from \eqref{eq:f and g} using \eqref{EqUsefulOne}. Note that $\mu\le^+\nu$ and $\mu\notin\mcB_X^+$ implies that $\nu\notin\mcB_X^+$. 

\ssk

$\bullet$ Formula \eqref{EqThirdFormula} comes from identity \eqref{EqUsefulOne} by rewriting the terms ${\sfg}_x(D^k\nu)$ for $\nu\in\mcB^+\setminus\mcB_X^+$ in an appropriate form. As a preliminary remark, notice that applying $\sfg_{yx}\otimes\sfg_x$ to the defining identity \eqref{def of D^ktau} for the $D^k\nu$, we have
\begin{align*}
\sfg_y(\nu) = \sum_{\sigma\leq^+\nu,\, \sigma\notin\mcB_X^+} \sfg_{yx}(\sigma)\,\sfg_x(\nu/^+\sigma) + \sum_k\sfg_x(D^k\nu)\,\frac{(y-x)^k}{k!}.
\end{align*}
Since one has
$$
\partial_y^k\sfg_{yx}(\sigma)\big|_{y=x} = 0,
$$
for any $x\in\bbR^d$, whenever $|k| < |\sigma|$, one then has
\begin{align}\label{eq:g_x(D^k tau)}
{\sf g}_x(D^k\nu) = \unit_{|k|<|\nu|}\, \partial_y^k\bigg\{\sfg_y(\nu) - \underset{|\sigma| \leq |k|}{\sum_{\sigma<^+\nu,\,\sigma\notin\mcB_X^+}} \sfg_{yx}(\sigma)\,\sfg_x(\nu/^+\sigma)\bigg\}\Big|_{y=x}.
\end{align}
At the same time, one has
\begin{align*}
\sfg_y(\nu) &= \sum_{\substack{\mu\le^+\nu}} (\sfg_x^{-1}*\sfg_x)(\nu/^+\mu)\,\sfg_y(\mu) = \sum_{\substack{\mu\le^+\nu \\ \mu\notin\mcB_X^+}} (\sfg_x^{-1}*\sfg_x)(\nu/^+\mu)\,\sfg_y(\mu) \\
&= \sum_{\substack{\mu\le^+\sigma\le^+\nu \\ \mu,\sigma\notin\mcB_X^+}} \sfg_x(\nu/^+\sigma)\,\sfg_x^{-1}(\sigma/^+\mu)\sfg_y(\mu).
\end{align*}
In the second equality, we use assumption \textbf{\textsf{(B-2)}} to derive $(\sfg_x^{-1}*\sfg_x)(\nu/^+X^k)=0$ for any $k\in\bbN^d$.
In the third equality, we use that $\mu\le^+\sigma$ and $\mu\notin\mcB_X^+$ implies $\sigma\notin\mcB_X^+$. Furthermore, since $\mu\le^+\sigma\notin\mcB_X^+$, $\mu\in\mcB_X^+$, and $|\sigma|\le |k|$ implies $\mu<^+\sigma$ (hence $|\mu|<|k|$), we have
\begin{align*}
\sum_{\substack{\mu\le^+\sigma\le^+\nu \\ \mu,\sigma\notin\mcB_X^+,\,|\sigma|\le|k|}} \sfg_x(\tau/^+\sigma)\,\sfg_x^{-1}(\sigma/^+\mu)\,\sfg_y(\mu) = \sum_{\substack{\sigma\le^+\nu \\ \sigma\notin\mcB_X^+,\,|\sigma|\le|k|}} \sfg_x(\nu/^+\sigma)\,\sfg_{yx}(\sigma)+p_{<|k|}(y),
\end{align*}
where $p_{<|k|}$ is a polynomial of degree less than $|k|$, hence $\partial^k_yp_{<|k|}=0$. We thus obtain from formula \eqref{eq:g_x(D^k tau)}, that
$$
\sfg_x(D^k\nu) = {\bf1}_{|k|<|\nu|}\partial_y^k \bigg\{\sum_{\substack{\mu\le^+\sigma\le^+\nu \\ \mu,\sigma\notin\mcB_X^+,\,|\sigma|>|k|}} \sfg_x(\nu/^+\sigma)\,\sfg_x^{-1}(\sigma/^+\mu)\,\sfg_y(\mu)\bigg\}\Big|_{y=x}.
$$
Inserting this expression in formula \eqref{EqUsefulOne} one gets,
\begin{align*}
\sff_x(D^k\tau)
&= \sum_{\nu\le^+\tau,\,\nu\notin\mcB_X^+} \sfg_x^{-1}(\tau/^+\nu)\,\sfg_x(D^k\nu)   \\
&= \partial_y^k \bigg\{\sum_{\substack{\mu\le^+\sigma\le^+\nu\le^+\tau \\ \mu,\sigma,\nu\notin\mcB_X^+,\,|\sigma|>|k|}} \sfg_x^{-1}(\tau/^+\nu)\,\sfg_x(\nu/^+\sigma)\,\sfg_x^{-1}(\sigma/^+\mu)\,\sfg_y(\mu)\bigg\}\Big|_{y=x}   \\
&= \partial_y^k \bigg\{\sum_{\substack{\mu\le^+\sigma\le^+\nu\le^+\tau \\ \mu\notin\mcB_X^+,\,|\sigma|>|k|}} \sfg_x^{-1}(\tau/^+\nu)\,\sfg_x(\nu/^+\sigma)\,\sfg_x^{-1}(\sigma/^+\mu)\,\sfg_y(\mu)\bigg\}\Big|_{y=x}   \\
&= \partial_y^k \bigg\{\sum_{\substack{\mu\le^+\sigma\le^+\tau \\ \mu,\sigma\notin\mcB_X^+,\,|\sigma|>|k|}} (\sfg_x*\sfg_x^{-1})(\tau/^+\sigma)\,\sfg_x^{-1}(\sigma/^+\mu)\,\sfg_y(\mu)\bigg\}\Big|_{y=x}   \\
&= {\bf1}_{|k|<|\tau|}\partial_y^k \bigg\{\sum_{\mu\le^+\tau,\, \mu\notin\mcB_X^+}\sfg_x^{-1}(\tau/^+\mu)\,\sfg_y(\mu)\bigg\}\Big|_{y=x}.
\end{align*}
In the third line, we can omit the condition $\sigma,\nu\notin\mcB_X^+$ because of $\mu\notin\mcB_X^+$.
In the last line, we use that $(\sfg_x*\sfg_x^{-1})(\tau/^+\sigma)=1$ if and only if $\sigma=\tau$.
\end{Dem}


\section{From paracontrolled systems to models and modelled distributions}
\label{SectionPCtoRS}

We prove the main results of this work in this section. Theorem \ref{mainthm1} gives a parametrization of the space of models by `bracket' data in paracontrolled representations. The main part of the work consists in building a $\sf g$-map from a paracontrolled representation for it on a minimal subset of a linear basis of $T^+$. Assumption \textbf{\textsf{(C)}} below gives a structural assumption on $T^+$ that identifies this minimal set. The proofs of Corollary \ref{CorDensity} on the density of smooth models, Corollary \ref{CorLyonsVictoir} and Corollary \ref{CorExtension} on extension problems, are proved in Section \ref{SubsectionCorollaries}.

\ssk

Theorem \ref{mainthm2} provides a parametrization of the space of modelled distributions of regularity $\gamma$, for a fixed $\gamma\in\bbR$, by a product of H\"older spaces. It is proved in Section \ref{SubsectionPCtoModelledDistributions}. On a technical level, one brings back the proof of Theorem \ref{mainthm2} to an extension problem for the $\sf g$-map from the Hopf algebra $T^+$ to a larger Hopf algebra $T^+_{\bsF}$. This allows to see Theorem \ref{mainthm2} as a corollary of Theorem \ref{mainthm1} under the additional assumption \textbf{\textsf{(D)}}.

\ssk

Unlike the other assumptions, assumption \textbf{\textsf{(D)}} is about a basis $\mcB$ of $T$ rather than about $T$ itself. It is thus possible that a given basis satisfies assumption \textbf{\textsf{(D)}} whereas another does not. This flexibility is at the heart of the proof of Theorem \ref{mainthm3}, dealing with the case of BHZ regularity structures, investigated in Section \ref{SubsectionBHZ-RS}. Those regularity structures introduced by Bruned, Hairer and Zambotti in \cite{BHZ} provide the universal model of regularity structures associated with a subcritical singular stochastic PDE.

\bigskip

\subsection{From paracontrolled systems to models}
\label{SubsectionPCtoModels}

The following claim is the same as Corollary 15 in \cite{BH}, with the modified paraproduct $\sfP^m$ in the role of $\sfP$. Recall from Theorem \ref{PropDefnBracket} the definition of the reference distributions $\Brac{\sigma}^{m,\sf M}$, in the paracontrolled representation of the $\sf \Pi$ operator of a model $\sf M$, using the modified paraproduct ${\sf P}^m$.

\medskip

\begin{prop}   \label{modifiedCorUniqueExtension}
Let $\scrT$ be a regularity structure satisfying assumptions \textbf{\textsf{(A)}} and \textbf{\textsf{(B)}}. Pick $m\in\bbN$, and assume we are given a map ${\sf g}:\bbR^d\to G^+$, such that conditions \textsf{\textit{(a)}} and \textsf{\textit{(b)}} in Definition \ref{defi:weighted model} are satisfied. Then for any family $\Big(\Brac{\tau}\in \mcC_\rap^{|\tau|}(\bbR^d)\Big)_{\tau\in\mcB_\bullet,\, |\tau|<0}$, there exists a unique model $\sfM=(\sfg,\sfPi)\in\scrM_\rap(\scrT,\bbR^d)$ such that 
\begin{align} \label{modifiedCorUniqueExtension:formula}
\sfPi\tau = \sum_{\sigma<\tau} \sfP^m_{\sfg(\tau/\sigma)}\Brac{\sigma}^{m,\sfM} + \Brac{\tau}, \qquad \forall\,\tau\in\mcB_\bullet, \,\vert\tau\vert<0.
\end{align}
The map 
$$
\left(\sfg, \Big(\Brac{\tau}\in \mcC_\rap^{|\tau|}(\bbR^d)\Big)_{\tau\in\mcB_\bullet,\, |\tau|<0}\right) \mapsto {\sf M} \in \scrM_\rap(\scrT,\bbR^d)
$$ 
is continuous.
\end{prop}

\medskip

Note that the distributions $\Brac{\sigma}^{m,\sfM}$ in \eqref{modifiedCorUniqueExtension:formula} are recursively defined by application of Theorem \ref{PropDefnBracket} to the subspace $\bigoplus_{\beta<|\tau|}T_\beta$. If $\sigma\in\mcB_\bullet$ with $\vert\sigma\vert<0$, then $\Brac{\sigma}^{m,\sfM}=\Brac{\sigma}$.

\medskip

\begin{Dem}
$\bullet$ Recall there is no other element than $\unit$ of zero homogeneity in the present setting, and pick a basis vector $\tau\in\mcB_\bullet$ with $|\tau|<0$, and assume that $\sf(g,\Pi)$ is a model on the sector $T_{<|\tau|}$. Set for all $x\in\bbR^d$
$$
\bsh_\tau(x) := \sum_{\sigma<\tau}\, {\sf g}_x(\tau/\sigma)\sigma;
$$
this defines a modelled distribution in $\mcD_\rap^{|\tau|}(T,{\sf g})$. Then the bound $({\sf\Pi}_x^{\sf g}\tau)_{x\in\bbR^d}\in D_{\rap}^{|\tau|}$
is equivalent to that ${\sf \Pi}\tau$ is one of the reconstructions of $\bsh_\tau$. From the version of Theorem \ref{thm:MD to PD} with the modified paraproduct ${\sf P}^m$, the distribution 
$$
\reR\bsh_\tau
=\sum_{\sigma<\tau} {\sf P}^m_{{\sf g}(\tau/\sigma)}\Brac{\sigma}^{m,{\sf M}} + \Brac{\bsh_\tau}^{m,{\sf M}}
$$
is a reconstruction of $\bsh_\tau$. Since 
$$
{\sf \Pi}\tau - \reR\bsh_\tau = \Brac{\tau} - \Brac{\bsh_\tau}^{m,{\sf M}} \in \mcC_\rap^{|\tau|}(\bbR^d),
$$
the distribution ${\sf \Pi}\tau$ appears then as another reconstruction of $\bsh_\tau$. 

\medskip

$\bullet$ If one picks now a basis vector $\mu\in\mcB$, with $|\mu|>0$, then $\bsh_\mu\in\mcD_\rap^{|\mu|}(T,{\sf g})$ has a unique reconstruction, equal to ${\sf \Pi}\mu$, that is characterized by the data 
$$
\Big({\sf\Pi}_x^{\sf g}\sigma, {\sf g}_x(\mu/\sigma) ; x\in\bbR^d, \sigma<\mu\Big),
$$ 
from the defining property of a reconstruction. An elementary induction then shows the existence of a unique extension of $\sf \Pi$ to $T$ that satisfies the property ${\sf \Pi}\tau = \textsf{\textbf{R}}\bsh_\tau$, for every $\tau\in\mcB$ with positive homogeneity.
\end{Dem}

\medskip

The fact that this statement holds not only for the paraproduct $\sf P$ but also for the modified paraproduct ${\sf P}^m$ will play a pivotal role in the proof of Lemma \ref{lem:model on T^M} below. The proof of Proposition \ref{modifiedCorUniqueExtension} makes it clear that the above parametrization of the set of $\sf \Pi$ maps is related to the non-uniqueness of the reconstruction map on the set of modelled distributions of negative regularity exponent. This statement leaves us with the task of giving a parametrization of the set of characters $\sf g$ on $T^+$ by their paracontrolled representation. We need for that purpose to make the following assumptions on the Hopf algebra $(T^+,\Delta^+)$ and the basis $\mcB^+$ of $T^+$. Recall that $D^k:T_\alpha^+\to T_{\alpha-|k|}^+$, is a linear map satisfying the recursive rules from Lemma \ref{lem:recursive rules}. Recall that a pre-order $\pord$ is a reflexive transitive binary relation. Write $\sigma\spord\tau$ if $\sigma\pord\tau$ and $\tau\ntrianglelefteq\sigma$.

\medskip

\noindent \textbf{\textsf{Assumption (C)}}{\it 
\begin{enumerate}
\item[\textbf{\textsf{(1)}}] There exists a finite subset $\mcG_\circ^+$ of $\mcB_\circ^+$ such that $\mcB_\circ^+$ is of the form
$$
\mcB_\circ^+ = \bigsqcup_{\tau\in\mcG_\circ^+} \Big\{D^k\tau\ ;\, k\in\mathbb{N}^d,\ |\tau|-|k|>0\Big\}.
$$   
	\item[\textbf{\textsf{(2)}}] There exists a preorder $\pord$ on the set $\mcB^+$ such that, for each $\tau\in\mcG_\circ^+$, 
	the coproduct $\Delta^+\tau$ is of the form
\begin{align}\label{eq:canonical coproduct}
\Delta^+\tau=\tau\otimes\unit+\sum_{\sigma<^+\tau,\, \sigma\notin\mcB_X^+}\sigma\otimes(\tau/^+\sigma) + \sum_{k}\frac{X^{k}}{k!}\otimes D^k\tau,
\end{align}
with $\sigma\in\mcB^+(\tau^-)$ and $\tau/^+\sigma\in\textrm{\emph{span}}\big(\mcB^+(\tau^-)\big)$,  for each $\sigma$ in the above sum,
where for each $\tau\in\mcB^+$, denote by $\mcB^+(\tau^-)$ the submonoid of $\mcB^+$ generated by
$$
\big\{X_1,\dots,X_d\big\}\cup\bigsqcup_{\sigma\in\mcG_\circ^+,\, \sigma\spord\tau} \Big\{D^k\sigma\ ;\, k\in\mathbb{N}^d,\ |\sigma|-|k|>0\Big\}.
$$

\item[\textbf{\textsf{(3)}}] For any element $\sigma\in\mcB^+\setminus\mcB_X^+$ such that there exists $\tau\in\mcG_\circ^+$ and $\sigma\le^+\tau$, the homogeneity of $\sigma$ is non-integer.
\end{enumerate}   }

\medskip

Note the \emph{disjoint} union in the description of $\mcB_\circ^+$. Assumption \textbf{\textsf{(C-1)}} identifies a set of generators, modulo the action of the $D$ operator. Assumption \textbf{\textsf{(C-2)}} provides a useful induction structure. Assumption \textbf{\textsf{(C-3)}} is a part of assumption \textbf{\textsf{(B')}} and it is used at the end of the proof of Theorem \ref{mainresult:reconst of g}. If one understands the coproduct $\Delta^+$ as giving the elementary pieces of any given element, assumption \textbf{\textsf{(C)}} as a whole provides an inductive description of $\mcB^+$.

\smallskip

As discussed in Section \ref{SubsectionBHZ-RS}, BHZ regularity structure satisfies assumption \textbf{\textsf{(C)}}. Indeed, we can choose $\mcG_\circ^+$ as a set of all conforming trees of the form $I_0^{\mathfrak{t}}(\tau)$, and the operator $D^k$ appears as the form $D^kI_0^{\mathfrak{t}}(\tau)=I_k^{\mathfrak{t}}(\tau)$.
In the BHZ regularity structure, one of the examples of $\pord$ is the binary relation based on the scale of graphs.
Since $\sigma$ and $\tau/^+\sigma$ in \eqref{eq:canonical coproduct} are subtree and quotient tree of $\tau$ respectively, it follows from definition of $\Delta^+$ that $\sigma,\tau/^+\sigma\spord\tau$.
The last assumption is true, if the types $\{\mathfrak{t}\}$ are assigned rationally independent homogeneities $\{|\frak{t}|\}$.
See Theorem \ref{BHZ-ABC} for details. We leave now this special setting and come back to our general setting.

\medskip

\begin{lem}\label{extension of assumption C2}
Denote $T^+(\tau^-):=\text{\rm span}\big(\mcB^+(\tau^-)\big)$.
For any $\tau\in\mcG_\circ^+$ and any $\sigma\in\mcB^+(\tau^-)$, one has
\begin{align}\label{eq:extension of assumption C2}
\Delta^+\sigma\in T^+(\tau^-)\otimes T^+(\tau^-).
\end{align}
\end{lem}

\medskip

\begin{Dem}
By multiplicativity it is sufficient to show the case $\sigma=D^k\eta\in\mcB^+(\tau^-)$ with $\eta\in\mcG_\circ^+$ and $k\in\bbN^d$.
If $k=0$, \eqref{eq:extension of assumption C2} follows because of the transitivity of $\spord$.
By the formula \eqref{defi:D^ktau} for $\Delta^+D^k\eta$ and the property $D^k:T^+(\tau^-)\to T^+(\tau^-)$ that is proved by Lemma \ref{lem:recursive rules}, $D^k\eta$ also satisfies \eqref{eq:extension of assumption C2}.
\end{Dem}

\medskip

Recall from formula \eqref{eq:g_x(D^k tau)} that if we are given characters $({\sf g}_x)_{x\in\bbR^d}$ on $T^+$ as in Definition \ref{defi:weighted model}, then
\begin{equation} \label{EqGxDkTau}
{\sf g}_x(D^k\tau) = \unit_{|k|<|\tau|}\, \partial_y^k\bigg\{\sfg_y(\tau) - \underset{|\sigma| \le |k|}{\sum_{\sigma<^+\tau,\,\sigma\notin\mcB_X^+}} \sfg_{yx}(\sigma)\,\sfg_x(\tau/^+\sigma)\bigg\}\Big|_{y=x}.
\end{equation}
The induction structure from assumption \textbf{\textsf{(C-2)}} restricts the above sum and shows that the family of all ${\sf g}_x(D^k\tau)$ is uniquely determined by the preceding formula. It follows then from assumption \textbf{\textsf{(C-1)}} that the character $\sf g$ on $T^+$ is entirely determined by the datum of the ${\sf g}(\tau)$, for $\tau\in\mcG_\circ^+$. We have in particular, if $\tau\in\mcG_\circ^+$ is minimal (i.e., there is no $\sigma\in\mcG_\circ^+$ such that $\sigma\spord\tau$) then
$$
\sfg_y(\tau) = \sfg_{yx}(\tau) + \sum_{\vert k\vert<\vert\tau\vert}\frac{(y-x)^k}{k!}\,\sfg_x(D^k\tau),
$$
since $\mcB^+(\tau^-)=\mcB_X^+$, so for $\vert k\vert<\vert\tau\vert$, one has
\begin{equation} \label{EqDefnGD}
\sfg_x(D^k\tau) = \partial^k_y\sfg_y(\tau)\big\vert_{y=x},
\end{equation}
and
$$
\sff_x(D^k\tau) = \sfg_x(D^k\tau),
$$
and 
\begin{equation} \label{EqDefnG2D}
\sfg_{yx}(D^k\tau) = \sfg_y(D^k\tau) - \sum_\ell \frac{(y-x)^\ell}{\ell !}\,\sfg_x(D^{k+\ell}\tau).
\end{equation}

\medskip

Recall that, given a concrete regularity structure $\scrT$,
$$
\scrT^+ = \big((T^+,\Delta^+),(T^+,\Delta^+)\big)
$$
is also a concrete regularity structure, and that for a $\sf g$ map as in Definition \ref{defi:weighted model} one defines a model ${\sf M^g = (g,\Pi^g)}$ on $\scrT^+$ setting
$$
\big({\sf \Pi^g}\tau\big)(y) = {\sf g}_{y}(\tau).
$$   

\medskip

\begin{thm}   \label{mainresult:reconst of g}
Let $\scrT$ stand for a concrete regularity structure satisfying assumptions \textbf{\textsf{(A-C)}}. Then, for any family $\Big(\Brac{\tau}\in C^{|\tau|}_\rap(\bbR^d)\Big)_{\tau\in\mcG^+_\circ}$, there exists a unique model $\sfM^{\sf g} = (\sfg,{\sf\Pi}^\sfg)$ on $\scrT^+$ such that 
\begin{equation} \label{EqgTau}
{\sf g}(\tau) = \underset{\sigma\in\mcB^+\setminus \mcB_X^+}{\sum_{\sigma<^+\tau}} \sfP_{\sfg(\tau/^+\sigma)}\Brac{\sigma}^{{\sf M}^{\sfg}} + \Brac{\tau}, \qquad \forall\,\tau\in\mcG^+_\circ.
\end{equation}
The map 
\begin{equation} \label{EqMap}
\Big(\Brac{\tau}\in C^{|\tau|}_\rap(\bbR^d)\Big)_{\tau\in\mcG^+_\circ} \mapsto {\sf M^g}\in\scrM_\rap(\scrT^+,\bbR^d)
\end{equation}
is locally Lipschitz continuous.
\end{thm}

\medskip

Note that one uses the paraproduct $\sf P$ and the brackets $\Brac{\cdot}^{\sf M^g}$ in the statement. The modified paraproduct ${\sf P}^m$ is only used in the proof of Lemma \ref{lem:model on T^M}, where we construct a model on an intermediate regularity structure introduced along the proof. The injectivity of the map \eqref{EqMap} is elementary, so Theorem \ref{mainresult:reconst of g} and Proposition \ref{modifiedCorUniqueExtension}, with Theorem \ref{PropDefnBracket}, prove all together Theorem \ref{mainthm1}.

\medskip

The remaining of this section is dedicated to proving Theorem \ref{mainresult:reconst of g}. The proof is done by induction on the preorder $\pord$.

\medskip

\textbf{\textsf{$\bullet$ Initialisation of the induction.}} 
If $\tau\in\mcG_\circ^+$ is an minimal element, then set 
$$
\sfg(\tau) := \Brac{\tau},
$$
and define $\sfg(D^k\tau)$ and $\sfg_{yx}(D^k\tau)$ by \eqref{EqDefnGD} and \eqref{EqDefnG2D}. It is clear on these formulas that they define elements of the spaces $\mcC^{\vert\tau\vert-\vert k\vert}_\rap(\bbR^d)\subset L^\infty_\rap(\bbR^d)$ and $\mcC^{\vert\tau\vert-\vert k\vert}_{(2),\rap}(\bbR^d\times\bbR^d)$, respectively.

\medskip

\textbf{\textsf{$\bullet$ Induction step.}} Fix $\tau\in\mcG_\circ^+$ and assume that $\sfg$ has been constructed on the submonoid $\mcB^+(\tau^-)$ as a continuous function of the bracket data -- so all the functions $\Brac{\sigma}^{{\sf M}^\sfg}$ and ${\sfg}(\tau/^+\sigma)$ make sense as elements of their natural spaces. Define ${\sf g}(\tau)$ by identity \eqref{EqgTau}, and define ${\sf g}(D^k\tau)$ by \eqref{EqGxDkTau}, for all $k\in\bbN^d$ with $\vert k\vert<\vert\tau\vert$. The induction step consists in proving that $\sfg_x(D^k\tau)\in L_\rap^\infty(\bbR^d)$ and ${\sfg}_{yx}(D^k\tau)\in\mcC^{|\tau|-|k|}_{(2),\rap}(\bbR^d\times\bbR^d)$, as one can use for $\alpha,\beta$ non-negative the inclusions
$$
L_\slow^\infty(\bbR^d)\cdot L_\rap^\infty(\bbR^d)\subset L_\rap^\infty(\bbR^d)
$$ 
and 
$$
\mcC_{(2)}^\alpha(\bbR^d\times\bbR^d)\cdot\mcC_{(2),\rap}^\beta(\bbR^d\times\bbR^d)\subset\mcC_{(2),\rap}^{\alpha+\beta}(\bbR^d\times\bbR^d),
$$ 
to get the regularity properties of $\sfg_x(\mu\,D^k\tau)$ and $\sfg_{yx}(\mu\,D^k\tau)$, for $\mu\in\mcB^+(\tau^-)$.

\medskip

Choose $m\in\bbN$, with $m>\vert\tau\vert$. We introduce a regularity structure $\scrT^m(\tau)$ with Hopf algebra part $T^+(\tau^-)$ and $T$-space defined as follows.
Consider the formal symbols
$$
\sigma^{(m)}
$$
indexed by $\sigma\in\mcB^+\setminus\mcB_X^+$, with homogeneity
$$
\big|\sigma^{(m)}\big| := |\sigma| - m.
$$
Set
\begin{align*}
T^m(\tau) := \textrm{span}
\Big(\big\{\sigma^{(m)}\ ;\,\sigma<^+\tau,\, \sigma\notin\mcB_X^+\big\}
\cup\big\{\tau^{(m)}\big\}\Big),
\end{align*}
so all elements of $T^m(\tau)$ have negative homogeneity. 
Lemma \ref{extension of assumption C2} ensures that we can define a coassociative coproduct 
$$
\delta: T^m(\tau)\to T^m(\tau)\otimes T^+(\tau^-)
$$ 
setting
$$
\delta\big(\sigma^{(m)}\big) := \sum_{\mu\leq^+\sigma,\,\mu\notin\mcB_X^+} \mu^{(m)}\otimes(\sigma/^+\mu)
$$
for each basis element of $T^m(\tau)$.
Lemma \ref{extension of assumption C2} also ensures that
$$
\Delta^+\big(T^+(\tau^-)\big) \subset T^+(\tau^-)\otimes T^+(\tau^-),
$$ 
so
$$
\scrT^m(\tau) := \Big((T^+(\tau^-),\Delta^+), (T^m(\tau),\delta)\Big)
$$ 
is a concrete regularity structure.

\smallskip

We build a model $\sf (g,\Lambda)$ on $\scrT^m(\tau)$, from ${\sf g} : T^+(\tau^-)\rightarrow\bbR$ given by an induction assumption and an operator ${\sf\Lambda}:T^m(\tau)\to\mathcal{S}'(\mathbb{R}^d)$ defined by
$$
{\sf\Lambda}(\sigma^{(m)}) := |\nabla|^m {\sf g}(\sigma),
$$
where $|\nabla|^m$ is the Fourier multiplier operator $|\nabla|^m\zeta = \mcF^{-1}\big(|\cdot|^m\mcF\zeta\big)$. 
The pair $\sf(g,\Lambda)$ turns out to be a model by Lemma \ref{lem:model on T^M} below. Then formula \eqref{EqThirdFormula} giving $\sff_x(D^k\sigma)$ can be interpreted in terms of that model, under the form of identities
$$
\sff_x(D^k\sigma) = {\bf J}^{k,m}\Big({\sf \Lambda}^\sfg_x(\sigma^{(m)})\Big)(x)
$$
for operators ${\bf J}^{k,m}$ on distributions defined below. The identity 
$$
{\sf \Lambda}^\sfg_x = {\sf \Lambda}^\sfg_y \circ\widehat{\sfg_{yx}}^\delta,
$$
where $\widehat{\sfg}^\delta := (\textrm{Id}\otimes \sfg)\delta$, is then used crucially to obtain estimates on $\sff_x(D^k\sigma)$, that eventually give informations on $\sfg_x(D^k\tau)$ and $\sfg_{yx}(D^k\tau)$ via formulas \eqref{formula:g_x and f_x} and \eqref{formula:g_yx and f}.

\medskip

\begin{lem}\label{lem:model on T^M}
The pair $({\sf g},{\sf\Lambda})$ is a rapidly decreasing model on the regularity structure $\scrT^m(\tau)$.
\end{lem}

\medskip

\begin{Dem}
Since we have the identity
$$
{\sf\Lambda}(\sigma^{(m)}) = |\nabla|^m\sfg(\sigma) = \sum_{\mu<\sigma,\,\mu\notin\mcB_X^+}\sfP^m_{\sfg(\sigma/^+\mu)}|\nabla|^m \Brac{\mu}^\sfg + |\nabla|^m\Brac{\sigma}^\sfg,
$$
for all $\sigma\in\mcB^+\setminus\mcB_X^+$ with $\sigma\le^+\tau$, from the intertwining relation defining ${\sf P}^m$ and the induction assumption, the operator $\sf\Lambda$ is the unique model on $\scrT^m(\tau)$ associated by Proposition \ref{modifiedCorUniqueExtension} to the inputs 
$$
\Brac{\sigma^{(m)}}:=|\nabla|^m\Brac{\sigma}^\sfg\in \mcC_\rap^{|\sigma|-m}(\bbR^d),
$$
since all elements of $T^m(\tau)$ have negative homogeneity.
\end{Dem}

\medskip

Note that it follows from identity \eqref{EqThirdFormula} in Lemma \ref{LemUsefulIdentities} that the model $\sfPi$ and the function $\sff(D^k\sigma)$ are related by the identity

\begin{equation}  \label{formula:f_xD^k} 
\begin{split}
\sff_x(D^k\sigma) &= \partial_y^k\bigg\{\sum_{\mu\leq^+\sigma,\,\mu\notin\mcB_X^+}\sfg_x^{-1}(\sigma/^+\mu)\sfg_y(\mu)\bigg\}\Big|_{y=x}   \\
&= \partial_y^k\bigg\{|\nabla_y|^{-m} \sum_{\mu\leq^+\sigma,\,\mu\notin\mcB_X^+}\sfg_x^{-1}(\sigma/^+\mu) {\sf\Lambda}\big(\mu^{(m)}\big)(y)\bigg\}\Big|_{y=x}  \\
&= \partial_y^k\bigg\{|\nabla_y|^{-m} {\sf\Lambda}^{\sfg}_x\big(\sigma^{(m)}\big)(y)\bigg\}\Big|_{y=x}   \\
&=: \sum_j\, \bfJ_j^{k,m}\Big({\sf\Lambda}^{\sfg}_x\big(\sigma^{(m)}\big)\Big)(x),
\end{split} 
\end{equation}
where the operators $\bfJ_j^{k,m}$ are defined by
$$
\bfJ_j^{k,m}(\zeta) := \partial^k|\nabla|^{-m}\Delta_j\zeta,
$$
for an appropriate distribution $\zeta\in\mcS'(\bbR^d)$. If $j\ge0$, since the Fourier transform of $\Delta_j\zeta$ is supported on an annulus, the function $\bfJ_j^{k,m}(\zeta)$ is always well-defined; this is not the case of $\bfJ_{-1}^{k,m}(\zeta)$. However, we only use in this section distributions $\zeta$ of the form $\zeta=|\nabla|^m\xi$ (where such $\xi$ is unique in the class of rapidly decreasing functions), so $\bfJ_{-1}^{k,m}(\zeta)=\partial^k\Delta_{-1}\xi$, in our setting.

\medskip

\begin{lem}  \label{lem:I_j^k}
Under assumptions \textbf{\textsf{(A-C)}}, for any $\sigma\in\mcB^+\setminus\mcB_X^+$ with $\sigma\le^+\tau, k\in\bbN^d$, and $a\in\bbN$, we have
\begin{align*}
\Big|\bfJ_j^{k,m}\big({\sf\Lambda}^{\sfg}_x(\sigma^{(m)})\big)(x)\Big| &\lesssim|x|_*^{-a}\,2^{-j(|\sigma|-|k|)},   \\
\Big|\bfJ_j^{k,m}\big({\sf\Lambda}^{\sfg}_x(\sigma^{(m)})\big)(y)\Big| &\lesssim|y|_*^{-a}\sum_{\mu\le^+\sigma,\,\mu\notin\mcB_X^+}|y-x|^{|\sigma|-|\mu|}2^{-j(|\mu|-|k|)}.
\end{align*}
Consequently, $\sff_x(D^k\sigma)\in L_\rap^\infty$.
\end{lem}

\medskip

\begin{Dem}
For the first estimate, since 
$$
\bfJ_{-1}^{k,m}\big({\sf\Lambda}^{\sfg}_x(\sigma^{(m)})\big)(x) = \sum_{\mu\leq^+\sigma,\,\mu\notin\mcB_X^+}\sfg_x^{-1}(\sigma/^+\mu)\partial_x^k\Delta_{-1}\sfg_x(\mu) \in L_\rap^\infty,
$$ 
by assumption, it is sufficient to consider the case $j\ge0$. By the property of $\rho_j$, there exists a smooth function $\tilde\rho$ supported on an annulus, and such that setting $\tilde\rho_j(\cdot) := \tilde\rho(2^{-j}\cdot)$, one has $\widetilde{\rho}_j\rho_j=\rho_j$. Set 
$$
\widetilde{Q}_j^{k,m} := \partial^k|\nabla|^{-m}(\mcF^{-1}\tilde{\rho}_j),
$$ 
and note the scaling property
$$
\widetilde{Q}_j^{k,m}(\cdot) = 2^{j(d+|k|-m)}\widetilde{Q}_0^{k,m}(2^j\cdot).
$$
We now use the fact that $\sf (g,\Lambda)$ is a model to write
\begin{align*}
\bfJ_j^{k,m}\big({\sf \Lambda}^{\sfg}_x(\sigma^{(m)})\big)(x)
&=\int \widetilde{Q}_j^{k,m}(x-y)\Delta_j\big({\sf\Lambda}^{\sfg}_x(\sigma^{(m)})\big)(y)dy   \\
&=\int \widetilde{Q}_j^{k,m}(x-y)\Delta_j\big({\sf\Lambda}^{\sfg}_y\circ\widehat{\sfg_{yx}}^\delta(\sigma^{(m)})\big)(y)dy   \\
&=\sum_{\mu\le^+\sigma} \int \widetilde{Q}_j^{k,m}(x-y)\,\sfg_{yx}(\sigma/^+\mu)\,\Delta_j\big({\sf\Lambda}^{\sfg}_y(\mu^{(m)})\big)(y)dy.
\end{align*}
Recall that $|x+y|_*\leq |x|_* |y|_*$, for all $x,y\in\bbR^d$. By Lemma \ref{lem:model on T^M}, for any $a\in\bbN$ we have
\begin{align*}
&|x|_*^a\Big|\bfJ_j^{k,m}\big({\sf\Lambda}^{\sfg}_x(\sigma^{(m)})\big)(x)\Big| \\
&\lesssim\sum_{\mu\le^+\sigma}\int |x-y|_*^a \big|\widetilde{Q}_j^{k,m}(x-y)\big| |y-x|^{|\sigma|-|\mu|}\,|y|_*^a\,\Big|\Delta_j\big({\sf\Lambda}^{\sfg}_y(\mu^{(m)})\big)(y)\Big|\,dy   \\
&\lesssim\sum_{\mu\le^+\sigma}2^{-j(|\mu|-m)} \int |z|_*^a\big|\widetilde{Q}_j^{k,m}(z)\big| |z|^{|\sigma|-|\mu|}\,dz   \\
&\lesssim\sum_{\mu\le^+\sigma}2^{-j(|\mu|-m)}2^{j(|k|-m-|\sigma|+|\mu|)} \int |z|_*^a\big|\widetilde{Q}_0^{k,m}(z)\big| |z|^{|\sigma|-|\mu|}\,dz   \\
&\lesssim2^{-j(|\sigma|-|k|)}.
\end{align*}
We get the second estimate from the first using once again the fact that ${\sf (g,\Lambda)}$ is a model, writing
$$
\bfJ_j^{k,m}\big({\sf\Lambda}^{\sfg}_x(\sigma^{(m)})\big)(y) = \bfJ_j^{k,m}\Big({\sf\Lambda}^{\sfg}_y\big(\widehat{{\sf g}_{yx}}^\delta (\sigma^{(m)})\big)\Big)(y) = \sum_{\mu\le^+\sigma,\,\mu\notin\mcB_X^+}\sfg_{yx}(\sigma/^+\mu)\,\bfJ_j^{k,m}\big({\sf\Lambda}^{\sfg}_y(\mu^{(m)})\big)(y).
$$
\end{Dem}

\medskip

We can now prove that $\sfg_x(D^k\tau)\in L_\rap^\infty(\bbR^d)$ and ${\sfg}_{yx}(D^k\tau)\in\mcC^{|\tau|-|k|}_{(2),\rap}(\bbR^d\times\bbR^d)$, and close the induction step. We use the formulas from Lemma \ref{LemUsefulIdentities} for that purpose. First, since 
$$
\sfg_x(D^k\tau) = \sum_{\sigma\le^+\tau,\,\sigma\notin\mcB_X^+} \sfg_x(\tau/^+\sigma)\,\sff_x(D^k\sigma),
$$
with $\sfg_x(\tau/^+\sigma)\in L_\slow^\infty(\bbR^d)$ and $\sff_x(D^k\sigma)\in L_\rap^\infty(\bbR^d)$, from Lemma \ref{lem:I_j^k}, we have indeed $\sfg_x(D^k\tau)\in L_\rap^\infty(\bbR^d)$. Second, one can rewrite the identity
$$
\sfg_{yx}(D^k\tau) = \sum_{\sigma\le^+\tau,\,\sigma\notin\mcB_X^+} \sfg_{yx}(\tau/^+\sigma)\,\sff_y(D^k\sigma) - \sum_{\ell}\frac{(y-x)^\ell}{\ell!}\,\sff_x(D^{k+\ell}\tau),
$$
from Lemma \ref{LemUsefulIdentities}, using identity \eqref{formula:f_xD^k} for the $\sff$-terms. This gives for $\sfg_{yx}(D^k\tau)$ the formula
$$
\sum_j 
\Bigg\{
\sum_{\substack{\sigma\le^+\tau,\, \sigma\notin\mcB_X^+\\|k|<|\sigma|}} \sfg_{yx}(\tau/^+\sigma)\,\bfJ_j^{k,m}\big({\sf\Lambda}^{\sfg}_y(\sigma^{(m)})\big)(y) - \sum_{|k+\ell|<|\tau|}\frac{(y-x)^\ell}{\ell!}\,\bfJ_j^{k+\ell, m}\big({\sf\Lambda}^{\sfg}_x(\tau^{(m)})\big)(x)   
\Bigg\}
$$
\begin{equation*} \begin{split}
=: \sum_j\sfg_{yx}^j(D^k\tau).
\end{split} \end{equation*}
Given $x,y\in\bbR^d$, set $j_0=-1$, if $|y-x|\ge2$, and pick otherwise $j_0\ge-1$ such that $|y-x|\simeq2^{-j_0}$. One uses the first estimate from Lemma \ref{lem:I_j^k} to bound above the sum over $j\ge j_0$
\begin{equation} \label{EqEstimate} \begin{split}
&|x|_*^a\sum_{j\ge j_0}\big|\sfg_{yx}^j(D^k\tau)\big| \\
&\lesssim \sum_{j\ge j_0}\sum_{\substack{\sigma\le^+\tau,\,\sigma\notin\mcB_X^+\\|k|<|\sigma|}} |y-x|^{|\tau|-|\sigma|}\,2^{-j(|\sigma|-|k|)} + \sum_{j\ge j_0}\sum_{|k+\ell|<|\tau|}|y-x|^{|\ell|}\,2^{-j(|\tau|-|k|-|\ell|)} \\
&\lesssim \sum_{\substack{\sigma\le^+\tau\\|k|<|\sigma|}} |y-x|^{|\tau|-|\sigma|}\,2^{-j_0(|\sigma|-|k|)} + \sum_{|k+\ell|<|\tau|}|y-x|^{|\ell|}\,2^{-j_0(|\tau|-|k|-|\ell|)}   \\
&\lesssim|y-x|^{|\tau|-|k|}.
\end{split} \end{equation}
To consider the sum over $j<j_0$, assume now that $|y-x|<2$. Then, since $\sf (g, \Lambda)$ is a model and
\begin{align*}
{\sf\Lambda}^{\sfg}_x(\tau^{(m)}) = {\sf\Lambda}^{\sfg}_y\big(\widehat{{\sf g}_{yx}}^{\delta}\tau^{(m)}\big) = \sum_{\sigma\le^+\tau}\sfg_{yx}(\tau/^+\sigma){\sf\Lambda}^{\sfg}_y(\sigma^{(m)}),
\end{align*}
we have for $\sfg_{yx}^j(D^k\tau)$ the formula
\begin{align*}
\sfg_{yx}^j(D^k\tau)
&=\bfJ_j^{k,m}\big({\sf\Lambda}^{\sfg}_x(\tau^{(m)})\big)(y) - \sum_{\substack{\sigma\le^+\tau,\, \sigma\notin\mcB_X^+\\|k|>|\sigma|}} \sfg_{yx}(\tau/^+\sigma)\bfJ_j^{k,m}\big({\sf\Lambda}^{\sfg}_y(\sigma^{(m)})\big)(y) \\
&\quad- \sum_{|k+\ell|<|\tau|}\frac{(y-x)^\ell}{\ell!}\bfJ_j^{k+\ell, m}\big({\sf\Lambda}^{\sfg}_x(\tau^{(m)})\big)(x)   \\  
&= \lceil b\rceil\sum_{|k'|=\lceil b \rceil}\frac{(y-x)^{k'}}{{k'}!}\int_0^1(1-t)^{\lceil b\rceil}\bfJ_j^{k+{k'}}\big({\sf\Lambda}^{\sfg}_x(\tau^{(m)})\big)\big(x+t(y-x)\big)dt \\
&\quad- \sum_{\substack{\sigma\le^+\tau,\,\sigma\notin\mcB_X^+\\|k|>|\sigma|}}\sfg_{yx}(\tau/\sigma)\,\bfJ_j^{k,m}\big({\sf\Lambda}^\sfg_y(\sigma^{(m)})\big)(y)
\end{align*}
where $b := |\tau|-|k|$, by the multivariable Taylor remainder formula. Note that $|\tau|,|\sigma|\notin\bbN$ in the above formula, by assumption \textbf{\textsf{(C-3)}}. Since $|y-x|<2$, $|x+t(y-x)|_*\simeq|x|_*$. It follows then from Lemma \ref{lem:I_j^k} that $\sum_{-1\leq j<j_0}\big|\sfg_{yx}^j(D^k\tau)\big|$ is bounded above by
\begin{align*}
&\sum_{j<j_0}\sum_{|{k'}|=\lceil b\rceil}\sum_{\sigma\le^+\tau,\,\sigma\notin\mcB_X^+} |y-x|^{|{k'}|+|\tau|-|\sigma|}|x|_*^{-a}\,2^{-j(|\sigma|-|k|-|{k'}|)} \\
&\quad+ \sum_{j<j_0}\sum_{\substack{\sigma\le^+\tau,\,\sigma\notin\mcB_X^+\\|k|>|\sigma|}} |y-x|^{|\tau|-|\sigma|}|y|_*^{-a}\,2^{-j(|\sigma|-|k|)}   \\
&\lesssim |x|_*^{-a}\sum_{\sigma\le^+\tau,\,\sigma\notin\mcB_X^+} |y-x|^{|{k'}|+|\tau|-|\sigma|}\,2^{-j_0(|\sigma|-|k|-|{k'}|)} \\
&\quad+ |y|_*^{-a}\sum_{\substack{\sigma\le^+\tau,\,\sigma\notin\mcB_X^+\\|k|>|\sigma|}}|y-x|^{|\tau|-|\sigma|}\,2^{-j_0(|\sigma|-|k|)}   \\
&\lesssim\big(|x|_*^{-a}+|y|_*^{-a}\big)|y-x|^{|\tau|-|k|}.
\end{align*}
Together with inequality \eqref{EqEstimate}, the preceding upper bound tells us that ${\sfg}_{yx}(D^k\tau)\in\mcC^{|\tau|-|k|}_{(2),\rap}(\bbR^d\times\bbR^d)$. This closes the induction step.   \hfill $\rhd$

\bigskip

\noindent \textbf{\textsf{Remarks -- 1. On branched rough paths.}} The setting of $\bbR^\ell$-valued branched rough paths provides an example of regularity structure where Theorem \ref{mainthm2} applies, giving an alternative point of view on the results of Tapia and Zambotti in \cite{TZ18}. The Hopf algebra $(T^+,\Delta^+)$ is in that case the Butcher-Connes-Kreimer Hopf algebra. We recall the details for the reader as it also sets the scene for part of the results of Section \ref{SubsectionCorollaries}.

The set $T^+$ is the free commutative unital algebra generated by the set $\mcB^+_\circ$ of non-planar rooted trees with node decorations in a finite set $\{1,\dots,\ell\}$ and no decoration on the edges. The empty tree plays the role of the unit in $T^+$. A product of decorated trees is called a forest, so generic elements of $T^+$ are linear combinations of forests. The splitting map $\Delta^+$ is the algebra morphism defined on trees as follows. Given a labelled rooted decorated tree $\tau$, denote by $\textsf{Sub}(\tau)$ the set of subtrees of $\tau$ with the same root as $\tau$, and induced decoration. Given such a subtree $s$, we obtain a collection $\tau_1,\dots,\tau_n$ of decorated rooted trees by removing $s$ and all the adjacent edges to $s$ from $\tau$, and keeping the node decoration inherited from $\tau$. Write $\tau/s$ for the monomial $\tau_1\dots\tau_n$. One defines a linear multiplicative map $\Delta^+ : T^+ \rightarrow T^+\otimes T^+$, defining it on decorated trees by the formula
\[
\Delta^+\tau = \sum_{s\in\textsf{Sub}(\tau)} s\otimes (\tau/s),
\] 
An explicit formula for the antipode was first given by Connes and Kreimer in their celebrated work \cite{CK}; see \cite{CHV} for a simple and enlighting proof. Each node decoration $i\in\{1,\dots,\ell\}$ is assigned a homogeneity $\alpha\in(0,1)$,
and each decorated tree $\tau$ is equipped with the homogeneity $\alpha(\sharp\tau)$, where $\sharp\tau$ denotes the number of nodes contained in $\tau$. The homogeneity of a forest is the sum of the homogeneities of its decorated trees. It is elementary to check that $(T^+,\Delta^+)$ is indeed a Hopf algebra. To avoid polynomials and derivatives, we consider the subalgebra of trees with homogeneities smaller than $1$. Thus assumptions \textbf{\textsf{(A-C)}} without polynomials $X^k$ and derivatives $D^k$ hold, but it does not matter here.
Branched rough paths are $\sf g$-maps on $(T^+,\Delta^+)$ over a fixed time interval $[0,T]$ in place of $\bbR^d$. Theorem \ref{mainresult:reconst of g} applies then in this setting and provides a parametrization of the set of branched rough paths by the product space $\prod_{\tau\in\mcB_\circ^+} C^{\vert\tau\vert}([0,T], \bbR^\ell)$, in accordance with Tapia and Zambotti's main result, Theorem 1.2 and Corollary 1.3 in \cite{TZ18}. Our parametrization is different from their identification of the space of branched rough paths as a principal homogeneous space over the preceding product of H\"older spaces.  

(Theorem \ref{mainthm1} cannot be applied in a finite region $[0,T]$ directly. To overcome this point, we extend a function $f\in C([0,T], \bbR^\ell)$ to $[-T,T]$ symmetrically, and extend it to $\big[(2n-1)T, (2n+1)T\big]$ for any $n\in\bbZ$ periodically. Then for any $\alpha\in(0,1)$, the H\"older space $C^\alpha([0,T],\bbR^\ell)$ is identified with the space
$$
C_{p,s}^\alpha(\bbR, \bbR^\ell) := \Big\{ f\in C^\alpha(\bbR,\bbR^\ell)\,;\, f(t)=f(-t),\ f(t+2T)=f(t)\ \text{for any $t\in \bbR$} \Big\}.
$$
Note that Littlewood-Paley blocks $\Delta_i$ preserve the symmetry and periodicity, so such spaces are closed under paradifferential operators (paraproduct, its two-parameter extension, etc.) used in this paper. Hence we can apply Theorem \ref{mainthm1} to such spaces.)

\medskip

\textbf{\textsf{2. On the signature of arbitrary models.}} We mentioned in the introduction that admissible models on regularity structures built from integration operators have a well-defined signature -- that is a unique extension to the full regularity structure with elements of arbitrary large positive homogeneity. This comes from the fact that such models are determined uniquely by the definition of the $\sf \Pi$ map on elements of the regularity structure of negative homogeneity. Extending a regularity structure with additional elements of positive homogeneity the initial datum of the restriction of $\sf \Pi$ on the elements of negative homogeneity still defines a unique admissible model on the extended regularity structure. Such an automatic extension result does not hold for general models, with unrelated $\sf g$ and $\sf \Pi$ maps. Indeed, Theorem \ref{mainresult:reconst of g} tells us that the set of $\sf g$-maps is parametrized by a set of functions indexed by $\mcG^+_\circ$. Embedding a regularity structure into a larger regularity structure will a priori embed the set $\mcG^+_\circ$ into a larger set, implying the non-uniqueness of an extension of the $\sf g$-map, from Theorem \ref{mainresult:reconst of g} again. The following statement follows nonetheless from Theorem \ref{mainthm1} while it is beyond the scope of Theorem 21 in \cite{BH}. See Section 4 in \cite{BH} for the definition of admissible model.

\medskip

\begin{cor}
Let $\scrT$ be a regularity structure satisfying assumptions \textbf{\textsf{(A-C)}}. Let $\scrT'\subset\scrT$ be a sub-regularity structure of $\scrT$ satisfying these assumptions as well, and such that $\scrT'$ contains all the elements of $\scrT$ of negative homogeneity. Then any admissible model on $\scrT'$ has a unique extension into an admissible model on $\scrT$.
\end{cor}

\medskip

Like with Lyons' extension theorem, it is important to notice that the extension map is a continuous map. So even in a stochastic setting where the construction of a model may require stochastic analysis arguments, once this is done, the extension of this model to a larger structure no longer involves probability arguments.

\bigskip

\subsection{From paracontrolled systems to modelled distributions}
\label{SubsectionPCtoModelledDistributions}

We prove Theorem \ref{mainthm2} in this section. Let $\scrT$ be a regularity structure satisfying assumptions \textbf{\textsf{(A-C)}}. Pick $\gamma\in\bbR$, and ${\sf M = (g, \Pi)}\in\scrM_\rap(\scrT,\bbR^d)$.

\medskip

The key observation is that proving Theorem \ref{mainthm2} is equivalent to an extension problem for the map $\sfg$. Consider indeed the commutative algebra $T_{\bsF}^+$ generated by $\mcB^+$ and new symbols
\begin{align*}
(\bsF_\tau)_{\tau\in\mcB,\, |\tau|<\gamma}.
\end{align*}
Define the homogeneity of the symbol $\bsF_\tau$ by
$$
|\bsF_\tau| := \gamma - |\tau|.
$$ 
The coproduct $\Delta^+_{\bsF}$ on $T_{\bsF}^+$ extending $\Delta^+$ and such that
\begin{align}\label{eq:def of Delta F/tau}
\Delta^+(\bsF_\tau) = (\bsF_\tau)\otimes\unit+\sum_{\tau\le\mu}(\mu/\tau)\otimes(\bsF_\mu),
\end{align} 
is coassociative and turns $T^+_{\bsF}$ into a Hopf algebra. It satisfies assumptions \textbf{\textsf{(A)}} and \textbf{\textsf{(B)}} with 
$$
\mcB^+_{{\bsF},\circ} := \mcB^+_\circ\cup\Big\{{\bsF}_\tau\,;\,\vert\tau\vert<\gamma\Big\}
$$
in the role of $\mcB^+_\circ$.  
Note that $T_{\bsF}^+$ does not satisfy assumption \textbf{\textsf{(C)}} in general, since the $D^k\bsF_\tau$ have no reason to be independent from the $\{\bsF_\mu\}_\mu$. The elementary proof of the next statement is left to the reader.

\medskip

\begin{lem}\label{lem: model on Ftau is equiv to MD}
Given a family $(f_\tau)_{\tau\in\mcB}$ of continuous functions on $\bbR^d$, set $\bsf := \sum_{\tau\in\mcB}f_\tau\tau$, and
$$
\sfg_x(\bsF_\tau) := f_\tau(x) .
$$
Then
$$
\big\langle \tau', \bsf(y) - \widehat{\sfg_{yx}}\bsf(x)\big\rangle = \sfg_{yx}(\bsF_\tau).
$$
\end{lem}

\medskip

Defining a modelled distribution $\bsf\in\mcD_\rap^\gamma(T,\sf g)$ is thus equivalent to extending the map $\sfg$ from $T^+$ to $T_{\bsF}^+$ in such a way that the extended map on $(T_{\bsF}^+,\Delta^+_{\bsF})$ still satisfies the regularity constraints from Definition \ref{defi:weighted model}.

\medskip

Recall from assumption \textbf{\textsf{(B)}} that either $\mu/\tau\in\textrm{span}(\mcB^+\setminus\mcB^+_X)$ or $\mu/\tau\in \text{\rm span}(\mcB_X^+)$, for $\tau,\mu\in\mcB$. 
If $\mu/\tau\in \text{\rm span}(\mcB_X^+)$, set
$$
\mu/\tau =: \sum_{k\in\bbN^d} c_\tau^\mu(k)\,\frac{X^k}{k!},
$$
and define
\begin{align}\label{eq:k-antiderivative}
D^k\bsF_\tau := \sum_{\substack{\tau\leq \mu \\ \mu/\tau\in\text{\rm span}(\mcB_X^+)}} c_\tau^\mu(k) \, \bsF_\mu.
\end{align}
Then we have 
\begin{align*}
\Delta^+\bsF_\tau = \bsF_\tau\otimes\unit + \sum_{\substack{\tau\leq\mu \\ \mu/\tau\in\textrm{span}(\mcB^+\setminus\mcB^+_X)} }(\mu/\tau)\otimes \bsF_\mu + \sum_{k\in\bbN^d}\frac{X^k}{k!}\otimes D^k\bsF_\tau.
\end{align*}

\medskip

\begin{thm} \label{thm:PD to MD}
Let a concrete regularity structure $\scrT$ satisfying assumptions $\textbf{\textsf{(A-C)}}$ be given, together with a family $\big(\Brac{f_\tau}\in\mcC_\rap^{\gamma-|\tau|}(\bbR^d)\big)_{\tau\in\mcB,\,\vert\tau\vert<\gamma}$. Assume that $\gamma-|\tau|\notin\bbN$ for any $\tau\in\mcB$ with $|\tau|<\gamma$. Pick a model ${\sf (g,\Pi)}\in\scrM_\rap(\scrT,\bbR^d)$. Define
\begin{align*}
f_\tau := \underset{\mu/\tau\in\text{\rm span}(\mcB^+\setminus\mcB^+_X)}{\sum_{\tau\leq\mu, \,|\mu|<\gamma}} \sfP_{f_\mu}\Brac{\mu/\tau}^\sfg + \Brac{f_\tau},
\end{align*}
and
\begin{align}   \label{usual condition on MD}
f_{\tau}^{(k)}(x) := \partial_y^k\bigg\{f_\tau(y) - 
\underset{\mu/\tau\in\text{\rm span}(\mcB^+\setminus\mcB^+_X)}{\sum_{\tau\leq\mu, \,|\mu|<\gamma,\, |\mu/\tau|\le|k|}}
\sfg_{yx}(\mu/\tau)\,f_\mu(x)\bigg\}\Big|_{y=x}.
\end{align}
If the \textbf{\textsf{structure conditions}}
\begin{align}  \label{structure condition}
f_{\tau}^{(k)} = \underset{\mu/\tau\in\text{\rm span}(\mcB_X^+)}{\sum_{\tau\leq\mu, \,|\mu|<\gamma}} c_\tau^\mu(k) \, f_\mu,
\end{align}
holds for any $\tau\in\mcB$ and $k\in\bbN^d$, then
$$
\bsf = \sum_{\tau\in\mcB}f_\tau\tau\in\mcD_\rap^{\gamma}(T, \sfg).
$$
\end{thm}

\medskip

The structure condition is reminiscent of a condition introduced by Martin and Perkowski in \cite{MP18} to give a characterisation of modelled distributions in terms of Besov type spaces. Given that we see $f_\tau$ as $\sfg({\bsF}_\tau)$, formula \eqref{usual condition on MD} is nothing but a formula for ${\sfg}(D^k{\bsF}_\tau)$ -- the analogue of formula \eqref{eq:g_x(D^k tau)} in the present setting.

\medskip

\begin{Dem}
Consider the extended Hopf algebra $^\textsf{free}T_{\bsF}^+$ freely generated by the symbols
$$
\{X_1,\dots,X_d\}\cup\mcB_\circ^+\cup\Big\{D^k(\bsF_\tau)\,;\, \tau\in\mcB,\ \gamma>|\tau|+|k|\Big\}.
$$
It satisfies assumptions \textbf{\textsf{(A-C)}}. 
By Theorem \ref{mainresult:reconst of g} giving a paracontrolled parametrization of the map $\sfg$ by its definition on the $\sfg(\tau)$, with $\tau\in\mcG^+_{\bsF,\circ} := \mcG^+_\circ\cup\big\{{\bsF}_\tau\,;\,\vert\tau\vert<\gamma\big\}$, there exists a unique model $\sfg$ on $^\textsf{free}T_{\bsF}^+$ that coincides with $\sfg$ on $T^+$, and such that
$$
\sfg(\bsF_\tau) := \sum_{\tau\le\mu,\,|\mu|<\gamma}\sfP_{\sfg(\bsF_\mu)}\Brac{\mu/\tau}^{\sfg} + \Brac{f_\tau},
$$
for all $\tau\in\mcB$ with $\vert\tau\vert<\gamma$. Since $T_{\bsF}^+$ is the quotient space of $^\textsf{free}T_{\bsF}^+$ by the relations \eqref{eq:k-antiderivative}, and
$$
\sfg\big(D^k\bsF_\tau\big) = \underset{\mu/\tau\in\text{\rm span}(\mcB_X^+)}{\sum_{\tau\leq\mu, \,|\mu|<\gamma}} c_\tau^\mu(k)\,\sfg(\bsF_\mu),
$$
from the structure condition \eqref{structure condition}, the map $\sfg$ is consistently defined on the quotient space, where it satisfies the estimates from Definition \ref{defi:weighted model}.
\end{Dem}

\medskip

One can get rid of the structure condition in some cases.

\medskip

\noindent \textbf{\textsf{Assumption (D) --}} {\it For any $\tau\in\mcB_\bullet$, there is no term of the form $\sigma\otimes X^k$ with $k\neq0$, in the expansion of $\Delta\tau$.}

\medskip

Under assumption \textbf{\textsf{(D)}}, we can show that, given $\tau\in\mcB$, the only $\mu\geq \tau$ such that $\mu/\tau$ has a non-null component on $X^k$ is $\mu=\underline{X}^k\tau$. Indeed, writing $\mu=\underline{X}^\ell\sigma$ ($\ell\in\bbN^d,\,\sigma\in\mcB_\bullet$), by the multiplicativity of $\Delta$ and by \eqref{EqDelta}, we have
\begin{align*}
(\text{\rm Id}\otimes P_X)\Delta\mu
&=(\text{\rm Id}\otimes P_X)(\Delta\underline{X}^\ell)(\Delta\sigma)\\
&=(\text{\rm Id}\otimes P_X)(\Delta\underline{X}^\ell)(\sigma\otimes{\bf1}_+)
=\sum_{k}\binom{\ell}{k}\underline{X}^{\ell-k}\sigma\otimes X^{k}.
\end{align*}
Then $\mu/\tau\in\text{\rm span}(X^k)$ if and only if $\tau=\underline{X}^{\ell-k}\sigma$, thus $\mu=\underline{X}^k\tau$.
Then \eqref{eq:k-antiderivative} takes the form
$$
D^k\bsF_{\underline{X}^{\ell-k}\sigma}=\frac{\ell!}{(\ell-k)!}\bsF_{\underline{X}^\ell\sigma},\qquad\sigma\in\mcB_\bullet.
$$
Moreover, this reduces to the formula
$$
D^k\bsF_\sigma = k! \,\bsF_{\underline{X}^k\sigma},\qquad\sigma\in\mcB_\bullet,
$$
hence the structure condition \eqref{structure condition} takes the simple form \eqref{EqSimpleStructureCondition} below. Note that the data in the next statement is indexed by $\mcB_\bullet$, unlike in the general case of Theorem \ref{thm:PD to MD} where it is indexed by $\mcB$.

\medskip

\begin{cor}  \label{thm:PD to MD2}
Let $\scrT$ be a regularity structure satisfying assumptions $\textbf{\textsf{(A-D)}}$, and a family $\big(\Brac{f_\tau}\in \mcC_\rap^{\gamma-|\tau|}(\bbR^d)\big)_{\tau\in\mcB_\bullet,\,\vert\tau\vert<\gamma}$ be given. Assume that $\gamma-|\tau|\notin\bbN$ for any $\tau\in\mcB$ with $|\tau|<\gamma$. Pick a model ${\sf (g,\Pi)}\in\scrM_\rap(\scrT,\bbR^d)$. Set, for $\tau\in\mcB_\bullet$ with $\vert\tau\vert<\gamma$,
\begin{align*}
f_\tau := \underset{\mu/\tau\in\text{\rm span}(\mcB^+\setminus\mcB^+_X)}{\sum_{\tau\leq\mu, \,|\mu|<\gamma}} \sfP_{f_\mu}\Brac{\mu/\tau}^\sfg + \Brac{f_\tau},
\end{align*}
and, for $\tau\in\mcB_\bullet, k\in\bbN^d\backslash\{0\}$ with $\vert k\vert+\vert\tau\vert<\gamma$,
\begin{equation}  \label{EqSimpleStructureCondition}
f_{\underline{X}^k\tau}(x) :=
\frac1{k!}\,
\partial_y^k\bigg\{f_\tau(y) - 
\underset{\mu/\tau\in\text{\rm span}(\mcB^+\setminus\mcB^+_X)}{\sum_{\tau\leq\mu, \,|\mu|<\gamma,\, |\mu/\tau|\le|k|}}
\sfg_{yx}(\mu/\tau)\,f_\mu(x)\bigg\}\Big|_{y=x}.
\end{equation}
Then 
$$
\bsf := \sum_{\sigma\in\mcB,\, \vert\sigma\vert<\gamma}f_\sigma\,\sigma = \underset{\vert\tau\vert+\vert k\vert<\gamma}{\sum_{\tau\in\mcB_\bullet,\, k\in\bbN^d}} f_{\underline{X}^k\tau}\,\underline{X}^k\tau\in\mcD_\rap^{\gamma}(T, \sfg).
$$
\end{cor}

\medskip

Corollary \ref{thm:PD to MD2} yields the homeomorphism result from Theorem \ref{mainthm2}. As stated in the introduction, we can see the further homeomorphism result
\begin{align*}
\scrM_{\rap}\ltimes\mcD_{\rap}^\gamma
&:=\left\{
({\sf M},\bsf)\,;\, {\sf M}\in\scrM_{\rap}(\scrT,\bbR^d),\ \bsf\in\mcD_{\rap}^\gamma(T,{\sf g}) \right\}\\
&\simeq 
\prod_{\sigma\in\mcG_\circ^+} \mcC_\rap^{|\sigma|}(\bbR^d) \times \prod_{\tau\in\mcB_\bullet,\,|\tau|<0} \mcC_\rap^{|\tau|}(\bbR^d)
\times \prod_{\tau\in\mcB_\bullet,\, |\tau|<\gamma} \mcC_\rap^{\gamma-|\tau|}(\bbR^d),
\end{align*}
where the left hand side has a topology induced by the metrics
\begin{align*}
&d_a\big(({\sf M},\bsf),({\sf M}',\bsf') \big)\\
&:=d_a({\sf M},{\sf M}')
+\sup_{\tau\in\mcB}\left\| \big\langle \tau', 
(\bsf(y)-\widehat{{\sf g}_{yx}}\bsf(x))-(\bsf'(y)-\widehat{{\sf g}'_{yx}}\bsf'(x))
\big\rangle \right\|_{\mcC_a^{\gamma-|\tau|}(\bbR^d\times\bbR^d)}.
\end{align*}

\medskip

Note that assumption \textbf{\textsf{(D)}} is an assumption about the basis $\mcB$ of $T$ we choose to work with, not about the regularity structure itself. It is thus possible that a given basis satisfies assumption \textbf{\textsf{(D)}} whereas another does not. This flexibility is at the heart of the proof of Theorem \ref{mainthm3} in the next section.

\bigskip

\subsection{Modelled distributions over BHZ regularity structures}
\label{SubsectionBHZ-RS}

Bruned, Hairer and Zambotti introduced in \cite{BHZ} class of regularity structures convenient for the study of singular stochastic PDEs. We call these structures \textbf{\textsf{BHZ regularity structures}}
$$
\scrT_{\tiny \textsc{BHZ}} = \Big((T_{\BHZ}^+,\Delta^+_{\BHZ}), (T_{\BHZ},\Delta_{\BHZ})\Big).
$$ 
Although the canonical basis of these concrete regularity structures do not satisfy assumption \textbf{\textsf{(D)}} the following result holds true.

\medskip

\begin{thm} \label{prop:reduction of antiderivative}
Assume that the set of homogeneities $\{|\mathfrak{t}|\}_{\mathfrak{t}\in\mathfrak{L}}\cup\{1\}$ is rationally independent, that is, the only tuple of integers $\{k_{\mathfrak{t}}\}_{\frak{t}}\cup\{k_1\}$ such that $\sum_{\mathfrak{t}}k_{\mathfrak{t}}|\mathfrak{t}|+k_1=0$ is the trivial soluion $k_{\mathfrak{t}}=k_1=0$. Then the canonical bases $\mcB_{\BHZ}^+$ and $\mcB_{\BHZ}$ satisfy assumptions \textbf{\textsf{(A-C)}}.
Moreover, one can construct a basis of $T_{\BHZ}$ that satisfies assumptions \textbf{\textsf{(A-D)}}.
\end{thm}

\medskip

The remaining of this section is dedicated to proving this statement. We recall first the elements of the construction of BHZ regularity structures that we need here. These concrete regularity structures are indexed by decorated rooted trees.

Any finite connected graph without loops and with a distinguished vertex is called a rooted tree. For any rooted tree $\tau$, denote by $N_\tau$ the node set, by $E_\tau$ the edge set, by $\varrho_\tau\in N_\tau$ the distinguished vertex, called root of $\tau$. Let also $\frkL$ be a finite set of types. (Edges will be interpreted differently depending on their type, when given any model on $\scrT_{\BHZ}$. Different types may for instance correspond to different convolution operators.) Denote by $\bsmcB$ the set of rooted decorated trees. Each $\tau\in\bsmcB$ is a rooted tree equipped with the type map $\frkt:E_\tau\to\frkL$ and with the decorations
\begin{itemize}
\item $\frkn:N_\tau\to\bbN^d$.
\item $\frko:N_\tau\to\bbZ^d\oplus\bbZ(\frkL)$.
\item $\frke:E_\tau\to\bbN^d$.
\end{itemize}
Equivalently, the set $\bsmcB$ is generated recursively by the application of the following operations -- see \cite[Section 4.3]{BHZ}.
\begin{itemize}
   \item One has $\bullet^k\in\bsmcB$ for any $k\in\bbN^d$, where $\bullet^k$ is a tree with only one node $\bullet$, with $\frkn(\bullet)=k$, and $\frak{o}(\bullet)=0\oplus0$.
   \item If $\tau,\sigma\in\bsmcB$ then $\tau\sigma\in\bsmcB$, where $\tau\sigma$ is called a \emph{tree product}; $\tau\sigma$ is a graph $\tau\sqcup\sigma$ divided by the equivalence relation $\sim$ on $N_\tau\sqcup N_\sigma$, where $x\sim y$ means $x=y$ or $x,y\in\{\varrho_\tau,\varrho_\sigma\}$. On the root $\varrho_{\tau\sigma}$, the decorations $\frkn(\varrho_{\tau\sigma})=\frkn(\varrho_\tau)+\frkn(\varrho_\sigma)$ and $\frko(\varrho_{\tau\sigma})=\frko(\varrho_\tau)+\frko(\varrho_\sigma)$ are given.
   \item For any $\frkt\in\frkL$ and $k\in\bbN^d$,
$$
\tau\in\bsmcB\quad\Rightarrow\quad I_k^\frkt(\tau)\in\bsmcB,
$$
where the tree $I_k^\frkt(\tau)$ is obtained by adding on $\tau$ one distinguished node $\varrho'$ and one edge $e=(\varrho_\tau,\varrho')$ of type $\frkt$, with decorations $\frke(e)=k$ and $\frko(\varrho')=0\oplus0$.
   \item For any $\alpha\in\bbZ^d\oplus\bbZ(\frkL)$, denote by $R_\alpha$ the operator on decorated rooted trees adding a value $\alpha$ on the decoration $\frko$ on $\varrho_\tau$. Assume
$$
\tau\in\bsmcB\quad\Rightarrow\quad R_\alpha(\tau)\in\bsmcB.
$$
By applying the operator $R_\alpha$ with various $\alpha$ on each step as above, one can see that, if $\tau\in\bsmcB$ then the same decorated tree with any other $\frko$-decorartion is also an element of $\bsmcB$.
\end{itemize}

\medskip

Each type $\frkt\in\frkL$ is assigned a nonzero real number $|\frkt|$, the collection of which satisfies the assumption of Theorem \ref{prop:reduction of antiderivative}. One assigns a homogeneity $\vert\frak{n}\vert, \vert\frak{o}\vert, \vert\frak{e}\vert, \vert\frak{t}\vert$ to the decorations and edge types of any decorated tree $\tau$, and set
\begin{align*}
\vert\tau\vert &:= \vert\frak{n}\vert + \vert\frak{o}\vert - \vert\frak{e}\vert + \vert\frak{t}\vert\\
&:=\sum_{n\in N_\tau}|\frak{n}(n)|+\sum_{n\in N_\tau}|\frak{o}(n)|-\sum_{e\in E_\tau}|\frak{e}(e)|+\sum_{n\in N_\tau}|\frak{t}(n)|,
\end{align*}
where $|a+\sum_{\frkt}a_\frkt\frkt|:=|a|+\sum_{\frkt}a_\frkt|\frkt|$ for $a+\sum_{\frkt}\frkt\in\bbZ^d\oplus\bbZ(\frkL)$.
A noise-type  object $\Theta$ is represented by $I_0^{\frak{t}}(\bullet^0)$, with $\frak{t}$ of negative homogeneity.

\ssk

With each subcritical singular stochastic PDE is associated a notion of conforming and strongly conforming decorated tree. The basis $\mcB_{\BHZ}$ of $T_{\BHZ}$ is made up of the set of elements of $\bsmcB$ that strongly conforms the rule (see Section 5 in \cite{BHZ}), and the basis $\mcB_{\BHZ}^+$ of $T_{\BHZ}^+$ is made up of the elements of the form
$$
\bullet^k\prod_{i=1}^NI_{k_i}^{\frak{t}_i}(\tau_i),
$$
where $k,k_i\in\bbN^d$, $\frak{t}_i\in\frak{L}$, $\tau_i\in\mcB_{\BHZ}$, and $|I_{k_i}^{\frak{t}_i}(\tau_i)|>0$.
Such a tree is said to conform the rule. We do not need more details here and refer the interested reader to Section 5 of \cite{BHZ}. We do not describe in particular the details of the definition of the splitting maps $\Delta_\BHZ$ and $\Delta^+_\BHZ$; we only record the following fact, where we write $\unit$ for $\bullet^0$, and $X^k$ for $\bullet^k$. 

\medskip

\begin{prop}\cite[Proposition 4.17]{BHZ}\label{recursive def of Delta}
The coproduct $\Delta=\Delta_{\BHZ}:T_\BHZ\to T_\BHZ\otimes T^+_\BHZ$, satisfies the following identities
\begin{align*}
&\Delta\unit=\unit\otimes\unit,\quad \Delta X_i=X_i\otimes\unit+\unit\otimes X_i,\quad \Delta(\tau\sigma)=(\Delta\tau)(\Delta\sigma),  \\
&\Delta I_k^\frkt(\tau) = \big(I_k^\frkt\otimes\iden\big)\Delta\tau + \sum_{\vert\ell\vert+\vert k\vert <\vert \tau\vert+\vert\frak{t}\vert}\frac{X^\ell}{\ell!}\otimes I_{k+\ell}^\frkt(\tau),\quad
\Delta R_\alpha(\tau) = \big(R_\alpha\otimes\iden\big)\Delta\tau.
\end{align*}
The coproduct $\Delta^+=\Delta_{\BHZ}^+:T^+_\BHZ\to T^+_\BHZ\otimes T^+_\BHZ$, satisfies the same identities with $\Delta$ in the right hand sides replaced by $\Delta^+$.
\end{prop}

\medskip

\begin{thm}\label{BHZ-ABC}
The bases $\mcB=\mcB_{\BHZ}$ and $\mcB^+=\mcB_{\BHZ}^+$ satisfy assumptions \textbf{\textsf{(A-C)}}.
\end{thm}

\medskip

\begin{Dem}
Assumption \textbf{\textsf{(A)}} is satisfied by setting
$$
\mcB_\circ^+:=\{I_k^{\frak{t}}(\tau)\in\mcB^+\,;\, \frak{t}\in\frak{L},\, k\in\bbN^d,\, \tau\in\mcB\},\qquad
\mcB_\bullet:=\{\tau\in\mcB\,;\, \frak{n}(\varrho_\tau)=0\},
$$
and $X^k=\underline{X}^k=\bullet^k$.
Assumption \textbf{\textsf{(B)}} follows because polynomial elements and non-polynomial elements are distinguished by the number of their edges. Indeed, $\sharp E_\tau=0$ if and only if $\tau\in\mcB_X^+=\mcB_{\underline{X}}$.
Assumption \textbf{\textsf{(C-1)}} is satisfied by setting
$$
\mcG_\circ^+:=\{I_0^{\frak{t}}(\tau)\in\mcB^+\,;\, \frak{t}\in\frak{L},\, \tau\in\mcB\}.
$$
Then $D^kI_0^{\frak{t}}(\tau)=I_k^{\frak{t}}(\tau)$ follows from Proposition \ref{recursive def of Delta}.
To check \textbf{\textsf{(C-2)}}, we define the binary relation on $\mcB^+$ by denoting $\sigma\pord\tau$ if
\begin{itemize}
\item $\sharp E_\sigma\le\sharp E_\tau$, or
\item $\sharp E_\sigma=\sharp E_\tau$ and $|\frak{n}_\sigma|\le|\frak{n}_\tau|$, where $\frak{n}_\tau$ (resp. $\frak{n}_\sigma$) denotes the $\frak{n}$-decoration given for $\tau$ (resp. $\sigma$).
\end{itemize}
This relation is transitive and satisfies the first condition of \textbf{\textsf{(C-2)}}.
The second one in \textbf{\textsf{(C-2)}} follows from the graphical definition of $\Delta^+$ -- see Section 2 in \cite{BHZ} for details. Essentially, we have the decomposition
$$
\Delta^+\tau=\sum\sigma\otimes(\tau/^+\sigma),
$$
where either of the following holds.
\begin{itemize}
\item $\sigma$ is the same graph as $\tau$ but with $\frak{n}_\sigma\le\frak{n}_\tau$.
$\tau/^+\sigma$ consists of only one node.
\item $\sigma$ is a strict subtree of $\tau$ such that $\varrho_\sigma=\varrho_\tau$, and $\tau/^+\sigma$ is a quotient graph of $\tau$ obtained by contracting the subgraph $\sigma$ into one node.
\end{itemize}
For the first case, if $\frak{n}_\sigma=\frak{n}_\tau$ then $\sigma=\tau$ as an element of $\bsmcB$, and if $\frak{n}_\sigma<\frak{n}_\tau$ then $\sigma\spord\tau$.
For the second case, if $\sharp N_\sigma=0$ then $\sigma$ is a polynomial and if $\sharp N_\sigma>0$ then $\sigma,\tau/^+\sigma\spord\tau$.
Hence the formula \eqref{eq:canonical coproduct} holds.
For the last assumption, since the set $\{|\mathfrak{t}|\}_{\mathfrak{t}\in\mathfrak{L}}\cup\{1\}$ is rationally independent, non-polynomial $\tau$ (hence $\tau$ has at least one edge) has non-integer homogeneity. Hence \textbf{\textsf{(C-3)}} holds.
\end{Dem}

\medskip

The canonical bases $\mcB_{\BHZ}$ of BHZ concrete regularity structures do not satisfy assumption \textbf{\textsf{(D)}} since one has
$$
\Delta I^{\frak{t}}_0(X_i\Theta) = I^{\frak{t}}_0(X_i\Theta)\otimes\unit + I^{\frak{t}}_0(\Theta)\otimes X_i + \sum_{\vert k\vert< \vert\Theta\vert+1+\vert\frak{t}\vert}\frac{X^k}{k!}\otimes I^{\frak{t}}_k(X_i\Theta),
$$
for any edge type $\frak{t}$ with positive homogeneity, but the second term in the right hand side contradicts to assumption \textbf{\textsf{(D)}}. 
Now we define another basis of $T_{\BHZ}$. Set 
$$
\bsT := \textrm{span}(\bsmcB).
$$
The tree product $(\tau,\sigma)\mapsto\tau\sigma$ and the operators $I_k^\frkt$ and $R_\alpha$ are linearly extended to $\bsT$. For any $\frkt\in\frkL$ and $k,\ell\in\bbN^d$, we define the new operator ${}_\ell I_k^\frkt:\bsT\to \bsT$, by
$$
{}_\ell I_k^\frkt(\tau) := \sum_{m\in\bbN^d}\binom{\ell}{m}X^m(-1)^{\ell-m}I_k^\frkt\big(X^{\ell-m}\tau\big).
$$
(An operator ${}_\ell I_k$ represents the convolution with a kernel $x^\ell(\partial^k K)(x)$. These operators also appeared in the very recent work \cite{HP19} of Hairer and Pardoux.) If $\tau$ is homogeneous, then ${}_\ell I_k^\frkt(\tau)$ is also homogeneous and
$$
\big|{}_\ell I_k^\frkt(\tau)\big| = |\frkt| - |k| + |\ell | + |\tau|.
$$

\medskip

\begin{lem}
Consider the subset $\widetilde{\bsmcB}_\bullet\subset \bsT$ generated by the following rules.
\begin{itemize}
\item $\unit\in\widetilde{\bsmcB}_\bullet$.
\item $\tau\in\widetilde{\bsmcB}_\bullet$ $\Rightarrow$ ${}_\ell I_k^\frkt(\tau)\in\widetilde{\bsmcB}_\bullet$.
\item $\tau\in\widetilde{\bsmcB}_\bullet$ $\Rightarrow$ $R_\alpha(\tau)\in\widetilde{\bsmcB}_\bullet$.
\item $\tau,\sigma\in\widetilde{\bsmcB}_\bullet$ $\Rightarrow$ $\tau\sigma\in\widetilde{\bsmcB}_\bullet$.
\end{itemize}
Set
$$
\widetilde{\bsmcB} := \Big\{X^k\tau\,;\, k\in\bbN^d, \tau\in\widetilde{\bsmcB}_\bullet\Big\}.
$$
Then $\widetilde{\bsmcB}$ is a linear basis of $\bsT$, and there exists a basis $\widetilde{\mcB}=\widetilde{\mcB}_{\BHZ}$ of $T_{\BHZ}$ such that $\widetilde{\mcB}\subset\widetilde{\bsmcB}$.
\end{lem}

\medskip

\begin{Dem}
Assume that $\tau\in\bsmcB$ is expanded by the basis $\widetilde{\bsmcB}$, that is, $\tau$ is of the form
$$
\tau=\sum_i a_i X^{k_i}\sigma_i
$$
with $a_i\in\bbR$, $k_i\in\bbN^d$, and $\sigma_i\in\widetilde{\bsmcB}_\bullet$. Since the commutative property $R_\alpha(X^k\cdot)=X^kR_\alpha(\cdot)$ holds by the definition, $R_\alpha(\tau)$ is also expanded by $\widetilde{\bsmcB}$. By the inversion formula
\begin{align*}
I_k^\frkt\big(X^\ell\sigma\big)=\sum_{m\in\bbN^d}\binom{\ell}{m}X^m(-1)^{\ell-m}\,{}_{\ell-m}I_k^\frkt(\sigma),
\end{align*}
$I_k^\frkt(\tau)$ is also expanded by $\widetilde{\bsmcB}$. Certainly, if $\tau,\sigma\in\textrm{span}(\widetilde{\bsmcB})$, then $\tau\sigma\in\textrm{span}(\widetilde{\bsmcB})$. We can conclude that $\bsT=\textrm{span}(\widetilde{\bsmcB})$ by the induction on the number of edges on $\tau$.

As in the definition of $\mcB_\BHZ$ from $\bsmcB$, one obtains $\widetilde\mcB$ by keeping only those elements from $\widetilde{\bsmcB}$ that strongly conforms. 
\end{Dem}

\medskip

The set $\widetilde{\bsmcB}$ can be encoded as a set of rooted decorated trees using different decorations from the preceding decorations. Each $\tau\in\widetilde{\bsmcB}_\bullet$ is represented by a rooted tree with $\frko$ and $\frke$ decorations, together with a new decoration
$$
\frkf:E_\tau\to\bbN^d.
$$
The map ${}_{\ell}I_k^\frkt:\tilde{\bsmcB}_\bullet\to\tilde{\bsmcB}_\bullet$, is defined as follows. For any $\tau\in\widetilde{\bsmcB}_\bullet$ with root $\varrho$, the tree ${}_{\ell}I_k^\frkt(\tau)$ is obtained by adding to $\tau$ one node $\varrho'$ and one edge $e:=(\varrho,\varrho')$, with decorations $\frke(e)=k$ and $\frkf(e)=\ell$. Each $\tau=X^k\sigma\in\widetilde{\bsmcB}$ is represented by a rooted tree with decorations $\frkn,\frko,\frke,\frkf$, where $\frkn$ vanishes at any node except the root, where it is equal to $k$. We call this tree representation of elements of $\widetilde{\bsmcB}$ the \textit{non-canonical representation}. Note that, under such exchange of representations, the shape of trees is preserved.

\medskip

\begin{thm}  \label{prop:reduction of antiderivative}
The basis $\widetilde{\mcB}$ of $T_{\BHZ}$ satisfies assumption \textbf{\textsf{(D)}}, where $\widetilde{\mcB}_\bullet=\widetilde{\bsmcB}_\bullet\cap\widetilde{\mcB}$.
\end{thm}

\medskip

\begin{Dem}
The proof is done by the induction on the number of edges on $\tau$ in its non-canonical representation. In fact, one can conclude a stronger claim; for any $\tau\in\widetilde{\mcB}_\bullet$, one has
\begin{align}\label{eq:stronger than D}
\Delta\tau=\sum_{\sigma\in\widetilde{\mcB}_\bullet,\, \eta\in \widetilde{\mcB}\setminus\text{\rm span}\{X^k\}_k} c_{\sigma\eta}^\tau\sigma\otimes\eta.
\end{align}
It is sufficient to show that, if the coproduct of $\tau\in\widetilde{\mcB}_\bullet$ has such a form, then ${}_{\ell}I_k^\frkt(\tau)$ also satisfies the same condition. To complete the proof, we compute explicitly the coproduct $\Delta( {}_{\ell}I_k^\frkt(\tau))$.
Since
\begin{align*}
\Delta I_k^\frkt\big(X^a\tau\big)
&=(I_k^\frkt\otimes\iden)\Delta\big(X^a\tau\big)
+\sum_{\ell\in\bbN^d}\frac{X^{\ell}}{\ell!}\otimes I_{k+\ell}^\frkt\big(X^a\tau\big)\\
&=\sum_{\sigma\le\tau,\, b\in\bbN^d}\binom{a}{b}
I_k^\frkt \big(X^b\sigma\big)\otimes X^{a-b}(\tau/\sigma)
+\sum_{\ell\in\bbN^d}\frac{X^{\ell}}{\ell!}\otimes I_{k+\ell}^\frkt\big(X^a\tau\big),
\end{align*}
we have
\begin{align*}
\Delta\big({}_{a} I_k^\frkt(\tau)\big)
&=\sum_{b\in\bbN^d}\binom{a}{b}
\big(\Delta X^b\big)(-1)^{a-b}\Delta I_k\big(X^{a-b}\tau\big)\\
&=\sum_{\sigma\le\tau,\, b,c,d\in\bbN^d}(-1)^{a-b}
\binom{a}{b}\binom{b}{c}\binom{a-b}{d}
X^c I_k^\frkt\big(X^d\sigma\big)
\otimes
X^{b-c} X^{a-b-d}(\tau/\sigma)\\
&\quad+\sum_{\ell,b,c\in\bbN^d}(-1)^{a-b}\binom{a}{b}\binom{b}{c}
X^c\frac{X^\ell}{\ell!}
\otimes
X^{b-c}I_{k+\ell}^\frkt\big(X^{a-b}\tau\big)\\
&=:({\sf i})+({\sf ii}).
\end{align*}
The term $({\sf ii})$ does not contain any terms of the form $\sigma\otimes X^k$ with $k\neq0$. The sum $(\sf i)$ is equal to

\begin{align*}
&\sum_{\substack{\sigma\le\tau \\ a=c+c'+d+d'}}(-1)^{d+d'}
\frac{a!}{c!c'!d!d'!}
X^c I_k^\frkt\big(X^d\sigma\big)
\otimes
X^{c'} X^{d'}(\tau/\sigma) \\
&=\sum_{\substack{\sigma\le\tau \\ a=\alpha+\beta}}\frac{a!}{\alpha!\beta!}
\Bigg(\sum_{\alpha=c+d}(-1)^{d}
\frac{\alpha!}{c!d!}
X^c I_k^\frkt\big(X^d\sigma\big)\Bigg)
\otimes
\Bigg(\sum_{\beta=c'+d'}(-1)^{d'}\frac{\beta!}{c'!d'!}
X^{c'} X^{d'}(\tau/\sigma)\Bigg)\\
&=\sum_{\substack{\sigma\le\tau \\ a=\alpha+\beta}}\binom{a}{\alpha}
{}_{\alpha} I_k^\frkt(\sigma)
\otimes
(X-X)^\beta(\tau/\sigma)\\
&=\sum_{\sigma\le\tau}
{}_{a} I_k^\frkt(\sigma)
\otimes
(\tau/\sigma)
=({}_{a} I_k^\frkt\otimes\iden)\Delta\tau.
\end{align*}
Since $\tau$ is assumed in the induction step to have a coproduct \eqref{eq:stronger than D}, hence  $\Delta( {}_{\ell} I_k^\frkt(\tau))$, enjoys the same property. 
\end{Dem}


\subsection{Density and extension corollaries} 
\label{SubsectionCorollaries}

Corollaries \ref{CorDensity}, \ref{CorLyonsVictoir}, \ref{CorExtension}, and  \ref{CorDensitySmoothModelled} are proved as follows. Note that Schwartz space $\mcS(\bbR^d)$ is dense in the space $C_{\rap}^\beta(\bbR^d)$ in the topology of $C_{\rap}^{\beta-\epsilon}(\bbR^d)$ for any $\epsilon>0$; for any $f\in C_{\rap}^\beta(\bbR^d)$, the function $e^{t\Delta}f$ belongs to $\mcS(\bbR^d)$ and satisfies
$$
\|e^{t\Delta}f-f\|_{C_a^{\beta-\epsilon}(\bbR^d)}\lesssim t^{2\epsilon}\|f\|_{C_a^\beta(\bbR^d)}
\xrightarrow{t\to0}0
$$
for any $a>0$. See e.g. Proposition 3.12 of \cite{MW}.

\medskip

\begin{Dem}[of Corollary \ref{CorDensity}]
By Theorem \ref{mainthm1}, the space $\scrM_{\rap}(\scrT,\bbR^d)$ is homeomorphic to the product space
\begin{align*}
\prod_{\sigma\in\mcG_\circ^+} \mcC_\rap^{|\sigma|}(\bbR^d) \times \prod_{\tau\in\mcB_\bullet,\,|\tau|<0} \mcC_\rap^{|\tau|}(\bbR^d).
\end{align*}
For any $\epsilon>0$, any elements of this space can be approximated by smooth elements in the topology of the same space with each exponent $|\tau|$ replaced by $|\tau|-\epsilon$. By the formulas \eqref{def:bracket g} and \eqref{def:bracket Pi}, it turns out that a smooth element of \eqref{EqParametrizationThm1} is transferred to a smooth model in $\scrM_{\rap}(\scrT,\bbR^d)$.
\end{Dem}

\ssk

The proof of Corollary \ref{CorDensitySmoothModelled} is completely parallel and left to the reader.

\bigskip

\begin{Dem}[of Corollaries \ref{CorLyonsVictoir} and \ref{CorExtension}]
For Corollary \ref{CorLyonsVictoir}, consider the algebra $T_+$ generated by the set $\mcB_\circ^+$ of rooted trees as in Remark 1 of the previous section. Given an $\bbR^\ell$-valued $\alpha$-H\"older function $h=(h_i)_{i=1}^d$, a lift of the control $h$ is a branched rough path $(H^\tau)_{\tau\in\mcB_\circ^+}$ such that $H^{\bullet_i}=h_i$, where $\bullet_i$ denotes a graph with only one node and with node decoration $i$. By Theorem \ref{mainthm1}, such a lift is transferred to an elements of the product space $\prod_{\tau\in\mcB_\circ^+} \mcC_\rap^{|\tau|}(\bbR)$ such that $\Brac{\bullet_i}=h_i$. A trivial extension is defined by $\Brac{\tau}=0$ if $\sharp\tau\ge2$, and the associated model is nothing but a trivial lift of $h$.

\smallskip

Corollary \ref{CorExtension} is proved by a similar argument. By admissibility, the set $\scrM(\scrT,\bbR^m)$ is homeomorphic to the space $\prod_{\tau\in\mcB_\bullet, |\tau|<0} \mcC_\rap^{|\tau|}(\bbR^m)$.
Given a multi-dimensional noise $(\zeta_j)_{j=1}^\ell$, a trivial extension is defined by
$$
\Brac{\tau}^{\sfg}=
\begin{cases}
\zeta_j, & \tau=\bullet_j,\\
0, & \text{otherwise}.
\end{cases}
$$
\end{Dem}

\bigskip

\appendix
\section{Concrete regularity structures}
\label{SectionAppendixConcreteRS}

We recall in this appendix the setting of concrete regularity structures introduced in \cite{BH}, and refer the reader to Section 2 of  \cite{BH} for motivations for the introduction of that setting.

\medskip

\begin{defn*}
A \textsf{\textbf{concrete regularity structure}} $\mathscr{T}=(T^+,T)$ is the pair of graded vector spaces   \vspace{-0.1cm}
$$
T^+ = \bigoplus_{\alpha\in A^+} T_\alpha^+, \qquad T = \bigoplus_{\beta\in A} T_\beta   \vspace{-0.1cm}
$$
such that the following holds.   \vspace{0.1cm}
\begin{itemize}
   \item The index set $A^{+}\subset \bbR_+$ contains the point $0$, and $A^++A^+ \subset A^+$; the index set $A\subset\bbR$ is bounded below, and both $A^+$and $A$ have no accumulation points in $\bbR$. Set 
$$
\beta_0 := \min A.
$$  
   
   \item The vector spaces $T_\alpha^+$ and $T_\beta$ are finite dimensional. \vspace{0.1cm}
   
   \item The set $T^+$ is an algebra with unit $\bf 1_+$, with an algebra morphism
   $$
   \Delta^+ : T^+\rightarrow T^+\otimes T^+,
   $$
   such that $\Delta^+{\bf 1}_+ = {\bf 1}_+\otimes {\bf 1}_+$, and, for $\tau\in T_\alpha^+$,
   \begin{equation}   \label{EqDefnDeltaPlus}
   \Delta^+\tau \in \left\{\tau\otimes {\bf 1}_+ + {\bf 1}_+\otimes \tau + \sum_{0<\beta<\alpha} T^+_\beta\otimes T^+_{\alpha-\beta}\right\},
   \end{equation}
   and $\Delta^+$ satisfies the coassociativity property
$$
(\Delta^+\otimes\text{\rm Id})\Delta^+ = (\text{\rm Id}\otimes\Delta^+)\Delta^+.
$$
That is, $T^+$ has a Hopf structure with coproduct $\Delta^+$ and counit ${\bf1}_+'$.  \vspace{0.1cm}

   \item One has $T_0^+ = \textrm{\emph{span}}({\bf 1}_+)$, and for any $\alpha,\beta\in A^+$, one has $T_\alpha^+ T_\beta^+ \subset T_{\alpha+\beta}^+$.   \vspace{0.1cm}
   
   \item One has a linear splitting map 
   $$
   \Delta : T \rightarrow T\otimes T^+,
   $$ 
   of the form
   \begin{equation}
   \label{EqDelta}  
   \Delta \tau \in \left\{\tau\otimes {\bf 1}_+ + \sum_{\beta<\alpha} T_\beta\otimes T^+_{\alpha-\beta}\right\}   
   \end{equation}
   for each $\tau\in T_\alpha$, with the right comodule property
   \begin{equation*}
   \big(\Delta\otimes \textrm{\emph{Id}}\big)\Delta = \big(\textrm{\emph{Id}}\otimes \Delta^+\big)\Delta.
   \end{equation*}
\end{itemize}
Let $\mcB_\alpha^+$  and $\mcB_\beta$ be bases of $T_\alpha^+$ and $T_\beta$, respectively. We assume $\mcB_0^+=\{{\bf1}_+\}$. Set 
$$
\mcB^+ := \bigcup_{\alpha\in A^+} \mcB_\alpha^+, \qquad \mcB := \bigcup_{\beta\in A} \mcB_\beta.
$$ 
An element $\tau$ of $T_\alpha^{(+)}$ is said to be \emph{homogeneous} and is assigned \textsf{\textbf{homogeneity}} $|\tau| := \alpha$. The homogeneity of a generic element $\tau\in T^{(+)}$ is defined as $|\tau| := \max\{\alpha\}$, such that $\tau$ has a non-null component in $T^{(+)}_\alpha$.
We denote by 
$$
\mathscr{T} := \big((T^+,\Delta^+), (T,\Delta)\big)
$$
a concrete regularity structure.
\end{defn*}

\medskip

One of the elementary and important examples is the Taylor polynomial ring. Consider symbols $X_1, \dots, X_d$ and set 
$$
T_X := \bbR[X_1,\dots, X_d].
$$ 
For a multi index $k=(k_i)_{i=1}^d\in\bbN^d$, we use the notation
$$
X^k := X_1^{k_1}\cdots X_d^{k_d}.
$$
We define the homogeneity $|X^k|=|k|:=\sum_ik_i$, and the coproduct
\begin{equation}  \label{EqCoproductPolynomials}
\Delta X_i=X_i\otimes\unit+\unit\otimes X_i.
\end{equation}
Then $\big((T_X,\Delta),(T_X,\Delta)\big)$ is a concrete regularity structure.

\ssk

The set $G^+$ of characters $g : T^+\to\bbR$, i.e. nonzero algebra morphisms, forms a group with the convolution product
$$
g_1*g_2:=(g_1\otimes g_2)\Delta^+.
$$

\vfill\pagebreak

\section{Technical estimates}
\label{SectionAppendix}

We provide in this appendix a number of technical estimates that are variations on the corresponding results from \cite{BH}. Proofs are given for completeness.

\medskip

\begin{lem}\label{lem:weighted scaling}
If $\alpha\ge0$ and $a\in\bbZ$, then
\begin{align*}
\int \big|P_i(x-y)\big| |x-y|^\alpha|y|_*^{-a}dy&\lesssim2^{-i\alpha}\,|x|_*^{-a},   \\
\int \big|Q_i(x-y)\big| |x-y|^\alpha|y|_*^{-a}dy&\lesssim2^{-i\alpha}\,|x|_*^{-a}.
\end{align*}
\end{lem}

\medskip

\begin{Dem}
Recall the inequalities in the beginning of Section \ref{SectionFunctionalSetting}. If $a\ge0$,
\begin{align*}
|x|_*^a\int \big|P_i(x-y)\big| |x-y|^\alpha|y|_*^{-a}dy &\lesssim\int \big|P_i(x-y)\big| |x-y|^\alpha|x-y|_*^a dy = \int \big|P_i(y)\big| |y|^\alpha|y|_*^a dy  \\
&= \int \big|P_0(y)\big| \Big|\frac{y}{2^i}\Big|^\alpha\Big|\frac{y}{2^i}\Big|_*^ady \le 2^{-i\alpha}\int \big|P_0(y)\big| |y|^\alpha |y|_*^ady\lesssim2^{-i\alpha}.
\end{align*}
If $a<0$,
\begin{align*}
\int \big|P_i(x-y)\big| |x-y|^\alpha |y|_*^{-a}dy &\lesssim|x|_*^{-a}\int \big|P_i(x-y)\big| |x-y|^\alpha|x-y|_*^{-a}dy \lesssim 2^{-i\alpha}\,|x|_*^{-a}.
\end{align*}
\end{Dem}

\medskip

As a consequence of Lemma \ref{lem:weighted scaling}, we have the inequality
\begin{align*}
\|\Delta_jf\|_{L_a^\infty}
&\le\sup_x |x|_*^a\int|Q_j(x-y)||f(y)|dy
\le\|f\|_{L_a^\infty}\sup_x |x|_*^a\int|Q_j(x-y)||y|_*^{-a}dy\\
&\lesssim\|f\|_{L_a^\infty}.
\end{align*}
for any $a\in\bbZ$. This ensures that $\textsc{S}$ maps $\mcC_a^\alpha(\bbR^d)$ to $\mcC_a^\infty(\bbR^d)$ for any $\alpha\in\bbR$.

\medskip

Recall the two-parameter extension of the paraproduct, used in \cite{BH}. For any distribution $\Lambda$ on $\bbR^d\times \bbR^d$, we define
\begin{align*}
\big(\bfQ_j \Lambda\big)(x) &:= \iint_{\bbR^d\times\bbR^d} P_j(x-y)Q_j(x-z) \Lambda(y,z)dydz,   \\
\big(\bfP \Lambda\big)(x) &:= \sum_{j\ge1} \big(\bfQ_j\Lambda\big)(x).
\end{align*}
If $\Lambda(y,z)$ is of the form $f(y)g(z)$, then $\bfP \Lambda=\sfP_fg$.

\medskip

\begin{prop}\cite[Proposition~8 (a)]{BH}
\label{prop:regularity of bfP}
Fix $a\in\bbZ$.
\begin{enumerate}
	\item For any $\Lambda\in\mcS'\big(\bbR^d\times\bbR^d\big)$ for which there exists $\alpha\in\bbR$ such that $\big\|\bfQ_j\Lambda\big\|_{L_a^\infty(\bbR^d)}\lesssim2^{-j\alpha}$, for all $j\ge1$, one has $\bfP\Lambda\in \mcC_a^\alpha(\bbR^d)$ and
$$
\|\bfP\Lambda\|_{\mcC_a^\alpha(\bbR^d)} \lesssim\sup_{j\ge1}2^{j\alpha} \big\|\bfQ_j\Lambda\big\|_{L_a^\infty(\bbR^d)}.
$$
	\item For any $\alpha>0$ and $F\in\mcC_{(2),a}^\alpha(\bbR^d\times\bbR^d)$, one has $\bfP F\in \mcC_a^\alpha(\bbR^d)$ and
$$
\|\bfP F\|_{\mcC_a^\alpha(\bbR^d)}\lesssim\trino{F}_{\mcC_{(2),a}^\alpha(\bbR^d\times\bbR^d)}.
$$
\end{enumerate}
\end{prop}

\medskip

\begin{Dem}
\textit{\textsf{(a)}} Since $\mathscr{F} P_j$ is supported in the annulus $\Big\{\lambda\in\bbR^d ; |\lambda|<2^j\times\frac23\Big\}$ and $\mathscr{F} Q_j$ is supported in the annulus $\Big\{\lambda\in\bbR^d ; 2^j\times\frac34 < |\lambda| < 2^j\times\frac83\Big\}$, the integral
\begin{align*}
\int Q_i(x-w)P_j(w-y)Q_j(w-z)dw
\end{align*}
vanishes if $|i-j|\ge5$. Hence $\Delta_i({\bfP} \Lambda)=\sum_{|i-j|\le4}\Delta_i(\bfQ_j\Lambda)$ and we have
$$
\|\Delta_i({\bfP} \Lambda)\|_{L^\infty_a}
\le\sum_{|i-j|\le4}\|\Delta_i(\bfQ_j\Lambda)\|_{L^\infty_a}
\lesssim\sum_{|i-j|\le4}\|\bfQ_j\Lambda\|_{L^\infty_a}
\lesssim\sum_{|i-j|\le4}2^{-\alpha j}
\lesssim2^{-\alpha i}.
$$
For \textit{\textsf{(b)}} it is sufficient to show that $\big\|\bfQ_jF\big\|_{L_a^\infty(\bbR^d)}\lesssim2^{-j\alpha}$. By Lemma \ref{lem:weighted scaling},
\begin{align*}
\big|\bfQ_jF(x)\big| &\lesssim \int \big|P_j(x-y)Q_j(x-z)\big| \,\big(|y|_*^{-a}+|z|_*^{-a}\big) |y-z|^\alpha\, dydz   \\
&\lesssim \int \big|P_j(x-y)Q_j(x-z)\big| \big(|y|_*^{-a}+|z|_*^{-a}\big) \Big(|x-y|^\alpha+|x-z|^\alpha\Big) \,dydz   \\
&\lesssim2^{-j\alpha}|x|_*^{-a}.
\end{align*}
\end{Dem}

\medskip

Recall from \cite{BB16} the definition of the operator 
$$
{\sfR}^\circ(f,g,h) := {\sf P}_f{\sf P}_gh - {\sf P}_{fg}h.
$$
This operator is continuous from $\mcC^\alpha(\bbR^d)\times \mcC^\beta(\bbR^d)\times \mcC^\gamma(\bbR^d)$ into $\mcC^{\alpha+\beta+\gamma}(\bbR^d)$, for any $\alpha,\beta\in[0,1]$ and $\gamma\in\bbR$ -- see Proposition 14 therein.

\medskip

\begin{prop}\cite[Proposition~10]{BH}
\label{prop:sum of R}
Consider a function $f\in L_\slow^\infty(\bbR^d)$ and a finite family $(a_k,b_k)_{1\le k\le N}$ in $L_\slow^\infty(\bbR^d)\times L_\slow^\infty(\bbR^d)$ such that
$$
f(y)-f(x) = \sum_{k=1}^N a_k(x)\big(b_k(y)-b_k(x)\big) + f_{yx}^\sharp,\quad x,y\in\bbR^d,
$$
with a remainder $f_{yx}^\sharp$. 
Let $\alpha>0$ and $\beta\in\bbR$ be given. Assume that either of the following assumptions holds.
\begin{enumerate}
\item $f\in L_\rap^\infty(\bbR^d)$, $a_kb_k\in L_\rap^\infty(\bbR^d)$, $f^\sharp\in\mcC_{(2),\rap}^\alpha(\bbR^d\times\bbR^d)$, and $g\in \mcC_\slow^\beta(\bbR^d)$.
\item $f^\sharp\in\mcC_{(2)}^\alpha(\bbR^d\times\bbR^d)$ and $g\in \mcC_\rap^\beta(\bbR^d)$.
\end{enumerate}
Then one has the estimate
$$
\sum_{k=1}^N {\sfR}^\circ\big(a_k,b_k,g\big)\in \mcC_\rap^{\alpha+\beta}(\bbR^d).
$$
\end{prop}

\medskip

\begin{Dem}
Recall from identity \eqref{EqDefnSmoothS} the definition of the operator $\textsc{S}$. As in the proof of Proposition 10 in \cite{BH}, we see that
\begin{align*}
\sum_k \sfR^\circ(a_k,b_k,g) = -\textsc{S}(\sfP_fg) + \sfP_f(\textsc{S} g) - \sum_k\sfP_{a_kb_k}(\textsc{S} g) - \bfP_{x,y}\Big(\big(\sfP_{f_{\cdot x}^\sharp}g\big)(y)\Big).
\end{align*}
The first three terms belong to $\mcC_\rap^\infty(\bbR^d)$, assuming either \textit{\textsf{(a)}} or \textit{\textsf{(b)}}. Consider the last term. Note that
\begin{align*}
\bfQ_j\Big(\big(\sfP_{f_{\cdot x}^\sharp}g\big)(y)\Big)(z) = \sum_{|i-j|\le4}\int P_j(z-x)Q_j(z-y) \big(S_if_{\cdot x}^\sharp\big)(y) \,(\Delta_ig)(y)\,dxdy.
\end{align*}
For case \textit{\textsf{(a)}}, there exists $b\in\bbN$ such that $|\Delta_ig(y)|\lesssim2^{-i\beta}|y|_*^b$. Since $f^\sharp\in\mcC_{a+b}^\alpha(\bbR^d\times\bbR^d)$ for any $a\in\bbN$, one has
\begin{align*}
|(S_if_{\cdot,x}^\sharp)(y)|
&\le\int \big|P_i(y-u)\big| \big|f_{ux}^\sharp\big|\,du
\lesssim \int \big|P_i(y-u)\big| |u-x|^\alpha \big(|u|_*^{-a-b}+|x|_*^{-a-b}\big)\,du   \\
&\lesssim \int \big|P_i(y-u)\big| \big(|u-y|^\alpha+|y-x|^\alpha\big)\big(|u|_*^{-a-b}+|x|_*^{-a-b}\big)\,du   \\
&\lesssim \big(|x|_*^{-a-b}+|y|_*^{-a-b}\big)\big(2^{-i\alpha}+|y-x|^\alpha\big)
\end{align*}
by Lemma \ref{lem:weighted scaling}. Hence we have
\begin{align*}
&\Big|\bfQ_j\Big(\big(\sfP_{f_{\cdot x}^\sharp}g\big)(y)\Big)(z)\Big|   \\
&\lesssim \sum_{|i-j|\le4}\int \big|P_j(z-x)\big|\big|Q_j(z-y)\big| \big| \big(S_if_{\cdot x}^\sharp)(y)\big| \big|(\Delta_ig)(y)\big|\,dxdy   \\
&\lesssim \sum_{|i-j|\le4}\int \big|P_j(z-x)\big| \big|Q_j(z-y)\big| \big(|x|_*^{-a-b} + |y|_*^{-a-b}\big)\,|y|_*^b \,\big(2^{-i\alpha}+|y-x|^\alpha\big)\,2^{-i\beta}\,dxdy   \\
&\lesssim\sum_{|i-j|\le4}\int \big|P_j(z-x)||Q_j(z-y)\big| \big(|x|_*^{-a-b}|y|_*^b + |y|_*^{-a}\big)\Big(2^{-i\alpha}+|z-x|^\alpha+|z-y|^\alpha\Big)\,2^{-i\beta}\,dxdy   \\
&\lesssim\sum_{|i-j|\le4}|z|_*^{-a} \,\big(2^{-i\alpha} + 2^{-j\alpha}\big)\,2^{-i\beta} \lesssim |z|_*^{-a}\,2^{-j(\alpha+\beta)}.
\end{align*}
For case \textit{\textsf{(b)}}, since $\big|\Delta_ig(y)\big| \lesssim 2^{-i\beta}\,|y|_*^{-a}$ for any $a\in\bbN$, and
\begin{align*}
|(S_if_{\cdot,x}^\sharp)(y)|\le
\int \big|P_i(y-u)\big|\,\big|f_{ux}^\sharp\big|\,du &\lesssim\int \big|P_i(y-u)\big| |u-x|^\alpha\,du \lesssim 2^{-i\alpha}+|y-x|^\alpha,
\end{align*}
we have
\begin{align*}
&\Big|\bfQ_j\Big(\big(\sfP_{f_{\cdot x}^\sharp}g\big)(y)\Big)(z)\Big| 
\lesssim \sum_{|i-j|\le4} \int \big|P_j(z-x)\big|\, \big|Q_j(z-y)\big| \,\big| (S_if_{\cdot x}^\sharp)(y)\big| \,\big|(\Delta_ig)(y)\big|\,dxdy   \\
&\quad\lesssim \sum_{|i-j|\le4} \int \big|P_j(z-x)\big| \,\big|Q_j(z-y)\big| \,|y|_*^{-a} \,\big(2^{-i\alpha} + |y-x|^\alpha\big)\,2^{-i\beta}\,dxdy\\
&\quad\lesssim \sum_{|i-j|\le4} \int \big|P_j(z-x)\big|\,\big|Q_j(z-y)\big|\, |y|_*^{-a} \,\Big(2^{-i\alpha} + |z-x|^\alpha+|z-y|^\alpha\Big)\,2^{-i\beta}\,dxdy   \\
&\quad\lesssim \sum_{|i-j|\le4} |z|_*^{-a} \,\big(2^{-i\alpha}+2^{-j\alpha}\big) \,2^{-i\beta} \lesssim |z|_*^{-a}\,2^{-j(\alpha+\beta)}.
\end{align*}
By Proposition \ref{prop:regularity of bfP}, we are done.
\end{Dem}

\medskip

\begin{prop} \cite[Proposition~9]{BH} \label{prop:reconstruction}
Let $\gamma\in\bbR$ and $\beta_0\in\bbR$ be given together with a family $\Lambda_x$ of distributions on $\bbR^d$, indexed by $x\in\bbR^d$. Assume one has
$$
\sup_{x\in\bbR^d}|x|_*^a\|\Lambda_x\|_{\mcC^{\beta_0}}< \infty
$$
for any $a\in\bbZ$ and one can decompose $(\Lambda_y - \Lambda_x)$ under the form
   \begin{equation*}
   \Lambda_y - \Lambda_x = \sum_{\ell=1}^L c^\ell_{yx}\,\Theta_x^\ell   
   \end{equation*}
for $L$ finite, $\bbR^d$-indexed distributions $\Theta_x^\ell$, and real-valued coefficients $c_{yx}^\ell$ depending measurably on $x$ and $y$.
Assume that for each $\ell$ there exists $\beta_\ell<\gamma$ such that either of the following conditions holds.
\begin{enumerate}
\item $\Theta^\ell\in D_{\rap}^{\beta_\ell}$ and $c^\ell\in\mcC_{(2)}^{\gamma-\beta_\ell}(\bbR^d\times\bbR^d)$.
\item $\Theta^\ell\in D^{\beta_\ell}$ and $c^\ell\in\mcC_{(2),\rap}^{\gamma-\beta_\ell}(\bbR^d\times\bbR^d)$.
\end{enumerate}

\noindent Write ${\bfP}(\Lambda)$ for ${\bfP}_{y,z}\big(\Lambda_y(z)\big)$ below.  \vspace{0.15cm}

\begin{enumerate}
\renewcommand{\labelenumi}{\textbf{\textsf{(\roman{enumi})}}}
\item If $\gamma>0$, then there exists a unique function $\lambda\in \mcC_{\rap}^\gamma(\bbR^d)$ such that 
$$
\Big\{\big(\bfP(\Lambda)-\lambda\big) - \Lambda_x\Big\}_{x\in\bbR^d}\in D_\rap^\gamma.
$$
\item If $\gamma<0$, then
$$
\Big\{\bfP(\Lambda)-\Lambda_x\Big\}_{x\in\bbR^d}\in D_\rap^\gamma.
$$
\end{enumerate}
Consequently, $\bfP(\Lambda)\in \mcC_\rap^{\beta_0}(\bbR^d)$. If furthermore $\Lambda\in D_\rap^\gamma$, then $\bfP(\Lambda)\in \mcC_\rap^\gamma(\bbR^d)$.
\end{prop}

\medskip

\begin{Dem}
In view of Proposition 9 in \cite{BH}, it is sufficient to show that
\begin{align}\label{eq:key of reconst}
\sup_{x\in\bbR^d}|x|_*^a\,\Big|\Delta_j\big(\bfP(\Lambda)-\Lambda_x\big)(x)\Big| \lesssim 2^{-j\gamma}.
\end{align}
We write for that purpose
\begin{align*}
&\Delta_i\big(\bfP(\Lambda)-\Lambda_x\big)(x)\\
&=\sum_{j\ge1,\, |i-j|\le4}\iiint Q_i(x-y)P_j(y-u)Q_j(y-v)\,(\Lambda_u-\Lambda_x)(v)\,dydudv - \Delta_i\textsc{S}(\Lambda_x)(x)\\
&=:A+B+C,
\end{align*}
where
\begin{align*}
A
&=\sum_{|i-j|\le4}\iiint Q_i(x-y)P_j(y-u)Q_j(y-v)\,(\Lambda_u-\Lambda_y)(v)\,dydudv\\
&=\sum_{|i-j|\le4}\sum_{\ell=1}^L \iiint Q_i(x-y)P_j(y-u)Q_j(y-v)\,c_{uy}^\ell\Theta_y(v)\,dydudv
\end{align*}
and
\begin{align*}
B
&=\sum_{|i-j|\le4}\iiint Q_i(x-y)P_j(y-u)Q_j(y-v)\,(\Lambda_y-\Lambda_x)(v)\,dydudv\\
&=\sum_{|i-j|\le4}\sum_{\ell=1}^L \iiint Q_i(x-y)P_j(y-u)Q_j(y-v)\,c_{yx}^\ell\Theta_x(v)\,dydudv
\end{align*}
For the C term,
\begin{align*}
\sup_x|x|_*^a\,\big|\Delta_i\textsc{S}(\Lambda_x)(x)\big| &\lesssim2^{-ir}\sup_x|x|_*^a\,\big\|\textsc{S}(\Lambda_x)\big\|_{\mcC^r}   \\
&\lesssim 2^{-ir}\sup_x|x|_*^a\,\|\Lambda_x\|_{\mcC^{\beta_0}}\lesssim2^{-ir}
\end{align*}
for any $r>0$. For the $A$ term, for any $a\in\bbZ$ we have
\begin{align*}
&|A|
\le\sum_{|i-j|\le4}\sum_{\ell=1}^L \iint \big|Q_i(x-y)\big| \big|P_j(y-u)\big| |c_{uy}^\ell| |\Delta_j\Theta_y(y)|\,dydu\\
&\lesssim
\left\{
\begin{aligned}
&\sum_{|i-j|\le4}\sum_{\ell=1}^L \iint \big|Q_i(x-y)\big| \big|P_j(y-u)\big|
|u-y|^{\gamma-\beta_\ell}|y|_*^{-a}2^{-j\beta_\ell} dydu,
&&\text{if {\textit{\textsf{(a)}}}}\\
&\sum_{|i-j|\le4}\sum_{\ell=1}^L \iint \big|Q_i(x-y)\big| \big|P_j(y-u)\big|
\big(|u|_*^{-a}+|y|_*^{-a}\big)|u-y|^{\gamma-\beta_\ell}2^{-j\beta_\ell} dydu,
&&\text{if {\textit{\textsf{(b)}}}}
\end{aligned}
\right.\\
&\lesssim\sum_{|i-j|\le4}
\int\big|Q_i(x-y)\big||y|_*^{-a}2^{-j\gamma}dy
\lesssim |x|_*^{-a}2^{-i\gamma}.
\end{align*}
The $B$ term has the same estimate by a similar argument. So estimate \eqref{eq:key of reconst} follows from Lemma \ref{lem:weighted scaling}.

\ssk

\textbf{\textsf{(i)}}
If $\gamma>0$, the estimate \eqref{eq:key of reconst} implies that the sum
$$
\lambda(x):=\sum_{j\ge-1}\Delta_j\big(\bfP(\Lambda)-\Lambda_x\big)(x)
$$
defines an element $\lambda$ of $\mcC_{\rap}^\gamma(\bbR^d)$.
To show it, we follow the argument in Section 6 of \cite{GIP}.
We decompose $\lambda=\lambda^{\le j+1}+\lambda^{>j+1}$, where
$$
\lambda^{\le j+1}(x):=\sum_{i\le j+1}\Delta_i\big(\bfP(\Lambda)-\Lambda_x\big)(x)
=S_{j+3}\big(\bfP(\Lambda)-\Lambda_x\big)(x).
$$
We consider $\Delta_j\lambda=\Delta_j\lambda^{\le j+1}+\Delta_j\lambda^{>j+1}$.
For the second term, by the estimate \eqref{eq:key of reconst} one has
$$
\|\Delta_j\lambda^{>j+1}\|_{L_a^\infty}
\lesssim\|\lambda^{>j+1}\|_{L_a^\infty}
\lesssim\sum_{i>j+1}2^{-i\gamma}\lesssim2^{-j\gamma}.
$$
For the first term, since $\Delta_jS_{j+3}=\Delta_j$, one has
\begin{align*}
\Delta_j\lambda^{\le j+1}(y)
&=\int Q_j(y-x)S_{j+3}\big(\bfP(\Lambda)-\Lambda_x\big)(x)dx\\
&=\int Q_j(y-x)S_{j+3}\Big(\bfP(\Lambda)-\Lambda_y+\sum_\ell c_{yx}^\ell\Theta_x^\ell\Big)(x)dx\\
&=\Delta_j\big(\bfP(\Lambda)-\Lambda_y\big)(y)+\sum_\ell\int Q_j(y-x)c_{yx}^\ell (S_{j+3}\Theta_x)(x)dx
\end{align*}
Similarly to above, we can show that $|\Delta_j\lambda^{\le j+1}(y)|\lesssim|y|_*^{-a}2^{-j\gamma}$ for any $a\in\bbZ$.
In the end we have $\|\Delta_j\lambda\|_{L_a^\infty}\lesssim2^{-j\gamma}$, hence $\lambda\in\mcC_\rap^\gamma(\bbR^d)$.

\ssk

Since $\sum_{i\ge-1}\Delta_i\big(\bfP(\Lambda)-\Lambda_x-\lambda\big)(x)=0$ by definition, we have
\begin{align*}
\big| S_i\big(\bfP(\Lambda)-\Lambda_x-\lambda\big)(x) \big|
\le\sum_{j\ge i-1}\big| \Delta_j\big(\bfP(\Lambda)-\Lambda_x-\lambda\big)(x) \big|
\lesssim|x|_*^{-a}\sum_{j\ge i-1}2^{-j\gamma}\lesssim|x|_*^{-a}2^{-i\gamma}
\end{align*}
for any $a\in\bbZ$.

\ssk

\textbf{\textsf{(ii)}}
If $\gamma<0$, then directly from \eqref{eq:key of reconst},
\begin{align*}
\big| S_i\big(\bfP(\Lambda)-\Lambda_x\big)(x) \big|
\le\sum_{j< i-1}\big| \Delta_j\big(\bfP(\Lambda)-\Lambda_x\big)(x) \big|
\lesssim|x|_*^{-a}\sum_{j< i-1}2^{-j\gamma}\lesssim|x|_*^{-a}2^{-i\gamma}
\end{align*}
for any $a\in\bbZ$.
\end{Dem}

\medskip

\begin{cor}\label{canonical reconstruction}
Given a concrete regularity structure $\scrT$ satisfying assumptions \textbf{\textsf{(A)}} and \textbf{\textsf{(B)}} and given a rapidly decreasing model $\sfM=(\sfg,\sfPi)$, we define the map $\textsf{\textbf{R}} : \mcD_\rap^\gamma(T,\sfg)\to \mcC_\rap^{\beta_0}$, by
$$
\textsf{\textbf{R}}\bsf = \bfP_{x,y}\Big(\big(\sfPi_x^\sfg\bsf(x)\big)(y)\Big).
$$
Then one has
$$
\Big(\textsf{\textbf{R}}\bsf-\sfPi_x^\sfg\bsf(x)\Big)_{x\in\bbR^d} \in D_\rap^\gamma.
$$
\end{cor}

\medskip

\begin{Dem}
Let $\Lambda_x=\sfPi_x^\sfg\bsf(x)$. Since
$$
\Lambda_y-\Lambda_x = \sum_{\tau\in\mcB} \big\langle\tau',\widehat{\sfg_{xy}}\bsf(y)-\bsf(x)\big\rangle\,\sfPi_x^\sfg\tau
$$
one gets conditions {\textit{\textsf{(a)}}} and {\textit{\textsf{(b)}}} of Proposition \ref{prop:reconstruction} from the definition of a model.
\end{Dem}

\medskip

\begin{Dem}[\textbf{\textsf{of Theorem \ref{PropDefnBracket}}}]
We prove the case $m=0$ here for simplicity. For general $m$, the proof is at the end of this appendix, after we introduce the modified paraproducts.

\smallskip

Consider the first formula \eqref{def:bracket g}. First we show that, for each $\tau\in\mcB^+$ we have
\begin{align}\label{predef:bracket g}
\sfg(\tau)=\sum_{\unit<^+\nu<^+\tau,\, \nu\in\mcB^+}\sfP_{\sfg(\tau/^+\nu)}[\nu]^\sfg+[\tau]^\sfg,
\end{align}
where $[\nu]^\sfg\in \mcC_\rap^{|\nu|}(\bbR^d)$, if $\nu\in\mcB^+\setminus\mcB_X^+$, and $[\nu]^\sfg\in \mcC_\slow^\infty(\bbR^d)$, if $\nu\in\mcB_X^+$. If $\tau=X^k$, then since $\Delta^+X^k=\sum_{\ell}\binom{k}{\ell}X^\ell\otimes X^{k-\ell}$ we have
$$
\sfg(X^k)=\sum_{0<\ell<k}\binom{k}{\ell}\sfP_{\sfg(X^\ell)}[X^{k-\ell}]^\sfg+[X^k]^\sfg.
$$
We see for instance that $[\textbf{\textsf{1}}]^{\sf g}=1$, then $[X]^{\sf g}=x$, since ${\sf g}_x(X)=x$, and since ${\sf g}_x(X^2)=x^2$, one has 
$$
x^2 = 2{\sf P}_xx + [X^2]^{\sf g}.
$$ 
We recognize $[X^2]^{\sf g}={\sf \Pi}(x,x)$. More generally, since $\sfg_x(X^k)=x^k$ is a function belonging to $C_\slow^\infty(\bbR^d)$, by an induction we have $[X^k]^\sfg\in C_\slow^\infty(\bbR^d)$. Now let $\tau\in\mcB^+\setminus\mcB_X^+$. Recall the formula obtained in \cite{BH};
\begin{align*}
&[\tau]^\sfg = \textsc{S}\sfg(\tau)+\bfP_{x,y}\big(\sfg_{yx}(\tau)\big)   \\
&+ \sum_{n=1}^\infty(-1)^{n-1}\sum_{\unit<^+\sigma_{n+1}<^+\cdots<^+\sigma_1<^+\tau} {\sfR}^\circ\Big(\sfg(\tau/^+\sigma_1)\cdots\sfg(\sigma_{n-1}/^+\sigma_n),\sfg(\sigma_n/^+\sigma_{n+1}),[\sigma_{n+1}]^\sfg\Big).
\end{align*}
This is obtained from the expansion formula obtained in \cite{BH};
\begin{align}\label{expansion of g(tau/sigma)}
\begin{aligned}
&\sfg_y(\tau/^+\sigma)-\sfg_x(\tau/^+\sigma)\\
&=\sum_{n=1}^\infty(-1)^{n-1}\sum_{\sigma<^+\sigma_n<^+\cdots<^+\sigma_1<^+\tau}
\sfg_x(\tau/^+\sigma_1)\cdots\sfg_x(\sigma_{n-1}/^+\sigma_n)\Big(\sfg_y(\sigma_n/^+\sigma)-\sfg_x(\sigma_n/^+\sigma)\Big)\\
&\quad+\sfg_{yx}(\tau/^+\sigma)
\end{aligned}
\end{align}
with $\sigma={\bf1}_+$ and by definition of the ${\sf R}^\circ$ operator. Since $\tau\in\mcB^+\setminus\mcB_X^+$, we have $\textsc{S}\sfg(\tau)\in \mcC_\rap^\infty(\bbR^d)$ and $\bfP_{x,y}(\sfg_{yx}(\tau))\in \mcC_\rap^{|\tau|}(\bbR^d)$. For the ${\sfR}^\circ$ terms, we apply Proposition \ref{prop:sum of R} to \eqref{expansion of g(tau/sigma)}. If $\sigma\in\mcB_X^+$, then since $\tau/^+\sigma\in\textrm{span}(\mcB^+\setminus\mcB_X^+)$, by assumption \textbf{\textsf{(B-2)}}, we have $\sfg_x(\tau/^+\sigma)\in L_{\rap}^\infty(\bbR^d)$ and $\sfg_{yx}(\tau/^+\sigma)\in\mcC_{(2),\rap}^{|\tau|-|\sigma|}(\bbR^d\times\bbR^d)$. For the sum over $\sigma<^+\sigma_n<^+\cdots<^+\sigma_1<^+\tau$, we can see that at least one element among
$$
\sfg(\tau/^+\sigma_1),\quad \dots,\quad\sfg(\sigma_{n-1}/^+\sigma_n),\quad\sfg(\sigma_n/^+\sigma)
$$
belongs to $L_\rap^\infty(\bbR^d)$. Indeed, if $\sigma_n\notin\mcB_X^+$ then $\sfg(\sigma_n/^+\sigma)\in L_\rap^\infty(\bbR^d)$. Otherwise, if $\sigma_{n-1}\notin\mcB_X^+$ then $\sfg(\sigma_{n-1}/^+\sigma_n)\in L_\rap^\infty(\bbR^d)$. Since $\tau\notin\mcB_X^+$, for at least one $i$ we have $\sfg(\sigma_i/^+\sigma_{i+1})\in L_\rap^\infty(\bbR^d)$. Since $L_\slow^\infty(\bbR^d)\cdot L_\rap^\infty(\bbR^d)\subset L_\rap^\infty(\bbR^d)$, we can apply Proposition \ref{prop:sum of R}-\textit{\textsf{(a)}} to get
$$
\sum_{n=1}^\infty(-1)^{n-1}\sum_{\sigma<^+\sigma_n<^+\cdots<^+\sigma_1<^+\tau} \sfR^\circ\Big(\sfg(\tau/^+\sigma_1)\cdots\sfg(\sigma_{n-1}/^+\sigma_n),\sfg(\sigma_n/^+\sigma),[\sigma]^\sfg\Big)\in \mcC_\rap^{|\tau|}(\bbR^d).
$$
If $\sigma\notin\mcB_X^+$, since $\sfg_{yx}(\tau/^+\sigma)\in\mcC_{(2)}^{|\tau|-|\sigma|}(\bbR^d\times\bbR^d)$ and $[\sigma]^\sfg\in \mcC_\rap^{|\sigma|}(\bbR^d)$ we can apply Proposition \ref{prop:sum of R}-\textit{\textsf{(b)}} to get the same estimate. Hence we obtain the required estimates in the formula \eqref{predef:bracket g}.

\ssk

To get \eqref{def:bracket g} from \eqref{predef:bracket g}, it is sufficient to show
\begin{align}\label{diff in C_rap^infty}
\Brac{\tau}^\sfg-[\tau]^\sfg\in \mcC_\rap^\infty(\bbR^d)
\end{align}
for any $\tau\in\mcB^+\setminus\mcB_X^+$. Assume that all $\nu\in\mcB^+\setminus\mcB_X^+$ with $|\nu|<|\tau|$ satisfy \eqref{diff in C_rap^infty}. Then we have
\begin{align*}
\Brac{\tau}^\sfg-[\tau]^\sfg &= \sum_{\unit<^+\nu<^+\tau}\sfP_{\sfg(\tau/^+\nu)}[\nu]^\sfg-\sum_{\unit<^+\nu<^+\tau,\,\nu\notin\mcB_X^+}\sfP_{\sfg(\tau/^+\nu)}\Brac{\nu}^\sfg   \\
&= \sum_{\unit<^+\nu<^+\tau,\,\nu\notin\mcB_X^+}\sfP_{\sfg(\tau/^+\nu)}\big([\nu]^\sfg-\Brac{\nu}^\sfg\big) + \sum_{k\neq0}\sfP_{\sfg(\tau/^+X^k)}[X^k]^\sfg.
\end{align*}
The first term belongs to $\mcC_\rap^\infty(\bbR^d)$ by assumption. For the second term, since $[X^k]^\sfg\in \mcC_\slow^\infty(\bbR^d)$ and $\sfg(\tau/^+X^k)\in L_\rap^\infty(\bbR^d)$, we can complete the proof.

\ssk

One can obtain formula \eqref{def:bracket Pi} in the similar way. As above, we define the quantity $[\tau]^{\sf M}$ for each $\tau\in\mcB$ by
\begin{align*}
{\sf\Pi}\tau=\sum_{\nu<\tau,\, \nu\in\mcB}\sfP_{\sfg(\tau/\nu)}[\nu]^{\sf M}+[\tau]^{\sf M}.
\end{align*}
Then we can show that $[\nu]^{\sf M}\in \mcC_\rap^{|\nu|}(\bbR^d)$, if $\nu\in\mcB\setminus\mcB_{\underline{X}}$, and $[\nu]^{\sf M}\in \mcC_\slow^\infty(\bbR^d)$, if $\nu\in\mcB_{\underline{X}}$.
The only difference is that, for $\tau\in\mcB\setminus\mcB_{\underline{X}}$, we use the formula obtained in \cite{BH};
\begin{align*}
&[\tau]^\sfg = \textsc{S}({\sf\Pi}\tau)+\bfP_{x,y}\big(({\sf\Pi}_x^{\sfg}\tau)(y)\big)   \\
&+ \sum_{n=1}^\infty(-1)^{n-1}\sum_{\sigma_{n+1}<\cdots<\sigma_1<\tau} {\sfR}^\circ\Big(\sfg(\tau/\sigma_1)\cdots\sfg(\sigma_{n-1}/\sigma_n),\sfg(\sigma_n/\sigma_{n+1}),[\sigma_{n+1}]^{\sf M}\Big).
\end{align*}
and use Proposition \ref{prop:reconstruction} to get $\bfP_{x,y}\big(({\sf\Pi}_x^{\sf g}\tau)(y)\big)\in \mcC_\rap^{|\tau|}(\bbR^d)$.
Since the property \eqref{diff in C_rap^infty} also holds for the operator $\Brac{\cdot}^{\sf M}-[\cdot]^{\sf M}$, we can conclude \eqref{def:bracket Pi}.
\end{Dem}

\bigskip

\begin{Dem}[\textbf{\textsf{of Theorem \ref{thm:MD to PD}}}]
\eqref{eq:para system} is proved by a similar argument as Theorem \ref{PropDefnBracket}.
See Theorem 14 in \cite{BH} for details.
More easily, it is useful to consider the extended algebra $T_{\bsF}^+$ defined in Section \ref{SubsectionPCtoModelledDistributions}.
Since a modelled distribution $\bsf\in\mcD^\gamma(T,\sfg)$ defines a $\sfg$-part of the model on $T_{\bsF}^+$ by Lemma \ref{lem: model on Ftau is equiv to MD}, we have
$$
f_\sigma=\sfg(\bsF_\sigma)=\sum_{\substack{\sigma<\mu \\ \mu/\sigma\in\text{\rm span}(\mcB^+\setminus\mcB_X^+)}}{\sf P}_{\sfg(\bsF_\mu)}\Brac{\mu/\tau}^{\sfg}+\Brac{\bsF_\sigma}^{\sfg}.
$$
Thus $\Brac{f_\sigma}^{\sfg}=\Brac{\bsF_\sigma}^{\sfg}\in\mcC_{\rap}^{\gamma-|\sigma|}(\bbR^d)$.

\ssk

As for \eqref{eq:para reconst}, a similar interpretation is useful. Consider a symbol $\bsF$ and an extended model space $T_{\bsF}:=T\oplus\text{\rm span}(\bsF)$. Giving the homogeneity $|\bsF|:=\gamma$ and the coproduct formula
$$
\Delta\bsF=\bsF\otimes{\bf1}_++\sum_{\tau\in\mcB,\,|\tau|<\gamma}\tau\otimes(\bsF_\tau),
$$
the pair $(T_{\bsF}^+,T_{\bsF})$ turns out to be a regularity structure. (It is not difficult to check that $T_{\bsF}$ is a comodule over $T_{\bsF}^+$ by using \eqref{eq:def of Delta F/tau}.)
For given a reconstruction $\reR\bsf$, we can define the model on $T_{\bsF}$ by setting ${\sf\Pi}\bsF:=\reR\bsf$. Indeed, similarly to Lemma \ref{lem: model on Ftau is equiv to MD}, we can show that
$$
{\sf\Pi}_x^{\sfg}\bsF=\reR\bsf-{\sf\Pi}_x^{\sfg}\bsf(x).
$$
Then \eqref{eq:para reconst} follows from \eqref{def:bracket Pi} in Theorem \ref{PropDefnBracket}.
\end{Dem}

\bigskip

To complete the proof of Theorem \ref{PropDefnBracket}, we define here the two-parameter extension $\bfP^m$ of the modified paraproduct $\sfP^m$. Note that, there is an annulus $\mcA\subset\bbR^d$ such that the Fourier transform of the function
$$
x\mapsto P_j(x-y)\,Q_j(x-z)
$$
is contained in $2^j\mcA$ (independently to $y,z$). Let $\chi$ be a smooth function on $\bbR^d$ supported in a larger annulus $\mcA'$ and such that $\chi\equiv1$ on $\mcA$. Letting $R_j=\mcF^{-1}\big(\chi(2^{-j}\cdot)\big)$, we have
$$
(\bfQ_j \Lambda)(x)=\iiint_{\bbR^d\times\bbR^d\times\bbR^d} R_j(x-w)P_j(w-y)Q_j(w-z) \Lambda(y,z)\,dydzdw.
$$
For $m\in\bbZ$, set
\begin{equation*} \begin{split}
Q_j^{-m} &:= \mcF^{-1}\big(|\cdot|^{-m}\rho_j\big),   \\
R_j^m &:= \mcF^{-1}\big(|\cdot|^m\chi(2^{-j}\cdot)\big);
\end{split} \end{equation*}
then they are smooth functions such that $Q_j^{-m} = |\nabla|^{-m}Q_j$ and $R_j^m = |\nabla|^m R_j$.

\medskip

\begin{defn}
For any $m\in\bbN$ and any two-variable distribution $\Lambda$ on $\bbR^d\times\bbR^d$, define
\begin{align*}
(\bfQ_j^m \Lambda)(x) &:= \iiint_{\bbR^d\times\bbR^d\times\bbR^d} R_j^m(x-w)P_j(w-y)Q_j^{-m}(w-z) \Lambda(y,z)\,dydzdw,   \\
(\bfP^m \Lambda)(x) &:= \sum_{j\ge1}(\bfQ_j^m\Lambda)(x).
\end{align*}
\end{defn}

\medskip

If necessary, we emphasize the integrated variables by writing
$$
\bfP^m\Lambda =\bfP_{y,z}^m\big(\Lambda(y,z)\big).
$$
For the special case $\Lambda(y,z)=f(y)g(z)$, we have the consistency relation
\begin{align*}
\bfP^m\Lambda = {\sfP}^m_fg.
\end{align*}
All the above estimates in this appendix still hold for these modified operators. Indeed, because of the scaling properties
$$
Q_j^{-m}(x)=2^{j(d-m)}Q_0^{-m}(2^jx),\quad
R_j^m(x)=2^{j(d+m)}R_0^m(2^jx),
$$
we can show the following analogue of Lemma \ref{lem:weighted scaling}; for any $\alpha\ge0$ and $a\in\bbZ$, one has
\begin{align}\label{lem:weighted scaling 2}
\begin{aligned}
\int \big|R_i^m(x-y)\big| |x-y|^\alpha|y|_*^{-a}dy&\lesssim2^{-i(\alpha+m)}\,|x|_*^{-a},   \\
\int \big|Q_i^{-m}(x-y)\big| |x-y|^\alpha|y|_*^{-a}dy&\lesssim2^{-i(\alpha-m)}\,|x|_*^{-a}.
\end{aligned}
\end{align}
Thus we can repeat the argument in this appendix as follows.
\begin{itemize}
\item Proposition \ref{prop:regularity of bfP}-\textsf{\textit{(a)}} still holds, since $\Delta_i(\bfP^m\Lambda)=\sum_{|i-j|\le N}\Delta_i(\bfQ_j^m\Lambda)$ for some integer depending only on the support of $\chi$.
\item Proposition \ref{prop:regularity of bfP}-\textsf{\textit{(b)}} still holds, since by the scaling property,
\begin{align*}
\big|\bfQ_j^mF(x)\big| &\lesssim \int \big|R_j^m(x-w)P_j(w-y)Q_j^{-m}(w-z)\big| \,|F(y,z)|\, dydzdw   \\
&\lesssim \int \big|R_j^m(x-w)P_j(w-y)Q_j^{-m}(w-z)\big| \big(|y|_*^{-a}+|z|_*^{-a}\big) \Big(|w-y|^\alpha+|w-z|^\alpha\Big) \,dydzdw   \\
&\lesssim 2^{-j(\alpha-m)}\int |R_j^m(x-w)| |w|_*^{-a}dw\\
&\lesssim 2^{-j\alpha}|w|_*^{-a}.
\end{align*}
\item Proposition \ref{prop:sum of R} still holds if ${\sf R}^{\circ}$ is replaced by
$$
{\sf R}^m(f,g,h):={\sf P}_f^m{\sf P}^m_gh-{\sf P}^m_{fg}h,
$$
by a parallel argument using \eqref{lem:weighted scaling 2}.
\end{itemize}
Consequently, we can repeat the proof of Theorem \ref{PropDefnBracket} for any $m\in\bbN$.

\bigskip

\section{The slowly growing setting}  
\label{AppendixSlow}

In applications of regularity structures to the study of singular stochastic PDEs set in the entire space $\bbR^d$ usually involve noises that do not have rapid decrease at infinity, but rather have moderate growth at infinity. Our results can be formulated as follows in this slightly modified setting.

\ssk

We define the spaces $\scrM_{\slow}$ and $\mcD_{\slow}^\gamma$ of slowly growing models and modelled distributions, respectively, by replacing `rap' in definitions in Section \ref{SubsectionFromModelsToPCSystems} by `slow'. We can repeat the same arguments to obtain the variations of Theorems \ref{mainthm1} and \ref{mainthm2} with the spaces $\scrM_{\slow}$ and $\mcD_{\slow}^\gamma$, respectively.
All we need is to consider the weight $|x|_*^{-a}$ for {\it some} sufficiently large $a$, instead of {\it any} $a\in\bbZ$.
Precisely, we need the following minor modifications of the arguments.
\begin{itemize}
\item Proposition \ref{prop:sum of R} still holds under the assumption $f^\sharp\in\mcC_{(2),\slow}^\alpha(\bbR^d\times\bbR^d)$ and $g\in\mcC_{\slow}^\beta(\bbR^d)$, instead of rapid decrease assumptions.
\item Proposition \ref{prop:reconstruction} still holds under the assumption
$$
\sup_{x\in\bbR^d}|x|_*^a\|\Lambda_x\|_{\mcC^{\beta_0}}<\infty
$$
for {\it some} $a\in\bbZ$, and for any $\ell$, $\Theta^\ell\in D_{\slow}^{\beta_\ell}$ and $c^\ell\in\mcC_{(2),\slow}^{\gamma-\beta_\ell}(\bbR^d\times\bbR^d)$.
\item Lemma \ref{lem:I_j^k} still holds for {\it some} $a\in\bbZ$, instead of any $a$.
\end{itemize}
Details are left to readers. We end this appendix by writing the precise statements of main theorems.

\medskip

\begin{thm}
Let $\scrT$ be a concrete regularity structure satisfying assumptions \textbf{\textsf{(A-C)}}. Then one can construct a locally Lipschitz continuous map
\begin{equation} \begin{split}
\scrM_\slow(\scrT, \bbR^d) &\to
\prod_{\sigma\in\mcB^+\setminus\mcB_X^+} \mcC_\slow^{|\sigma|}(\bbR^d) \times \prod_{\tau\in\mcB\setminus\mcB_{\underline{X}}} \mcC_\slow^{|\tau|}(\bbR^d)   \\
{\sf (g,\Pi)} &\mapsto \Big(\Brac{\sigma}^{\sf M}, \Brac{\tau}^{\sf g}\,;\, \sigma\in\mcB^+\setminus\mcB_X^+, \tau\in\mcB\setminus\mcB_{\underline{X}} \Big)
\end{split} \end{equation}
by giving paracontrolled representations of $\sf g$ and $\sf \Pi$, for $({\sf g}, {\sf \Pi})\in \scrM_\rap(\scrT, \bbR^d)$. Furthermore, $\scrM_\slow(\scrT, \bbR^d)$ is locally bi-Lipschitz homeomorphic to the product space
\begin{align*}
\prod_{\sigma\in\mcG_\circ^+} \mcC_\slow^{|\sigma|}(\bbR^d) \times \prod_{\tau\in\mcB_\bullet,\,|\tau|<0} \mcC_\slow^{|\tau|}(\bbR^d).
\end{align*}
\end{thm}

\begin{thm}
Let a concrete regularity structure $\scrT$ satisfy assumptions \textbf{\textsf{(A-D)}}. Pick $\gamma\in\bbR\setminus\{0\}$ such that $\gamma-|\tau|\notin\bbN$ for any basis element $\tau$ of $T$ with $|\tau|<\gamma$, and $\sfM=(\sfg,\sfPi)\in\scrM_\slow(\scrT, \bbR^d)$. Then one can construct a locally Lipschitz continuous map
$$
\mcD_\slow^\gamma(T,\sfg) \to \prod_{\tau\in\mcB,\, |\tau|<\gamma} \mcC_\slow^{\gamma-|\tau|}(\bbR^d)
$$
by giving a paracontrolled representation of elements in $\mcD_\slow^\gamma(T,\sfg)$. Furthermore, $\mcD_\slow^\gamma(T, \sfg)$ is locally bi-Lipschitz homeomorphic to the product space
$$
\prod_{\tau\in\mcB_\bullet,\, |\tau|<\gamma} \mcC_\slow^{\gamma-|\tau|}(\bbR^d).
$$
\end{thm}

\bigskip
\bigskip

\noindent \textcolor{gray}{$\bullet$} {\sf I. Bailleul} -- Univ. Rennes, CNRS, IRMAR - UMR 6625, F-35000 Rennes, France.   \\
\noindent {\it E-mail}: ismael.bailleul@univ-rennes1.fr   

\medskip

\noindent \textcolor{gray}{$\bullet$} {\sf M. Hoshino} --  Faculty of Mathematics, Kyushu University, Japan   \\
{\it E-mail}: hoshino@math.kyushu-u.ac.jp

\end{document}